\documentclass{amsart}
\title{Comparison of types for inner forms of $\mathrm{GL}_{N}$}
\author{Yuki Yamamoto}
\usepackage{mathrsfs}
\usepackage{amssymb}
\usepackage[all]{xy}
\newtheorem{defn}{Definition}[section]
\newtheorem{rem}[defn]{Remark}
\newtheorem{thm}[defn]{Theorem}
\newtheorem{prop}[defn]{Proposition}
\newtheorem{lem}[defn]{Lemma}
\newtheorem{cor}[defn]{Corollary}

\newenvironment{prf}{Proof. }{\hfill$\Box$}

\DeclareMathOperator{\Hom}{Hom}
\DeclareMathOperator{\GL}{GL}
\DeclareMathOperator{\id}{id}
\DeclareMathOperator{\End}{End}
\DeclareMathOperator{\Aut}{Aut}
\DeclareMathOperator{\M}{M}
\DeclareMathOperator{\Res}{Res}
\DeclareMathOperator{\Ind}{Ind}
\DeclareMathOperator{\cInd}{c--Ind}
\DeclareMathOperator{\ord}{ord}
\DeclareMathOperator{\Lie}{Lie}
\DeclareMathOperator{\Tr}{Tr}
\DeclareMathOperator{\Trd}{Trd}
\DeclareMathOperator{\Nrm}{N}
\DeclareMathOperator{\Nrd}{Nrd}
\DeclareMathOperator{\Gal}{Gal}
\DeclareMathOperator{\Cent}{Cent}
\DeclareMathOperator{\Stab}{Stab}
\DeclareMathOperator{\Latt}{Latt}
\DeclareMathOperator{\sr}{sr}
\newcommand{\Afra}{\mathfrak{A}}
\newcommand{\Bfra}{\mathfrak{B}}
\newcommand{\Bscr}{\mathscr{B}}
\newcommand{\gfra}{\mathfrak{g}}
\newcommand{\Kfra}{\mathfrak{K}}
\newcommand{\Lcal}{\mathcal{L}}

\newcommand{\ofra}{\mathfrak{o}}
\newcommand{\Pfra}{\mathfrak{P}}
\newcommand{\pfra}{\mathfrak{p}}
\newcommand{\Qfra}{\mathfrak{Q}}
\newcommand{\sfra}{\mathfrak{s}}
\newcommand{\Phibf}{\mathbf{\Phi}}
\newcommand{\dr}{\mathrm{d}}
\newcommand{\rbf}{\mathbf{r}}
\newcommand{\U}{\mathbf{U}}

\newcommand{\Z}{\mathbb{Z}}
\newcommand{\R}{\mathbb{R}}
\newcommand{\C}{\mathbb{C}}

\begin{document}
\maketitle

\begin{abstract}
Let $F$ be a non-archimedean local field, $A$ be a central simple $F$-algebra, and $G$ be the multiplicative group of $A$.  
To construct types for supercuspidal representations of $G$, simple types by S\'echerre--Stevens and Yu's construction are already known.  
In this paper, we compare these constructions.  
In particular, we show essentially tame supercuspidal representations of $G$ defined by Bushnell--Henniart are nothing but tame supercuspidal representations defined by Yu.  
\tableofcontents
\end{abstract}

\section{Introduction}
Let $F$ be a non-archimedean local field.  
Let $A$ be a finite dimensional central simple $F$--algebra.  
Let $V$ be a simple left $A$--module.  
Then $\End_{A}(V)$ is a central division $F$--algebra.  
Let $D$ be the opposite algebra of $\End_{A}(V)$.  
Then $V$ is also a right $D$--module and we have $A \cong \End_{D}(V)$.  
Let $G$ be the multiplicative group of $A$.  
Then we have $G \cong \GL_{m}(D)$.  

For supercuspidal representations of $G$, some constructions of types are known.  
For example, Bushnell--Kutzko\cite{BK1} constructed types, called simple types,  for any irreducible supercuspidal representations when $G=\GL_{N}(F)$. 
S\'echerre--Stevens\cite{SS} extended the construction of simple types to any irreducible supercuspidal representations of any inner form $G$ of $\GL_{N}(F)$.  
For some "tame" irreducible supercuspidal representations of general $p$-adic reductive groups, Yu\cite{Yu} obtained the construction of types.  
Then it is a natural question whether there exists some relationship between these constructions of types.  

In the case $G=\GL_{N}(F)$, Mayeux\cite{May} compared Bushnell--Kutzko's simple type and Yu's type.  
In both constructions, we can obtain pairs consisting of some subgroup in $G$ which is compact modulo center and an irreducible representation of this subgroup such that the compact induction of this representation to $G$ is irreducible and supercuspidal.  
To show the relationship between these constructions, Mayeux used the tameness of simple types.  

\begin{thm}[{{\cite[Corollary 11.1]{May}}}]
\begin{enumerate}
\item Let $(J, \lambda)$ be a simple type for some essentially tame supercuspidal representation $\pi$, and let $(\tilde{J}, \Lambda)$ be an extension of $(J, \lambda)$ such that $\pi \cong \cInd_{\tilde{J}}^{G} \Lambda$.  
Then there exists a Yu's datum $\Psi = \left(x, (G^{i})_{i=0}^{d}, (\rbf_{i})_{i=0}^{d}, (\Phibf_{i})_{i}^{d}, \rho \right)$ such that $J={}^{\circ}K^{d}(\Psi)$, $\tilde{J} = K^{d}(\Psi)$ and $\rho^{d}(\Psi) \cong \Lambda$.  
\item On the other hand, let $\Psi = \left(x, (G^{i})_{i=0}^{d}, (\rbf_{i})_{i=0}^{d}, (\Phibf_{i})_{i}^{d}, \rho \right)$ be a Yu's datum such that $\pi' = \cInd_{K^{d}}^{G} \rho^{d}(\Psi)$ is irreducible and supercuspidal.  
Then there exists a tame simple type $(J, \lambda)$ and its maximal extension $(\tilde{J}, \Lambda)$ such that $\tilde{J} = K^{d}(\Psi)$ and $\rho^{d}(\Psi) \cong \Lambda$.  
\item For $G=\GL_{N}(F)$, the set of essentially tame supercuspidal representations of $G$ is equal to the set of tame supercuspidal representations of $G$ defined by Yu.  
\end{enumerate}
\end{thm}

In this paper, we extend Mayeux's result to any inner forms of $\GL_{N}(F)$.  
In the case $G=\GL_{m}(D)$, we have to consider some differences between Bushnell--Kutzko's types and S\'echerre--Stevens's types.  
Simple types are constructed using the data from some 4-tuple called simple strata.  
Strictly speaking, by fixing a simple stratum $[\Afra, n, 0, \beta]$, we can construct the group $J=J(\beta, \Afra)$ and finitely many candidates of $J$-representations.  
When we choose a representation $\lambda$ from these candidates, we obtain a simple type $(J, \lambda)$.  
In considering a maximal extension $(\tilde{J}, \Lambda)$ of $(J, \lambda)$, the group $\tilde{J}$ is independent of the choice of $\lambda$ when $G=\GL_{N}(F)$.  
Then $\tilde{J}$ canonically coincides with some subgroup $K^{d}$ in $G$ determined by some Yu's datum.  
However, if $G \neq \GL_{N}(F)$, $\tilde{J}$ may depend on the choice of $\lambda$.  
When we relate a Yu's datum $\Psi$ to $(J, \lambda)$ in the same way as $\GL_{N}(F)$ case, the subgroup $K^{d}(\Psi)$ in $G$ can properly contain $\tilde{J}$.  
Therefore we have to fill the gap between $K^{d}(\Psi)$ and $\tilde{J}$.  

\begin{thm}[Theorem \ref{Main1}]
Let $(J, \lambda)$ be a simple type for some essentially tame supercuspidal representation $\pi$, and let $(\tilde{J}, \Lambda)$ be an extension of $(J, \lambda)$ such that $\pi \cong \cInd_{\tilde{J}}^{G} \Lambda$.  
Then there exists a Yu's datum $\Psi=(x, (G^{i})_{i=0}^{d}, (\rbf_{i})_{i=0}^{d}, (\Phibf_{d})_{i=0}^{d}, \rho)$ satisfying the following conditions:  
\begin{enumerate}
\item $J={}^{\circ}K^{d}(\Psi)$, 
\item $\tilde{J} \subset K^{d}(\Psi)$, and
\item $\rho_{d}(\Psi) \cong \cInd_{\tilde{J}}^{K^{d}(\Psi)} \Lambda$.  
\end{enumerate}
\end{thm}

\begin{thm}[Theorem \ref{Main2}]
Conversely, let $\Psi=\left( x, (G^{i}), (\rbf_{i}), (\Phibf_{i}), \rho \right)$ be a Yu datum of $G$.  
Then there exists a tame simple type $(J, \lambda)$ and a maximal extension $(\tilde{J}, \Lambda)$ of $(J, \lambda)$ such that
\begin{enumerate}
\item ${}^{\circ}K^{d}(\Psi) = J$, 
\item $K^{d}(\Psi) \supset \tilde{J}$, and
\item $\rho_{d}(\Psi) \cong \cInd_{\tilde{J}}^{K^{d}(\Psi)} \Lambda$.  
\end{enumerate}
\end{thm}

By these theorems, we obtain the following corollary.  

\begin{cor}[Corollary \ref{tame_and_esstame}]
For any inner form $G$ of $\GL_{N}(F)$, the set of essentially tame supercuspidal representations of $G$ is equal to the set of tame supercuspidal representations of $G$ defined by Yu \cite{Yu}.  
\end{cor}

We sketch the outline of this paper.  
First, in \S \ref{SecherreStevens} and \ref{Yu's_types}, we recall constructions of types.  
We explain simple types of $G$ by S\'echerre--Stevens in \S \ref{SecherreStevens} and Yu's construction of tame supercuspidal representations in \S \ref{Yu's_types}.  
Next, in \S \ref{TameSimple}-\ref{Factorization}, we prepare ingredients to compare two constructions.  
A class of simple types corresponding to Yu's type is defined in \S \ref{TameSimple}.  
In \S \ref{TTLS}, we determine tame twisted Levi subgroups in $G$.  
For some tame twisted Levi subgroup $G'$ in $G$ and some ``nice" $x \in \Bscr^{E}(G', F)$, we obtain another description of Moy--Prasad filtration on $G'(F)$ attached to $x$, using hereditary orders, in \S \ref{BroussousLemaire}.  
Then we can compare the groups which types are defined over.  
In \S \ref{Generic}, we discuss generic elements and generic characters.  
We see that generic characters have information on a defining sequence of some simple stratum.  
In \S \ref{Depth0} we show some lemmas on simple types of depth zero.  
These lemmas are used to take ``depth-zero" parts of types.  
In \S \ref{Factorization}, we represent a simple character with a tame simple stratum as a product of characters.  
This factorization is needed to construct generic characters.  
Finally, in \S \ref{SStoYu} and \ref{YutoSS}, we prove the main theorem.  
(Compact inductions of) tame simple types are constructed from a Yu datum in \S \ref{SStoYu}.  
Conversely, Yu's types are constructed from tame simple types in \S \ref{YutoSS}.  
\bigbreak

\noindent{\bfseries Acknowledgment}\quad
I am deeply grateful to my supervisor Naoki Imai for his enormous support and helpful advice.  
He also checked the draft of this paper and pointed out mistakes.  
I truly appreciate Arnaud Mayeux sending the former version of \cite{May}.  
I would like to thank Vincent S\'echerre for the insightful discussion with him.  
I also would like to thank Kazuki Tokimoto for giving helpful comments.  
I am supported by the FMSP program at Graduate School of Mathematical Sciences, the University of Tokyo.  
\bigbreak

\noindent{\bfseries Notation}\quad
In this paper, we consider smooth representations over $\C$.  
We fix a non-archimedean local field $F$.  
For a finite-dimensional division algebra $D$ over $F$, let $\ofra_{D}$ be the ring of integers, $\pfra_{D}$ be the maximal ideal of $\ofra_{D}$, and let $k_{D}$ be the residual field of $D$.  
We fix a smooth, additive character $\psi:F \to \C^{\times}$ of conductor $\pfra_{F}$.  
For a finite field extension $E/F$, let $v_{E}$ be the unique surjective valuation $E \to \Z \cup \{ \infty \}$.  
Moreover, for any element $\beta$ in some algebraic extension field of $F$, we put $\ord(\beta) = e(F[\beta]/F)^{-1}v_{F[\beta]}(\beta)$.  

If $K$ is a field and $G$ is a $K$-group scheme, then $\underline{\Lie}(G)$ denotes the Lie algebra functor and we put $\Lie(G)=\underline{\Lie}(G)(K)$.  
If a $K$-group scheme is denoted by a capital letter $G$, the Lie algebra functor of $G$ is denoted by the same small Gothic letter $\gfra$.  
We also denote by $\Lie^{*}(G)$ or $\gfra^{*}(K)$ the dual of $\Lie(G)=\gfra(K)$.  
For connected reductive group $G$ over $F$, we denote by $\Bscr^{E}(G, F)$ (resp. $\Bscr^{R}(G, F)$) the enlarged Bruhat--Tis building (resp. the reduced Bruhat--Tis building) of $G$ over $F$ defined in \cite{BT}, \cite{BT2}.  
If $x \in \Bscr^{E}(G, F)$, we denote by $[x]$ the image of $x$ via the canonical surjection $\Bscr^{E}(G, F) \to \Bscr^{R}(G, F)$.  
The group $G(F)$ acts on $\Bscr^{E}(G, F)$ and $\Bscr^{R}(G, F)$.  
For $x \in \Bscr^{E}(G, F)$, let $G(F)_{x}$ (reps. $G(F)_{[x]}$) denote the stabilizer of $x \in \Bscr^{E}(G, F)$ (resp. $[x] \in \Bscr^{R}(G, F)$).  
We denote by $\tilde{\R}$ the totally ordered commutative monoid $\R \cup \{ r+ \mid r \in \R \}$.  
When $G$ splits over some tamely ramified extension of $F$, for $x \in \Bscr^{E}(G, F)$ let $\{ G(F)_{x, r} \}_{r \in \tilde{\R}_{\geq 0}}$, $\{ \gfra(F)_{x,r} \}_{r \in \tilde{\R}}$ and $\{ \gfra^{*}(F)_{x,r} \}_{r \in \tilde{\R}}$ be the Moy--Prasad filtration \cite{MP1}, \cite{MP2} on $G(F)$, $\gfra(F)$ and $\gfra^{*}(F)$, respectively.  
Here, we have $\gfra^{*}(F)_{x,r} = \{ X^{*} \in \gfra^{*}(F) \mid X^{*}(\gfra(F)_{x, (-r)+}) \subset \pfra_{F} \}$ for $r \in \R$.  
If $G$ is a torus, Moy--Prasad filtrations are independent of $x$, and then we omit $x$.  

Let $G$ be a group, $H$ be a subgroup in $G$ and $\lambda$ be a representation of $H$.  
Then we put ${}^{g}H=gHg^{-1}$ for $g \in G$, and we define a ${}^{g}H$-representation ${}^{g} \lambda$ as ${}^{g}\lambda(h) = \lambda(g^{-1}hg)$ for $h \in {}^{g}H$.  
Moreover, we also put
\[
	I_{G}( \lambda ) = \{ g \in G \mid \Hom_{H \cap {}^{g}H}(\lambda, {}^{g}\lambda) \neq 0 \}.  
\]

\section{Simple types by S\'echerre--Stevens}
\label{SecherreStevens}

We recall the theory of simple types by S\'echerre--Stevens from \cite{S1}, \cite{S2}, \cite{S3}, \cite{SS}.  

\subsection{Lattices, hereditary orders}

Let $D$ be a finite-dimensional central division $F$-algebra.  
Let $V$ be a right $D$-module with $\dim_{F}V < \infty$.  
We put $A=\End_{D}(V)$, and then $A$ is a central simple $F$-algebra.  
Moreover, there exists $m \in \Z_{>0}$ such that $A \cong \M_{m}(D)$.  
Let $G$ be the multiplicative group of $A$, and then $G$ is isomorphic to $\GL_{m}(D)$.  
We also put $d = (\dim _{F} D)^{1/2}$ and $N=md$.  

An $\ofra_{D}$-submodule $\Lcal$ in $V$ is called a lattice if and only if $\Lcal$ is a compact open submodule.  

\begin{defn}[{{\cite[D\'efinition 1.1]{S1}}}]
For $i \in \Z$, let $\Lcal_{i}$ be a lattice in $V$.  
We say that $\Lcal = (\Lcal_{i})_{i \in \Z}$ is an $\ofra_{D}$-sequence if
\begin{enumerate}
\item $\Lcal_{i} \supset \Lcal_{j}$ for any $i<j$, and
\item there exists $e \in \Z_{>0}$ that $\Lcal_{i+e} = \Lcal_{i}\pfra_{D}$ for any $i$.  
\end{enumerate}
The number $e=e(\Lcal)$ is called the period of $\Lcal$.  
An $\ofra_{D}$-sequence $\Lcal$ is called an $\ofra_{D}$-chain if $\Lcal_{i} \supsetneq \Lcal_{i+1}$ for every $i$.  
An $\ofra_{D}$-chain $\Lcal$ is called uniform if $[\Lcal_{i} : \Lcal_{i+1}]$ is constant for any $i$.  
\end{defn}

An $\ofra_{F}$-subalgebra $\Afra$ in $A$ is called a hereditary $\ofra_{F}$-order if every left and right ideal in $\Afra$ is $\Afra$-projective.  

We explain the relationship between $\ofra_{D}$-sequences in $V$ and hereditary $\ofra_{F}$-orders in $A$ from \cite[1.2]{S1}.  
Let $\Lcal = (\Lcal_{i})$ be an $\ofra_{D}$-sequence in $V$.  
For $i \in \Z$, we put
\[
	\Pfra_{i}(\Lcal) = \{ a \in A \mid a \Lcal_{j} \subset \Lcal_{i+j}, \, j \in \Z \}.  
\]
Then $\Afra = \Pfra_{0}(\Lcal)$ is a hereditary $\ofra_{F}$-order with the radical $\Pfra(\Afra) = \Pfra_{1}(\Lcal)$.  
For every hereditary $\ofra_{F}$-order $\Afra$ in $V$, there exists an $\ofra_{D}$-chain $\Lcal$ in $V$ such that $\Afra = \Afra(\Lcal)$.  
If $\Lcal$ is a uniform $\ofra_{D}$-chain, $\Afra = \Afra(\Lcal)$ is called principal.  

For any $\ofra_{D}$-chain $\Lcal = (\Lcal_{i})$, let $\Kfra(\Lcal)$ be the set of $g \in G$ such that there exists $n \in \Z$  satisfying $g \Lcal_{i} = \Lcal_{i+n}$ for any $i$.  
For the hereditary $\ofra_{F}$-order $\Afra=\Afra(\Lcal)$, let $\Kfra(\Afra)$ be the set of $g \in G$ such that $g\Afra g^{-1}=\Afra$.  
Then we have $\Kfra(\Afra) = \Kfra(\Lcal)$ and $\Kfra(\Afra)$ is compact modulo center.  

For $g \in \Kfra(\Afra)$, there exists a unique element $n \in \Z$ such that $g\Afra = \Pfra(\Afra)^{n}$.  
The map $g \mapsto n$ induces a group homomorphism $v_{\Afra} : \Kfra(\Afra) \to \Z$.  
Let $\U(\Afra)$ be the kernel of $v_{\Afra}$.  
Then we have $\U(\Afra) = \Afra^{\times}$ and $\U(\Afra)$ is the unique maximal compact open subgroup in $\Kfra(\Afra)$.  
We put $\U^{0}(\Afra) = \U(\Afra)$ and $\U^{n}(\Afra) = 1 + \Pfra(\Afra)^{n}$ for $n \in \Z_{>0}$.  
We also put $e(\Afra|\ofra_{F}) = v_{\Afra}(\varpi_{F})$, and then we have $e(\Afra|\ofra_{F})=de(\Lcal)$ for an $\ofra_{D}$-chain $\Lcal$ in $V$ such that $\Afra = \Afra(\Lcal)$.  

Let $E$ be an extension field of $F$ in $A$.  
Since $A$ is a central simple $F$-algebra, the centralizer $B=\Cent_{A}(E)$ of $E$ in $A$ is a central simple $E$-algebra.  
On the other hand, $V$ is equipped with an $E$-vector space structure via $E \subset A$.  
Since the actions of $E$ and $D$ in $V$ are compatible, $V$ is also equipped with a right $D \otimes_{F} E$-module structure, and then we have $B = \Cent_{A}(E) = \End_{D \otimes_{F} E}(V)$.  

Let $\Afra$ be a hereditary $\ofra_{F}$-order in $A$.  
The order $\Afra$ is called $E$-pure if we have $E^{\times} \subset \Kfra(\Afra)$.  

\begin{prop}[{{\cite[Theorem 1.3]{Br}}}]
For an $E$-pure hereditary $\ofra_{F}$-order $\Afra$ in $A$, the subring $\Bfra = \Afra \cap B$ in $B$ is a hereditary $\ofra_{E}$-order in $B$ with the radical $\Qfra = \Pfra(\Afra) \cap B$.  
\end{prop}

For any finite extension field $E$ of $F$, we put $A(E)=\End_{F}(E)$, and then $A(E)$ is a central simple $F$-algebra.  
The field $E$ is canonically embedded in $A(E)$ as a maximal subfield.  
By \cite[1.2]{BK1}, there exists a unique $E$-pure hereditary $\ofra_{F}$-order $\Afra(E)=\{ x \in A(E) \mid x(\pfra_{E}^{i}) \subset \pfra_{E}^{i}, \, i \in \Z \}$ in $A(E)$, which is associated with the $\ofra_{F}$-chain $(\pfra_{E}^{i})_{i \in \Z}$.  
Then we have $v_{\Afra(E)}(\beta) = v_{E}(\beta)$ for $\beta \in E^{\times}$.  

For $\beta \in \bar{F}$, we put $n_{F}(\beta) = -v_{F[\beta]}(\beta) = -v_{\Afra(F[\beta])}(\beta)$ as in \cite[2.3.3]{S1}.  

Let $\Afra$ be a hereditary $\ofra_{F}$-order in $A$ with the radical $\Pfra$.  
For non-negative integer $i, j$ with $\lfloor j/2 \rfloor \leq i \leq j$, the map $1+x \mapsto x$ induces the group isomorphism
\[
	\U^{i+1}(\Afra)/\U^{j+1}(\Afra) \cong \Pfra^{i+1}/\Pfra^{j+1}.  
\]
If $i$ and $j$ are as above and $c \in \Pfra^{-j}$, we can define a character $\psi_{c}$ of $\U^{i+1}(\Afra)$ as
\[
	\psi_{c}(1+x) = \psi \circ \Trd_{A/F}(cx)
\]
for $1+x \in \U^{i+1}(\Afra)$.  
We have $\psi_{c} = \psi_{c'}$ if and only if $c-c' \in \Pfra^{-i}$.  

\subsection{Strata, defining sequences of simple strata}

\begin{defn}[{{\cite[\S 2.1, Remarque 4.1]{S3}}}]
\begin{enumerate}
\item A 4-tuple $[\Afra, n, r, \beta]$ is called a stratum in $A$ if $\Afra$ is a hereditary $\ofra_{F}$-order in $A$, $n$ and $r$ are non-negative integer with $n \geq r$, and $\beta \in \Pfra(\Afra)^{-n}$.  
\item A stratum $[\Afra, n, r, \beta]$ is called pure if the followings hold:  
	\begin{enumerate}
	\item $E=F[\beta]$ is a field.  
	\item $\Afra$ is $E$-pure.  
	\item $n>r$.  
	\item $v_{\Afra}(\beta) = -n$.  
	\end{enumerate}
\item A stratum $[\Afra, n, r, \beta]$ is called simple if one of the followings holds:  
	\begin{enumerate}
	\item $n=r=0$ and $\beta \in \ofra_{F}$.  
	\item $[\Afra, n, r, \beta]$ is pure, and $r < -k_{0}(\beta, \Afra)$, where $k_{0}(\beta, \Afra) \in \Z \cup {-\infty}$ is defined as in \cite[\S 2.1]{S3} such that $k_{0}(\beta, \Afra)=-\infty$ if and only if $\beta \in F$, and $v_{\Afra}(\beta) \leq k_{0}(\beta, \Afra)$ for $\beta \notin F$.  
	\end{enumerate}
\end{enumerate}
\end{defn}

\begin{rem}
In \cite[\S 2.1]{S3}, simple strata are assumed be pure.  
By adding strata satisfying (3)(a) to simple strata, we can regard simple types of depth zero as coming from simple strata.  
\end{rem}

\begin{defn}
Strata $[\Afra, n, r, \beta]$ and $[\Afra, n, r, \beta']$ in $A$ are called equivalent if $\beta-\beta' \in \Pfra(\Afra)^{-r}$.  
\end{defn}

\begin{thm}[{{\cite[Th\'eor\`eme 2.2]{S3}}}]
\label{existapp}
Let $[\Afra, n, r, \beta]$ be a pure stratum.  
Then there exists $\gamma \in A$ such that $[\Afra, n, r, \gamma]$ is simple and equivalent to $[\Afra, n, r, \beta]$.  
\end{thm}

For $\beta \in \bar{F}$, we put $k_{F}(\beta) = k_{0}(\beta, \Afra(F[\beta]))$ as in \cite[2.3.3]{S1}.  

\begin{prop}[{{\cite[Proposition 2.25]{S1}}}]
\label{relofk0}
Suppose $E=F[\beta]$ can be embedded in $A$.  
We fix an embedding $E \hookrightarrow A$.  
Let $\Afra$ be an $E$-pure hereditary $\ofra_{F}$-order in $A$.  
Then we have $k_{0}(\beta, \Afra) = e(\Afra|\ofra_{F})e(E/F)^{-1}k_{F}(\beta)$.  
\end{prop}

The following lemma is used later.  

\begin{lem}
\label{compofk0}
Let $E/F$ be a field extension in $A$, and let $\Afra$ be an $E$-pure hereditary $\ofra_{F}$-order in $A$.  
Then, we have $k_{0}(\gamma, \Afra) = e(\Afra|\ofra_{F})e(E/F)^{-1}k_{0}(\gamma, \Afra(E))$ for any $\gamma \in E$.  
\end{lem}

\begin{prf}
First, by Proposition \ref{relofk0} we have 
\[
	k_{0}(\gamma, \Afra(E)) = e(\Afra(E)|\ofra_{F})e(F[\gamma]/F)^{-1}k_{F}(\gamma).  
\]
On the other hand, we also have $e(\Afra(E)|\ofra_{F})=e(E/F)$ by definition of $\Afra(E)$.  
Then we obtain
\begin{eqnarray*}
	e(\Afra|\ofra_{F})e(E/F)^{-1}k_{0}(\gamma, \Afra(E)) & = & e(\Afra|\ofra_{F})e(E/F)^{-1}e(E/F)e(F[\gamma]/F)^{-1}k_{F}(\gamma) \\
	& = & e(\Afra|\ofra_{F})e(F[\gamma]/F)^{-1}k_{F}(\gamma) \\
	& = & k_{0}(\gamma, \Afra), 
\end{eqnarray*}
where the last equality also follows from Proposition \ref{relofk0}.  
\end{prf}

\begin{defn}
An element $\beta \in \bar{F}$ is called minimal if $\beta \in F$ or $k_{F}(\beta) = -v_{F}(\beta)$.  
\end{defn}

\begin{defn}
Let $[\Afra, n, r, \beta]$ be a simple stratum in $A$.  
A sequence $\left( [\Afra, n, r_i, \beta_i] \right)_{i=0}^{s}$ is called a defining sequence of $[\Afra, n, r, \beta]$ if
\begin{enumerate}
\item $\beta_0=\beta, r_0=r$, 
\item $r_{i+1}=-k_{0}(\beta_{i}, \Afra)$ for $i=0, 1, \ldots, s-1$, 
\item $[\Afra, n, r_{i+1}, \beta_{i+1}]$ is simple and equivalent to $[\Afra, n, r_{i+1}, \beta_{i}]$ for $i=0, 1, \ldots, s-1$, 
\item $\beta_{s}$ is minimal over $F$.  
\end{enumerate}
\end{defn}

By Theorem \ref{existapp}, for any simple stratum $[\Afra, n, r, \beta]$ there exists a defining sequence of $[\Afra, n, r, \beta]$, as in the case $A$ is split over $F$.  

\subsection{Simple characters}

Let $[\Afra, n, 0, \beta]$ be a simple stratum in $A$.  
Then we can define compact open subgroups $J(\beta, \Afra)$ and $H(\beta, \Afra)$ in $\U(\Afra)$ as in \cite[\S3]{S1}.  
The subgroup $H(\beta, \Afra)$ in $\U(\Afra)$ is also contained in $J(\beta, \Afra)$.  
For $i \in \Z_{\geq 0}$, we put $J^{i}(\beta, \Afra) = J(\beta, \Afra) \cap \U^{i}(\Afra)$ and $H^{i}(\beta, \Afra) = H(\beta, \Afra) \cap \U^{i}(\Afra)$.  

\begin{lem}[{{\cite[\S3.3]{S1}}}]
Let $[\Afra, n, 0, \beta]$ be a simple stratum in $A$.  
If $\beta$ is not minimal, let $([\Afra, n, r_{i}, \beta_{i}])_{i=0}^{s}$ be a defining sequence of $[\Afra, n, 0, \beta]$.  
\begin{enumerate}
\item $J^{i}(\beta, \Afra)$ is normalized by $\Kfra(\Afra) \cap B^{\times}$ for any $i \in \Z_{\geq 0}$.  
\item We have $J(\beta, \Afra) =\U(\Bfra) J^{1}(\beta, \Afra)$, where $\Bfra = \Afra \cap B$.  
\item If $\beta$ is minimal, we have 
\[
	J^{1}(\beta, \Afra) = \U^{1}(\Bfra)\U^{\lfloor (n+1)/2 \rfloor}(\Afra), \, H^{1}(\beta, \Afra) = \U^{1}(\Bfra) \U^{\lfloor n/2 \rfloor +1}(\Afra).
\]  
\item If $\beta$ is not minimal, we have 
\[
	J^{t}(\beta, \Afra) = J^{t}(\beta_{1}, \Afra), \, H^{t'+1}(\beta, \Afra) = H^{t'+1}(\beta_{1}, \Afra), 
\]
where $t = \lfloor (-k_{0}(\beta, \Afra)+1)/2 \rfloor$ and $t' = \lfloor -k_{0}(\beta, \Afra) /2 \rfloor$.  
Moreover, we also have 
\[
	J^{1}(\beta, \Afra) = \U^{1}(\Bfra) J^{t}(\beta_{1}, \Afra), \, H^{1}(\beta, \Afra) = \U^{1}(\Bfra) H^{t'+1}(\beta_{1}, \Afra).  
\]
\end{enumerate}
\end{lem}

\begin{defn}[{{\cite[Proposition 3.47]{S1}}}]
\label{defofsimpch}
Let $[\Afra, n, 0, \beta]$ be a simple stratum.  
We put $q= -k_{0}(\beta, \Afra)$.  
Let $0 \leq t < q$ and we put $t'= \max \{ t, \lfloor q/2 \rfloor \}$.  
If $\beta$ is not minimal over $F$, we fix a defining sequence $\left( [\Afra, n, r_i, \beta_i] \right)_{i=0}^{s}$ of $[\Afra, n, 0, \beta]$.  
The set of simple characters $\mathscr{C}(\beta, t, \Afra)$ consists of characters $\theta$ of $H^{t+1}(\beta, \Afra)$ satisfying the following conditions:  
\begin{enumerate}
\item $\Kfra(\Afra) \cap B^{\times}$ normalizes $\theta$.  
\item $\theta | _{H^{t+1}(\beta, \Afra) \cap \U(\Bfra)}$ factors through $\Nrd_{B/E}$.  
\item If $\beta$ is minimal over $F$, we have $\theta |_{H^{t+1}(\beta, \Afra) \cap \U^{\lfloor n/2 \rfloor +1}(\Afra)} = \psi_{\beta}$.  
\item If $\beta$ is not minimal over $F$, there exists $\theta' \in \mathscr{C}(\beta_{1}, t', \Afra)$ such that $\theta |_{H^{t'+1}(\beta, \Afra)} = \psi_{\beta-\beta_{1}} \theta'$.  
\end{enumerate}
\end{defn}

\begin{rem}
\label{existsimpch}
This definition is well-defined and independent of the choice of a defining sequence by \cite[D\'efinition 3.45, Proposition 3.47]{S1}.  
Moreover, for any simple stratum $[\Afra, n, 0, \beta]$ the set $\mathscr{C}(\beta, 0, \Afra)$ is nonempty by \cite[Corollaire 3.35, D\'efinition 3.45]{S1}.  
\end{rem}

We recall the properties of $\mathscr{C}(\beta, 0, \Afra)$ from \cite{S2}.  
For $\theta \in \mathscr{C}(\beta, 0, \Afra)$, there exists an irreducible $J^{1}(\beta, \Afra)$-representation $\eta_{\theta}$ containing $\theta$, unique up to isomorphism.  
We call $\eta_{\theta}$ the Heisenberg representation of $\theta$.  
We have $\dim \eta_{\theta} = \left( J^{1}(\beta, \Afra) : H^{1}(\beta, \Afra) \right) ^{1/2}$.  
Moreover, there exists an extension $\kappa$ of $\eta_{\theta}$ to $J(\beta, \Afra)$ such that $I_{G}(\kappa)=J^{1}B^{\times}J^{1}$.  
We call $\kappa$ a $\beta$-extension of $\eta_{\theta}$.  
If $\kappa$ is a $\beta$-extension of $\eta_{\theta}$, then any $\beta$-extension of $\eta_{\theta}$ is the form $\kappa \otimes (\chi \circ \Nrd_{B/E})$, where $\chi$ is trivial on $1+\pfra_{E}$ and $\chi \circ \Nrd_{B/E}$ is regarded a character of $J(\beta, \Afra)$ via the isomorphism $J(\beta, \Afra)/J^{1}(\beta, \Afra) \cong \U(\Bfra)/\U^{1}(\Bfra)$.  

\subsection{Maximal simple types}

We state the definition of maximal simple types.  
Recall that for a simple stratum $[\Afra, n, 0, \beta]$ we put $E=F[\beta]$, $B=\Cent_{A}(E)$ and $\Bfra = \Afra \cap B$.  
Since $B$ is a central simple $E$-algebra, there exist $m_{E} \in \Z$ and a division $E$-algebra $D_{E}$ such that $B \cong \M_{m_{E}}(D_{E})$.  

\begin{defn}[{{\cite[\S2.5, \S4.1]{S3}}}]
\label{defofsimpletype}
A pair $(J, \lambda)$ consisting a compact open subgroup $J$ in $G$ and an irreducible $J$-representation $\lambda$ is called a maximal simple type if there exists a simple stratum $[\Afra, n, 0, \beta]$ and irreducible $J$-representations $\kappa$ and $\sigma$ satisfying the following assertions:  
\begin{enumerate}
\item $\Bfra$ is a maximal hereditary $\ofra_{E}$-order in $A$, that is, $\Bfra \cong \M_{m_{E}}(\ofra_{D_{E}})$.  
\item $J = J(\beta, \Afra)$.  
\item $\kappa$ is a $\beta$-extension of $\eta_{\theta}$ for some $\theta \in \mathscr{C}(\beta, 0, \Afra)$.  
\item $\sigma$ is trivial on $J^{1}(\beta, \Afra)$, and when we regard $\sigma$ as a $\GL_{m_{E}}(k_{D_{E}})$-representation via the isomorphism
\[
	J(\beta, \Afra)/J^{1}(\beta, \Afra) \cong \U(\Bfra)/\U^{1}(\Bfra) \cong \GL_{m_{E}}(k_{D_{E}}), 
\]
$\sigma$ is a cuspidal representation of $\GL_{m_{E}}(k_{D_{E}})$.  
\item $\lambda \cong \kappa \otimes \sigma$.  
\end{enumerate}
\end{defn}

\begin{rem}
Let $(J, \lambda)$ be a maximal simple type associated with a simple stratum $[\Afra, 0, 0, \beta]$.  
Then we have $E=F[\beta]=F$, $B = \Cent_{A}(F) = A$ and $\Bfra = \Afra \cap A = \Afra$.  
Since $(J, \lambda)$ is a maximal simple type, $\Afra$ is a maximal hereditary $\ofra_{F}$-order in $A$.  
Moreover, we have $J(\beta, \Afra)=\U(\Afra)$ and $H^{1}(\beta, \Afra) = \U^{1}(\Afra)$.  
Let $\kappa$ and $\sigma$ be as in Definition \ref{defofsimpletype}.  
Since we have $\mathscr{C}(\beta, 0, \Afra) = \{ 1 \}$, there exists a character $\chi$ of $F^{\times}$ such that trivial on $1+\pfra_{F}$ and $\kappa = \chi \circ \Nrd_{A/F}$.  
Then $\kappa \otimes \sigma$ is trivial on $\U^{1}(\Afra)$ and cuspidal as a $\GL_{m}(k_{D})$-representation.  
Therefore $(J, \lambda) = (\U(\Afra), \kappa \otimes \sigma)$ is nothing but the maximal simple type of level 0, defined in \cite[\S2.5]{S3}.  
\end{rem}

\begin{thm}[{{\cite[Theorem 5.5(ii)]{GSZ} and \cite[Th\'eor\`eme 5.21]{S3}}}]
Let $\pi$ be an irreducible representation of $G$.  
Then $\pi$ is supercuspidal if and only if there exists a maximal simple type $(J, \lambda)$ such that $\lambda \subset \pi|_{J}$.  
\end{thm}

We recall the construction of irreducible supercuspidal representations of $G$ from maximal simple types.  
Let $(J, \lambda)$ be a maximal simple type associated with a simple stratum $[\Afra, n, 0, \beta]$.  
Let $\kappa$ and $\sigma$ be as in Definition \ref{defofsimpletype}.  
Since $\Bfra$ is maximal, we have $\Kfra(\Bfra) = \Kfra(\Afra) \cap B^{\times}$ by \cite[Lemme 1.6]{S2}, and then $\Kfra(\Bfra)$ normalizes $J(\beta, \Afra)$.  

We fix $g \in \Kfra(\Bfra)$ with $v_{\Bfra}(g) = 1$.  
Since $g$ normalizes $J(\beta, \Afra)$, we can consider the twist ${}^{g} \sigma$ of $\sigma$ by $g$.  
Let $l_{0}$ be the smallest positive integer such that ${}^{g^{l_{0}}}\sigma \cong \sigma$.  
Then $\tilde{J}(\lambda) = I_{G}(\lambda)$ is the subgroup in $G$ generated by $J$ and $g^{l_{0}}$.  

\begin{thm}[{{\cite[Th\'eor\`eme 5.2]{S3}, \cite[Corollary 5.22]{SS}}}]
\begin{enumerate}
\item For any maximal simple type $(J, \lambda)$, there exists an extension $\Lambda$ of $\lambda$ to $\tilde{J}(\lambda)$.  
\item Let $(\tilde{J}(\lambda), \Lambda)$ be as above.  
Then $\cInd_{\tilde{J}(\lambda)}^{G} \Lambda$ is irreducible and supercuspidal.  
\item For any irreducible supercuspidal representation $\pi$ of $G$, there exists an extension $(\tilde{J}(\lambda), \Lambda)$ of a maximal simple type $(J, \lambda)$ such that $\pi = \cInd_{\tilde{J}(\lambda)}^{G} \Lambda$.  
\end{enumerate}
\end{thm}

\subsection{Concrete presentation of open subgroups}

In the above, we define open subgroups $H^{1}(\beta, \Afra), J(\beta, \Afra)$ and $\tilde{J}(\lambda)$.  
In this subsection, we define another subgroup $\hat{J}(\beta, \Afra)$ and obtain the concrete presentation of some groups, which is used later.  

\begin{defn}
\label{defofJhat}
Let $[\Afra, n, 0, \beta]$ be a simple stratum with $\Bfra$ is maximal.  
Then we put $\hat{J}(\beta, \Afra) = \Kfra(\Bfra)J(\beta, \Afra)$.  
\end{defn}

\begin{rem}
\begin{enumerate}
\item Since $\Kfra(\Bfra)$ normalizes $J(\beta, \Afra)$, the set $\hat{J}(\beta, \Afra)$ is also a subgroup in $G$.  
We have $\Kfra(\Bfra) \cap J(\beta, \Afra) = \U(\Bfra)$, and then 
\[
	\hat{J}(\beta, \Afra)/J(\beta, \Afra) \cong \Kfra(\Bfra)/\U(\Bfra) \cong \Z.  
\]
\item Let $(J, \lambda)$ be a maximal simple type associated with $[\Afra, n, 0, \beta]$.  
Then we have $\tilde{J}(\lambda) \subset \hat{J}(\beta, \Afra)$.  
The group $\hat{J}(\beta, \Afra)$ only depends on $[\Afra, n, 0, \beta]$, while $\tilde{J}(\lambda)$ also depends on $\lambda$.  
\end{enumerate}
\end{rem}

We describe $H^{1}(\beta, \Afra), J(\beta, \Afra)$ and $\hat{J}(\beta, \Afra)$ concretely, using a defining sequence $([\Afra, n, r_{i}, \beta_{i}])_{i=0}^{s}$ of $[\Afra, n, 0, \beta]$.  
We put $B_{\beta_{i}} = \Cent_{A}(F[\beta_{i}])$ for $i=0, \ldots, s$.  

\begin{lem}
\label{presenofHJ}
Let $[\Afra, n, 0, \beta]$ be a maximal simple stratum of $A$ and $\left( [\Afra, n, r_i, \beta_i] \right) _{i=0}^{s}$ be a defining sequence of $[\Afra, n, 0, \beta]$.  
Then we have following concrete presentations of groups:  
\begin{enumerate}
\item $H^{1}(\beta, \Afra)=\left( B_{\beta_{0}}^{\times} \cap \U^{ \lfloor \frac{r_0}{2} \rfloor +1}(\Afra) \right) \cdots \left( B_{\beta_{s}} ^{\times} \cap \U^{ \lfloor \frac{r_{s}}{2} \rfloor +1}(\Afra) \right) \U^{\lfloor \frac{n}{2} \rfloor +1}(\Afra)$.  
\item $J(\beta, \Afra)=\U(\Bfra) \left( B_{\beta_{1}}^{\times} \cap \U^{ \lfloor \frac{r_1+1}{2} \rfloor}(\Afra) \right) \cdots \left( B_{\beta_{s}} ^{\times} \cap \U^{ \lfloor \frac{r_{s}+1}{2} \rfloor}(\Afra) \right) \U^{ \lfloor \frac{n+1}{2} \rfloor}(\Afra)$.  
\item $\hat{J}(\beta, \Afra)=\Kfra(\Bfra) \left( B_{\beta_{1}}^{\times} \cap \U^{ \lfloor \frac{r_1+1}{2} \rfloor}(\Afra) \right) \cdots \left( B_{\beta_{s}} ^{\times} \cap \U^{ \lfloor \frac{r_{s}+1}{2} \rfloor}(\Afra) \right) \U^{ \lfloor \frac{n+1}{2} \rfloor}(\Afra)$.  
\end{enumerate}
\end{lem}

\begin{prf}
We show (1) by induction on the length $s$ of a defining sequence.  
When $s=0$, that is, $\beta$ is minimal over $F$, then $H^{1}(\beta, \Afra) = \U^{1}(\Bfra) \U^{\lfloor n/2 \rfloor+1}(\Afra)$.  
Since we have $\U^{1}(\Bfra)=1+(B \cap \Pfra)=B \cap (1+\Pfra) = B \cap \U^{1}(\Afra)$ and $r_0=0$, the equality in (1) for minimal $\beta$ holds.  
Suppose $s > 0$, that is, $\beta$ is not minimal over $F$.  
Then $H^{1}(\beta, \Afra) = \U^{1}(\Bfra) H^{\lfloor r_1/2 \rfloor +1}(\beta_1, \Afra)$.  
By induction hypothesis, we have 
\[
	H^{1}(\beta_1, \Afra)=\U^1(\Bfra_{\beta_1}) \left( B_{\beta_{2}}^{\times} \cap \U^{ \lfloor \frac{r_2}{2} \rfloor +1}(\Afra) \right) \cdots \left( B_{\beta_{s}} ^{\times} \cap \U^{ \lfloor \frac{r_{s}}{2} \rfloor +1}(\Afra) \right) \U^{\lfloor \frac{n}{2} \rfloor +1}(\Afra).  
\]
Since $r_1 < r_2 < \ldots < r_s < n$, we have $\lfloor r_1/2 \rfloor + 1 \leq \lfloor r_2 /2 \rfloor + 1 \leq \ldots \leq \lfloor r_s /2 \rfloor + 1 \leq \lfloor n/2 \rfloor + 1$ and
\[
	 B_{\beta_{2}}^{\times} \cap \U^{ \lfloor \frac{r_2}{2} \rfloor +1}(\Afra) , \cdots , B_{\beta_{s}} ^{\times} \cap \U^{ \lfloor \frac{r_{s}}{2} \rfloor +1}(\Afra) , \U^{\lfloor \frac{n}{2} \rfloor +1}(\Afra) \subset \U^{\lfloor \frac{ r_{1}}{2} \rfloor +1} (\Afra).   
\]
Therefore we obtain
\begin{eqnarray*}
	H^{\lfloor \frac{r_{1}}{2} \rfloor +1}(\beta_1, \Afra) & = & \left( \U^1(\Bfra_{\beta_1}) \cap \U^{\lfloor \frac{r_{1}}{2} \rfloor +1} (\Afra) \right) \left( B_{\beta_{2}}^{\times} \cap \U^{ \lfloor \frac{r_2}{2} \rfloor +1}(\Afra) \right) \cdots \\
	 & & \hspace{80pt} \cdots \left( B_{\beta_{s}} ^{\times} \cap \U^{ \lfloor \frac{r_{s}}{2} \rfloor +1}(\Afra) \right) \U^{\lfloor \frac{n}{2} \rfloor +1}(\Afra) \\
	 & = & \left( B_{\beta_{1}}^{\times} \cap \U^{ \lfloor \frac{r_1}{2} \rfloor +1}(\Afra) \right) \cdots \left( B_{\beta_{s}} ^{\times} \cap \U^{ \lfloor \frac{r_{s}}{2} \rfloor +1}(\Afra) \right) \U^{\lfloor \frac{n}{2} \rfloor +1}(\Afra), 
\end{eqnarray*}
and the equality in (1) for non-minimal $\beta$ also holds.  

Similarly, we can show that 
\[
	J^{1}(\beta, \Afra)=\U^{1}(\Bfra) \left( B_{\beta_{1}}^{\times} \cap \U^{ \lfloor \frac{r_1+1}{2} \rfloor}(\Afra) \right) \cdots \left( B_{\beta_{s}} ^{\times} \cap \U^{ \lfloor \frac{r_{s}+1}{2} \rfloor}(\Afra) \right) \U^{ \lfloor \frac{n+1}{2} \rfloor}(\Afra).  
\]
Then (2) and (3) are deduced from the fact $J(\beta, \Afra)=\U(\Bfra) J^{1}(\beta, \Afra)$ and $\hat{J}(\beta, \Afra) = \Kfra(\Bfra)J(\beta, \Afra)$.  
\end{prf}

\section{Yu's construction of types for tame supercuspidal representations}
\label{Yu's_types}

In this section, we recall how to construct Yu's types from \cite{Yu}.  
Let $G$ be a connected reductive group over $F$.  

\subsection{Admissible sequences}

\begin{defn}
Let $(G^{i})=(G^{0}, \ldots, G^{d})$ be a sequence of group subscheme in $G$ over $F$.  
We call $(G^{i})$ is a tame twisted Levi sequence if $G^{0} \subset G^{1} \subset \cdots \subset G^{d}=G$ and there exists a tamely ramified extension $E$ of $F$ such that $G^{i} \times _{F} E$ is a split Levi subgroup in $G \times _{F} E$ for $i=0, \ldots, d$.  
\end{defn}

Let $\vec{G}=(G^{0}, \ldots, G^{d})$ be a tame twisted Levi sequence in $G$.  
Then there exist a maximal torus $T$ in $G^{0}$ over $F$ and a tamely ramified, finite Galois extension $E$ over $F$ such that $T \times_{F} E$ is split.  
For $i=0, \ldots, d$, we put $\Phi_{i}=\Phi(G^{i}, T; E) \cup \{0\}$.  
For $\alpha \in \Phi_{d} \setminus \{ 0\} = \Phi(G, T; E)$, we denote by $G_{\alpha}$ the root subgroup in $G_{E}$ defined by $\alpha$.  
Let $G_{\alpha}=T$ if $\alpha = 0$.  
Let $\gfra_{\alpha}$ be the Lie algebra of $G_{\alpha}$, which is a Lie subalgebra in $\gfra_{E}$, and let $\gfra_{\alpha}^{*}$ be its dual.  

Let $\vec{\rbf}=(\rbf_{0}, \rbf_{1}, \ldots, \rbf_{d}) \in \tilde{\R}^{d+1}$.  
Then we can define a map $f_{\vec{\rbf}} : \Phi_{d} \to \tilde{\R}$ by $f_{\vec{\rbf}}(\alpha) = \rbf_{i}$ if $i = \min \{ j \mid \alpha \in \Phi_{j} \}$.  

A sequence $\vec{\rbf}=(\rbf_{0}, \ldots, \rbf_{d}) \in \tilde{\R}^{d+1}$ is called an admissible sequence if and only if there exists $\nu \in \{ 0, \ldots, d \}$ such that 
\[
	0 \leq \rbf_{0} = \ldots = \rbf_{\nu}, \, \frac{1}{2}\rbf_{\nu} \leq \rbf_{\nu+1} \leq \ldots \leq \rbf_{d}.  
\]

Let $x$ be in the apartment $A(G, T, E) \subset \Bscr^{E}(G, E)$.  
Then we can determine the filtrations $\{ G_{\alpha}(E)_{x,r} \}_{r \in \tilde{\R}_{\geq 0}}$ on $G_{\alpha}(E)$, $\{ \gfra_{\alpha}(E)_{x,r} \}_{r \in \tilde{\R}}$ on $\gfra_{\alpha}(E)$, and $\{ \gfra_{\alpha}^{*}(E)_{x,r} \}_{r \in \tilde{\R}}$ on $\gfra_{\alpha}^{*}(E)$.  

We denote by $\vec{G}(E)_{x, \vec{\rbf}}$ the subgroup in $G(E)$ generated by $G_{\alpha}(E)_{x, f_{\vec{\rbf}}(\alpha)}$ ($\alpha \in \Phi_{d}$).  

By taking $x \in A(G, T, E) \cap \Bscr^{E}(G, E)$, we can determine a valuation on the root datum of $(G, T, E)$ in the sense of \cite{BT}.  
By restricting this valuation,  we can also define a valuation on the root datum of $(G^{i}, T, E)$.  
Then we can determine $x_{i} \in \Bscr^{E}(G^{i}, E)$ by the valuation, uniquely up to $X^{*}(G^{i}) \otimes \R$.  
When we take $x_{i}$ in such a way, we can determine an affine, $G^{i}(E)$-equivalent embedding $j_{i} : \Bscr^{E}(G^{i}, E) \to \Bscr^{E}(G, E)$ such that $j_{i}(x_{i})=x$.  
We identify $x_{i}$ with $x$ via $j_{i}$.  

To consider subgroups in $G(F)$, we also assume $x \in \Bscr^{E}(G, E)^{\Gal(E/F)}$, that is, $x \in \Bscr^{E}(G, F)$.  
Then we can determine the Moy--Prasad filtration on $G^{i}(F), \gfra^{i}(F)$ and $(\gfra^{i})^{*}(F)$ by $x$.  
We put $\vec{G}(F)_{x,\vec{\rbf}} = \vec{G}(E)_{x, \vec{\rbf}} \cap G(F)$.  

\begin{prop}[{{\cite[2.10]{Yu}}}]
The group $\vec{G}(F)_{x, \vec{\rbf}}$ is independent of the choice of $T$.  
If $\vec{\rbf}$ is increasing with $\rbf_{0} > 0$, then we have 
\[
	\vec{G}(F)_{x, \vec{\rbf}} = G^{0}(F)_{x, \rbf_{0}} G^{1}(F)_{x, \rbf_{1}} \cdots G^{d}(F)_{x, \rbf_{d}}.  
\]
\end{prop}

\subsection{Generic elements, generic characters}
Let $r \in \tilde{\R}_{>0}$ and let $r' \in \R$ with $(r'/2)+ \leq r \leq  r'+$.  
We put $G(F)_{x,r:r'}=G(F)_{x,r}/G(F)_{x,r'}$ and $\gfra(F)_{x,r:r'}=\gfra(F)_{x,r}/\gfra(F)_{x,r'}$.  
Then we have a group isomorphism $G(F)_{x,r:r'} \cong \gfra(F)_{x,r:r'}$.  

Let $S$ be a subgroup of $G(F)$ between $G(F)_{x,r/2+}$ and $G(F)_{x,r+}$, and let $\sfra$ be the sublattice of $\Lie(G)$ between $\gfra(F)_{x,r/2+}$ and $\gfra(F)_{x,r+}$ such that $\sfra/\gfra(F)_{x,r+} \cong S/G(F)_{x,r+}$.  

\begin{defn}
A character $\Phi$ of $S/G(F)_{x,r+}$ is realized by $X^{*} \in \Lie^{*}(G)_{x,-r}$ if $\Phi$ is equal to 
\[
	\xymatrix{
		S/G(F)_{x,r+} \cong \sfra/\gfra(F)_{x,r+} \ar[r]^-{X^{*}} & F \ar[r]^-{\psi} & \C^{\times} \\
	}.  
\]
\end{defn}

Let $G'$ be a tame twisted Levi subgroup in $G$.   
The Lie algebra $\Lie(G')$ and its dual $\Lie^{*}(G')$ are equipped with a $G'(F)$-action by conjugation.  
Let $Z(G')^{\circ}$ be the connected component of the center of $G$.   Then $\Lie^{*}(Z(G')^{\circ})$ is identified with the $G(F)$-fixed part of $\Lie^{*}(G')$.  
In this way we regard $\Lie^{*}(Z(G')^{\circ})$ as a subspace of $\Lie^{*}(G)$.  
  
To define generic characters of $G'$, we define generic elements in $\Lie^{*}(Z(G')^{\circ})$.  
Then we consider the conditions \textbf{GE1} and \textbf{GE2}.  

Let $E$ be a finite, tamely ramified extension of $F$ and $T$ be an $F$-torus in $G'$ such that $T \times_{F} E$ is maximal and split. 
Let $\alpha \in \Phi(G, T; \bar{F})$.   
Then the derivation $\dr \check{\alpha}$ is an $\bar{F}$-linear map from $\underline{\Lie} (\mathbb{G}_{m})(\bar{F}) \cong \bar{F}$ to $\Lie(T \times_{F} \bar{F})$.  
We obtain $H_{\alpha}=\dr \check{\alpha}(1)$ as an element in $\Lie(T \times_{F} \bar{F})$.  

Here, we recall the condition \textbf{GE1}.  
Let $X^{*} \in \Lie^{*}(Z(G')^{\circ})$.  
Then we can regard $X^{*} \in \Lie^{*}(G')$ as above.  
We put $X_{\bar{F}}^{*}=X^{*} \otimes_{F} 1 \in \Lie^{*}(G') \otimes_{F} \bar{F} = \Lie^{*}(G' \times_{F} \bar{F})$.  
Since $T \subset G'$, we have $H_{\alpha} \in \Lie(G' \times_{F}\bar{F}) = \Lie(G') \otimes_{F} \bar{F}$.  
Therefore we obtain $X_{\bar{F}}^{*} (H_{\alpha}) \in \bar{F}$.  

\begin{defn}
Let $X^{*} \in \Lie^{*}(Z(G')^{\circ})_{-r}$ for some $r \in \R$.  
We say $X^{*}$ satisfies \textbf{GE1} with depth $r$ if $\ord \left( X_{\bar{F}}^{*} (H_{\alpha}) \right) = -r$ for all $\alpha \in \Phi(G, T; \bar{F}) \setminus \Phi(G', T; \bar{F})$.  
\end{defn}

We also have to consider the condition \textbf{GE2} defined in \cite[\S 8]{Yu}.  
However, in our case if \textbf{GE1} holds, then \textbf{GE2} automatically holds.  

\begin{prop}[{{\cite[Lemma 8.1]{Yu}}}]
If the residual characteristic of $F$ is not a torsion prime for the root datum of $G$, then \textbf{GE1} implies \textbf{GE2}. 
\end{prop}

\begin{prop}[{{\cite[Corollary 1.13]{Ro}}}]
If a root datum is type A, then the set of torsion primes for the datum is empty.  
\end{prop}

From these propositions, we obtain the following corollary.  

\begin{cor}
\label{GE1toGE2}
If the root datum of $G$ is type A, then \textbf{GE1} implies \textbf{GE2}.    
\end{cor}

\begin{defn}
Let $X^{*} \in \Lie^{*}(Z(G')^{\circ})_{-r}$ for some $r \in \R$. 
The linear form $X^{*}$ is called $G$-generic of depth $r$ if and only if conditions \textbf{GE1} and \textbf{GE2} hold.  
\end{defn}

Eventually, we can define generic characters.  

\begin{defn}
A character $\Phi$ of $G'$ is called $G$-generic relative to $x$ of  depth $r \in \R_{>0}$ if $\Phi |_{G'(F)_{x,r+}}$ is trivial, $\Phi |_{G'(F)_{x,r}}$ is non-trivial, and there exists a $G$-generic element $X^{*} \in \Lie^{*}(Z(G')^{\circ})_{-r} \subset \Lie^{*}(G')_{x,-r}$ with depth $r$ such that $\Phi$ is realized by $X^{*}$ when $\Phi$ is regarded as a character of $G'(F)_{x,r:r+}$.  
\end{defn}

\subsection{Yu data}

Let $d \in \Z_{\geq 0}$.  

A 5-tuple $\Psi = \left( x, (G^{i})_{i=0}^{d}, (\rbf_{i})_{i=0}^{d}, (\Phibf_{i})_{i=0}^{d}, \rho \right)$ is called a Yu datum if $\Psi$ satisfies the following conditions:  

\begin{itemize}
\item The sequence $(G^{i})_{i=0}^{d}$ is a tame twisted Levi sequence such that $Z(G^{i})/Z(G)$ is anisotropic for $i=0, \ldots, d$ and 
\[
	G^{0} \subsetneq G^{1} \subsetneq \cdots \subsetneq G^{d} = G.  
\]
\item We have $x \in \Bscr^{E}(G^{0}, F) \cap A(G, T, E)$, where $T$ is a maximal $F$-torus in $G$ which splits over some tamely ramified extension of $F$.  
\item For $i=0, \ldots, d$, the number $\rbf_{i} \in \R$ such that 
\[
	0=\rbf_{-1} < \rbf_{0} < \ldots < \rbf_{d-1} \leq \rbf_{d}.  
\]
\item For $i=0, \ldots, d-1$, the character $\Phibf_{i}$ of $G^{i}(F)$ is $G^{i+1}$-generic relative to $x$ of depth $\rbf_{i}$.  
If $\rbf_{d-1} \neq \rbf_{d}$, the character $\Phibf_{d}$ of $G^{d}(F)$ is of depth $\rbf_{d}$.  
If $\rbf_{d-1} = \rbf_{d}$, the character $\Phibf_{d}$ of $G^{d}(F)$ is trivial.  
\item The irreducible representation $\rho$ of $G^{0}(F)_{[x]}$ is trivial on $G^{0}(F)_{x,0+}$ but nontrivial on $G^{0}(F)_{x}$, and $\cInd_{G^{0}(F)_{[x]}}^{G^{0}(F)} \rho$ is irreducible and supercuspidal.  
\end{itemize}

\subsection{Yu's construction}
\label{Yuconst}
In this subsection, we construct Yu's type by using some data from a Yu datum.  
Let $\Psi = \left( x, (G^{i})_{i=0}^{d}, (\rbf_{i})_{i=0}^{d}, (\Phibf_{i})_{i=0}^{d}, \rho \right)$ be a Yu datum.  

First, Yu constructed subgroups in $G$, which some representations are defined over.  
\begin{defn}
\label{defofKi}
For $i=0, \ldots, d$, let $\mathbf{s}_{i} = \rbf_{i}/2$.  
\begin{enumerate}
\item 
$\begin{array}{rcl}
	K_{+}^{i} & = & G^{0}(F)_{x, 0+} G^{1}(F)_{x,\mathbf{s}_{0}+} \cdots G^{i}(F)_{x, \mathbf{s}_{i-1}+} \\
	& = & (G^{0}, \ldots, G^{i})(F)_{x,(0+, \mathbf{s}_{0}+, \ldots, \mathbf{s}_{i-1}+)}.  
\end{array}$
\item 
$\begin{array}{rcl}
	{}^{\circ}K^{i} & = & G^{0}(F)_{x, 0} G^{1}(F)_{x, \mathbf{s}_0} \cdots G^{i}(F)_{x, \mathbf{s}_{i-1}} \\
	& = & G^{0}(F)_{x, 0} (G^{1}, \ldots, G^{i})(F)_{x, (\mathbf{s}_{0}, \ldots, \mathbf{s}_{i-1})}.  
\end{array}$
\item $K^{i}=G^{0}(F)_{[x]}G^{1}(F)_{x, \mathbf{s}_0} \cdots G^{i}(F)_{x, \mathbf{s}_{i-1}}=G^{0}(F)_{[x]} {}^{\circ}K^{i}$.  
\end{enumerate}
\end{defn}

\begin{prop}
For any $i=0, \ldots, d$, the groups $K_{+}^{i}$ and ${}^{\circ} K^{i}$ are compact, and $K^{i}$ is compact modulo center.  
\end{prop}

Yu also defined subgroups in $G(F)$, which ``fill the gap" between subgroups defined as above.  

\begin{defn}
For $i=1, \ldots, d$, 
\begin{enumerate}
\item $J^{i}=(G^{i-1}, G^{i})(F)_{x, (\rbf_{i-1}, \mathbf{s}_{i-1})}$, 
\item $J_{+}^{i}=(G^{i-1}, G^{i})(F)_{x, (\rbf_{i-1}, \mathbf{s}_{i-1}+)}$.  
\end{enumerate}
\end{defn}

Then, we have $K^{i}J^{i+1}=K^{i+1}$ and $K_{+}^{i}J_{+}^{i+1}=K_{+}^{i+1}$ for $i=0, \ldots, d-1$.  

Next, Yu defined characters $\hat{\Phibf}_{i}$ of $K_{+}^{d}$.  
The Lie algebra $\gfra(F)$ of $G(F)$ is equipped with a canonical $G(F)$-action.  
In particular, $Z(G^{i})^{\circ}(F)$ acts $\gfra(F)$ by restricting the $G(F)$-action.  
Then $Z(G^{i})^{\circ}(F)$-fixed part of $\gfra(F)$ is equal to the Lie algebra $\gfra^{i}(F)$ of $G^{i}(F)$.  
Moreover, we have a decomposition $\gfra(F) = \gfra^{i}(F) \oplus \mathfrak{n}^{i}(F)$ as a $Z(G^{i})^{\circ}(F)$-representation.  
This decomposition is well-behaved on the Moy--Prasad filtration:  we have $\gfra(F)_{x,s} = \gfra^{i}(F)_{x,s} \oplus \mathfrak{n}^{i}(F)_{x,s}$ for any $s \in \tilde{\R}$, where $\mathfrak{n}^{i}(F)_{x,s} \subset \mathfrak{n}^{i}(F)$.  
Let $\pi_{i}:\gfra(F) = \gfra^{i}(F) \oplus \mathfrak{n}^{i}(F) \to \gfra^{i}(F)$ be the projection.  
Then $\pi_{i}$ induces $\gfra(F)_{x,\mathbf{s}_{i}+:\rbf_{i}+} \to \gfra^{i}(F)_{x,\mathbf{s}_{i}+:\rbf_{i}+}$, and we obtain a group homomorphism 
\[
	\xymatrix{
		\tilde{\pi}_{i}:G(F)_{x,\mathbf{s}_{i}+} \ar[r] & G(F)_{x,\mathbf{s}_{i}+:\rbf_{i}+} \ar[r]^-{\pi_{i}} & G^{i}(F)_{x,\mathbf{s}_{i}+:\rbf_{i}+}.  
	}
\]
Here, Yu defined a character $\hat{\Phibf}_{i}$ of $K_{+}^{d}$ as
\begin{eqnarray*}
	\hat{\Phibf}_{i} |_{K_{+}^{d} \cap G^{i}(F)} & = & \Phibf_{i}, \\
	\hat{\Phibf}_{i} |_{K_{+}^{d} \cap G(F)_{x, \mathbf{s}_{i}+}} & = & \Phibf_{i} \circ \tilde{\pi}_{i}, 
\end{eqnarray*}
where $K^{d}_{+}$ is generated by $K_{+}^{d} \cap G^{i}(F)$ and $K_{+}^{d} \cap G(F)_{x, \mathbf{s}_{i}+}$ as we have
\begin{eqnarray*}
	K_{+}^{d} \cap G^{i}(F) & = & G^{0}(F)_{x, 0+} G^{1}(F)_{x,\mathbf{s}_{0}+} \cdots G^{i}(F)_{x, \mathbf{s}_{i-1}+} = K_{+}^{i}, \\
	K_{+}^{d} \cap G(F)_{x, \mathbf{s}_{i}+} & = & G^{i+1}(F)_{x, \mathbf{s}_{i}+} \cdots G^{d}(F)_{x, \mathbf{s}_{d-1}+}.  
\end{eqnarray*}

Using $\hat{\Phibf}_{i}$, Yu constructed a representation $\rho_{j}$ of $K_{j}$ for $j=0, \ldots, d$.  

\begin{lem}[{{\cite[\S 4]{Yu}}}]
\label{Phitilde}
Let $0 \leq i \leq d-1$.  
There is an irreducible representation $\tilde{\Phibf}_{i}$ of $K^{i} \ltimes J^{i+1}$ such that 
\begin{enumerate}
\item $\tilde{\Phibf}_{i} | _{1 \ltimes J_{+}^{i+1}}$ is $\hat{\Phibf}_{i} | _{J_{+}^{i+1}}$-isotypic, and
\item $\tilde{\Phibf}_{i} | _{K_{+}^{i} \ltimes 1}$ is $\mathbf{1}$-isotypic.  
\end{enumerate}
\end{lem}

\begin{lem}[{{\cite[\S 4]{Yu}}}]
Let $0 \leq i \leq d-1$.  
Let $\inf(\Phibf_{i})$ be the inflation of $\Phibf_{i} |_{K^{i}}$ to $K^{i} \ltimes J^{i+1}$, and let $\tilde{\Phibf}_{i}$ be as in Lemma \ref{Phitilde}.  
Then the $K^{i} \ltimes J^{i+1}$-representation $\inf(\Phibf_{i}) \otimes \tilde{\Phibf}_{i}$ factors through $K^{i} \ltimes J^{i+1} \to K^{i}J^{i+1}=K^{i+1}$.  
\end{lem}

\begin{defn}
We denote by $\Phibf'_{i}$ the $K^{i+1}$-representation $\inf(\Phibf_{i}) \otimes \tilde{\Phibf}_{i}$.  
\end{defn}

To obtain $\rho_j$ constructed by Yu, we use a little different way from Yu, by Hakim--Murnaghan.  

\begin{lem}[{{\cite[3.4]{HM}}}]
\begin{enumerate}
\item For $i=1, \ldots, d-1$, we have $K^{i} \cap J^{i+1} = G^{i}(F)_{x, \rbf_{i}} \subset J^{i}$.  
\item For $i=0, \ldots, d-1$, let $\mu$ be a $K^{i}$-representation which is trivial on $K^{i} \cap J^{i+1}$.  
Then we can obtain the inflation $\inf_{K^{i}}^{K^{i+1}} \mu$ of $\mu$ to $K^{i+1}$ via $K^{i+1}/J^{i+1} \cong K^{i}/(K^{i} \cap J^{i+1})$.  
The representation $\inf_{K^{i}}^{K^{i+1}} \mu$ is trivial on $J^{i+1}$, and also trivial on $K^{i+1} \cap J^{i+2}$ if $i<d-1$.  
\item If $i, \mu$ is as in (ii) and $i \leq j \leq d$, then we can also obtain the inflation $\inf_{K^{i}}^{K^{j}} \mu$ of $\mu$ to $K^{j}$ as $\inf_{K^{i}}^{K^{j}} \mu = \inf_{K^{j-1}}^{K^{j}} \circ \cdots \circ \inf_{K^{i}}^{K^{i+1}} \mu$.  
\end{enumerate}
\end{lem}

\begin{defn}[{{\cite[3.4]{HM}}}]
For $0 \leq i < j \leq d$, we put $\kappa_{i}^{j} = \inf_{K^{i+1}}^{K^{j}} \Phibf'_{i}$.  
For $0 \leq j \leq d$, we put $\kappa_{-1}^{j} = \inf_{K^{0}}^{K^{j}} \rho$ and $\kappa_{j}^{j} = \Phibf_{j}|_{K^{j}}$.  
And also, for $-1 \leq i \leq d$ we put $\kappa_{i}=\kappa_{i}^{d}$.  
\end{defn}

\begin{prop}[{{\cite[3.23]{May}}}]
Let $0 \leq j \leq d$.  
The representation $\rho_{j}$ constructed by Yu is isomorphic to
\[
	\kappa_{-1}^{j} \otimes \kappa_{0}^{j} \otimes \cdots \otimes \kappa_{j}^{j}.  
\]
In particular, 
\[
	\rho_{d} \cong \kappa_{-1} \otimes \kappa_{0} \otimes \cdots \otimes \kappa_{d}.  
\]
\end{prop}

Therefore, we obtain $\rho_{j}$ constructed by Yu.  

\begin{thm}[{{\cite[15.1]{Yu}}}]
The compactly induced representation $\cInd_{K^{j}}^{G^{j}(F)} \rho_{j}$ of $G^{j}(F)$ is irreducible and supercuspidal.  
\end{thm}

For later use, we recall the following proposition on the dimension of representation space of $\kappa_{i}$.  

\begin{prop}[{{\cite[3.24]{May}}}]
Let $0 \leq i < j \leq d$.  
Then the dimension of $\kappa_{i}^{j}$ is equal to the dimension of $\Phibf'_{i}$, which is also equal to $(J^{i+1}:J_{+}^{i+1})^{1/2}$.  
\end{prop}

\section{Tame simple strata}
\label{TameSimple}
In this section, we consider the class of simple strata corresponding to some Yu datum.  

\begin{defn}
\begin{enumerate}
\item A pure stratum $[\Afra, n, r, \beta]$ is called tame if $E=F[\beta]$ is a tamely ramified extension of $F$.  
\item A simple type $(J, \lambda)$ associated with a simple stratum $[\Afra, n, 0, \beta]$ is called tame if $[\Afra, n, 0, \beta]$ is tame.  
\end{enumerate}
\end{defn}

\begin{rem}
\label{esstame}
\begin{enumerate}
\item By \cite[(2.6.2)(4)(b), 2.7 Proposition]{BH}, the above definition is independent of the choice of simple strata.  
\item Essentially tame supercuspidal representations, defined in \cite[2.8]{BH}, are $G$-representations containing some tame simple types.  
\end{enumerate}
\end{rem}

As explained in \S \ref{SecherreStevens}, any simple strata has a defining sequence.  
Actually, if a simple stratum $[\Afra, n, 0, \beta]$ is tame, then we can show the existence of a ``nice" defining sequence of $[\Afra, n, 0, \beta]$.  
To take such as a defining sequence, we use several propositions.  

\begin{prop}[{{\cite[3.1 Corollary]{BH}}}]
\label{BHforapp}
Let $E$ be a finite, tamely ramified extension of $F$ and let $\beta \in E$ such that $E=F[\beta]$.  
Let $[\Afra(E), n, r, \beta]$ be a pure stratum in $A(E)$ with $r=-k_{F}(\beta) < n$.  
Then there exists $\gamma \in E$ such that $[\Afra(E), n, r, \gamma]$ is simple and equivalent to $[\Afra(E), n, r, \beta]$.  
Moreover, if $\iota : E \hookrightarrow A$ is an $F$-algebra inclusion and $[\Afra, n', r', \iota(\beta)]$ is a pure stratum of $A$ with $r'=-k_{0}(\iota(\beta), \Afra)$, then $[\Afra, n', r', \iota(\gamma)]$ is simple and equivalent to $[\Afra, n', r', \beta]$.  
\end{prop}

\begin{prop}[{{\cite[Proposition 4.4]{May}}}]
\label{Mforapp}
Assume $A \cong \M_{N}(F)$ for some $N$.  
Let $[\Afra, n, r, \beta]$ be a tame, pure stratum of $A$ with $r=-k_{0} (\beta, \Afra)$.  
Let $\gamma \in E=F[\beta]$ such that $[\Afra, n, r, \gamma]$ is simple and equivalent to $[\Afra, n, r, \beta]$.  
Then $[\Bfra_{\gamma}, r, r-1, \beta-\gamma]$ is simple.  
\end{prop}

By these propositions, we obtain the following proposition needed in our case.  

\begin{prop}
\label{appfortame}
Let $[\Afra, n, r, \beta]$ be a pure stratum of $A$ with $r=-k_{0}(\beta, \Afra)$.  
Then there exists an element $\gamma$ in $F[\beta]$ satisfying the following conditions:  
\begin{enumerate}
\item The stratum $[\Afra, n, r, \gamma]$ is simple and equivalent to $[\Afra, n, r, \beta]$.   
\item $\beta-\gamma$ is minimal over $F[\gamma]$.  
\item The equality $v_{\Afra}(\beta-\gamma)=k_{0}(\beta, \Afra)$ holds.  
\end{enumerate}
\end{prop}

\begin{prf}
By Proposition \ref{BHforapp}, there exists $\gamma$ satisfying (1).  
We show that $\gamma$ also satisfy (2) and (3).  
We apply Proposition \ref{Mforapp} to the case $A=A(E)$.  
Then the stratum $[\Bfra', -k_{0} \left( \beta, \Afra(E) \right), -k_{0} \left( \beta, \Afra(E) \right) -1, \beta-\gamma]$ is simple, where $\Bfra'=\Cent_{A(E)}(\gamma) \cap \Afra(E)$.  
Since this stratum is simple, $\beta-\gamma$ is minimal over $F[\gamma]$ and (2) is satisfied.  
To obtain (3), we calculate $v_{\Afra}(\beta-\gamma)$ and $k_{0}(\beta, \Afra)$.   
First, we have
\[
	v_{E}(\beta-\gamma) = v_{\ofra_{E}}(\beta-\gamma) = v_{\Bfra'}(\beta-\gamma) = -\left( -k_{0} \left( \beta, \Afra(E) \right) \right) = k_{0} \left( \beta, \Afra(E) \right) .  
\]
Then we obtain 
\[
	v_{\Afra}(\beta-\gamma) = \frac{e(\Afra|\ofra_{F})}{e(E/F)}v_{E}(\beta-\gamma) = \frac{e(\Afra|\ofra_{F})}{e(E/F)} k_{0} \left( \beta, \Afra(E) \right) = k_{0}(\beta, \Afra)
\]
and (3) is also satisfied.  
\end{prf}

Conversely, the following proposition is needed to construct a simple stratum in $A$ from some hereditary $\ofra_{F}$-order and elements in $A$ in \S \ref{YutoSS}.  

\begin{prop}[{{\cite[Proposition 4.2]{May}}}]
\label{Mforcons}
Assume $A \cong \M_{N}(F)$ for some $N$.  
Let $[\Afra, n, r, \beta]$ be a tame simple stratum.  
Let $[\Bfra_{\beta}, r, r-1, b]$ be a simple stratum, where $\Bfra_{\beta}=\Afra \cap \Cent_{A}(\beta)$.  
Suppose $b \notin F[\beta]$.  
Then we have $F[\beta + b]=F[\beta, b]$ and $[\Afra, n, r, \beta]$ is a pure stratum with $k_{0}(\beta+b, \Afra)=-r$.  
\end{prop}

\section{Tame twisted Levi subgroups of $G$}
\label{TTLS}

First, we show some subgroups in $G$ are tame twisted Levi subgroups.  

Let $V$ be a right $D$--module.  
Let $E/F$ be a field extension in $\End_{D}(V)$.  
Then $V$ can be equipped with the canonical right $D \otimes _{F} E$--module structure and we can define an $E$--scheme $\underline{\Aut}_{D \otimes _{F} E}(V)$ as
\[
\underline{\Aut}_{D \otimes_{F} E}(V)(C)=\Aut _{D \otimes_{F} C}(V \otimes_{E} C)
\]
for an $E$--algebra $C$.  

Let $E'/E/F$ be a field extension in $\End_{D}(V)$ such that $E'$ is a tamely ramified extension of $F$.  
We put $G=\underline{\Aut}_{D}(V)$, $H=\Res_{E/F} \underline{\Aut}_{D \otimes _{F} E}(V)$ and $H'= \Res_{E'/F} \underline{\Aut}_{D \otimes_{F} E'}(V)$.  
Then $H'$ is a closed subscheme in $H$ and $H$ is a closed subscheme in $G$.  

To show that $(H', H, G)$ is a tame twisted Levi sequence, we fix a maximal torus in $G$.  
We take a maximal subfield $L$ in $\End_{D}(V)$ such that $L$ is a tamely ramified extension of $E'$.  
We put $T=\Res_{L/F} \underline{\Aut}_{D \otimes _{F} L}(V)$.  

We put $I_{1} = \{ 1, \ldots, [E:F] \}$, $I_{2} = \{ 1, \ldots, [E':E] \}$ and $I_{3} = \{ 1, \ldots, [L:E'] \}$.  
Let $\sigma_1, \sigma_2, \ldots, \sigma_{[E:F]}$ be distinct elements in $\Hom_{F}(E, \bar{F})$.  
For $i \in I_{1}$, let $\sigma_{i,1}, \sigma_{i,2}, \ldots, \sigma_{i,[E':E]}$ be distinct elements in $\Hom_{F}(E', \bar{F})$ whose restrictions to $E$ are equal to $\sigma_{i}$.  
For $(i,j) \in I_{1} \times I_{2}$, let $\sigma_{i,j,1}, \sigma_{i,j,2}, \ldots, \sigma_{i,j,[L:E']}$ be distinct elements in $\Hom_{F}(L, \bar{F})$ whose restrictions to $E'$ are equal to $\sigma_{i,j}$.  
Then we have 
\begin{eqnarray*}
\Hom_{F}(E', \bar{F}) &=& \left\{ \sigma_{i,j} | (i,j) \in I_{1} \times I_{2} \right\} \\
\Hom_{F}(L, \bar{F}) &=& \left\{ \sigma_{i,j,k} | (i,j,k) \in I_{1} \times I_{2} \times I_{3} \right\}
\end{eqnarray*}
as $L/F$ is separable.  

Let $\tilde{L}$ be the Galois closure of $L/F$ in $\bar{F}$ and let $C$ be an $\tilde{L}$--algebra.  
Then we have an $L$--algebra isomorphism 
\begin{eqnarray*}
	L \otimes_{F} C &\cong& \prod_{(i,j,k) \in I_{1} \times I_{2} \times I_{3}} C_{i,j,k} \\
	l \otimes a & \mapsto & \left( \sigma_{i,j,k}(l)a \right)_{i,j,k}, 
\end{eqnarray*}
where $C_{i,j,k}=C$ for $(i,j,k) \in I_{1} \times I_{2} \times I_{3}$.  
Similarly, we have isomorphisms $E' \otimes_{F} C \cong \prod_{i,j} C_{i,j}$ and $E \otimes_{F} C \cong \prod_{i} C_{i}$ where $C_{i,j}$ and $C_{i}$ are also isomorphic to $C$.  
In the canonical inclusion $E \otimes C \hookrightarrow E' \otimes C$, the algebra $C_i$ is diagonally embedded in $\prod_{j}C_{i,j}$.  
And also, in the inclusion $E' \otimes C \hookrightarrow L \otimes C$, the algebra $C_{i,j}$ is diagonally embedded in $\prod_{k}C_{i,j,k}$.  
We put $V_{i,j,k}=V \otimes_{L} C_{i,j,k}$.  

\begin{prop}
\label{Levidecom}
\begin{enumerate}
\item Let $C$ be an extension field of $\tilde{L}$.  
Then we have a commutative diagram of $C$--schemes:  

\[
	\xymatrix{
		T \times_{F} C \ar[r]^-{\cong} \ar[d] & \prod_{i,j,k} \underline{\Aut}_{D \otimes_{F} C_{i,j,k} } (V_{i,j,k}) \ar[d] \\
		H' \times_{F} C \ar[r]^-{\cong} \ar[d] & \prod_{i,j} \underline{\Aut}_{D \otimes_{F} C_{i,j} } \left( \bigoplus_{k} V_{i,j,k}  \right) \ar[d] \\
		H \times_{F} C \ar[r]^-{\cong} \ar[d] & \prod_{i} \underline{\Aut}_{D \otimes_{F} C_{i} } \left( \bigoplus_{j,k} V_{i,j,k}  \right) \ar[d] \\
		G \times_{F} C \ar[r]^-{\cong} &  \underline{\Aut}_{D \otimes_{F} C } \left( \bigoplus_{i,j,k} V_{i,j,k}  \right).   \\
	}
\]

\item We have a commutative diagram of $C$-vector spaces:
\[
	\xymatrix{
		\Lie(T \times_{F} C) \ar[r]^-{\cong} \ar[d] & \prod_{i,j,k} \End_{D \otimes_{F} C_{i,j,k} } (V_{i,j,k}) \ar[d] \\
		\Lie(H' \times_{F} C) \ar[r]^-{\cong} \ar[d] & \prod_{i,j} \End_{D \otimes_{F} C_{i,j} } \left( \bigoplus_{k} V_{i,j,k}  \right) \ar[d] \\
		\Lie(H \times_{F} C) \ar[r]^-{\cong} \ar[d] & \prod_{i} \End_{D \otimes_{F} C_{i} } \left( \bigoplus_{j,k} V_{i,j,k}  \right) \ar[d] \\
		\Lie(G \times_{F} C) \ar[r]^-{\cong} & \End_{D \otimes_{F} C } \left( \bigoplus_{i,j,k} V_{i,j,k}  \right),  \\
	}
\]
where the vertical morphisms are all monomorphisms.  

\item Let $c \in L$, and let $m_{c} \in \Lie(T) = \End_{D \otimes L}(V)$ be the map $v \mapsto cv$ for $v \in V$.  
We put $m_{c, C}=m_{c} \otimes_{F} 1 \in \Lie(T) \otimes _{F} C = \Lie(T \times_{F} C)$.  
When we regard $m_{c,C}$ as an element in $ \End_{D \otimes_{F} C } \left( \bigoplus_{i,j,k} V_{i,j,k}  \right)$ via the morphisms in (ii), for $v_{i,j,k} \in V_{i,j,k}$ we have $m_{c,C}(v_{i,j,k})=\sigma_{i,j,k}(c)v_{i,j,k}$.  
\end{enumerate}
\end{prop}

\begin{prf}
This is proved similarly as the proof of \cite[Proposition 6.4]{May}.  
\end{prf}

\begin{cor}
\label{findTTLS}
The sequence $(H', H, G)$ is a tame twisted Levi sequence.  
Moreover, $Z(H')/Z(G)$ is anisotropic.  
\end{cor}

\begin{prf}
We put $C=\tilde{L}$, which is a finite, tamely ramified Galois extension of $F$.  
Since $L$ is a maximal $F$-subfield in $A$, the right $D \otimes _{F} L$-module $V$ is simple.  
Then for any $(i,j,k) \in I_{1} \times I_{2} \times I_{3}$ and $C$-algebra $\tilde{C}$, we have 
\begin{eqnarray*}
	\underline{\End}_{D \otimes _{F} C_{i,j,k}}(V_{i,j,k})(\tilde{C}) & = & \End_{D \otimes _{F} C_{i,j,k} \otimes_{C} \tilde{C}}(V \otimes_{L} C_{i,j,k} \otimes_{C} \tilde{C}) \\
	& \cong & \End_{D \otimes_{F} L \otimes_{L} \tilde{C}}(V \otimes_{L} \tilde{C}) \\
	& \cong & \End_{D \otimes_{F} L} (V) \otimes_{L} \tilde{C} \\
	& \cong & L \otimes_{L} \tilde{C} \cong \tilde{C} = \underline{\End}_{C}(C) (\tilde{C}).  
\end{eqnarray*}
Therefore we have $\underline{\End}_{D \otimes _{F} C_{i,j,k}}(V_{i,j,k}) \cong \underline{\End}_{C}(C)$ as $C$-schemes.  
We also have 
\begin{eqnarray*}
	\prod_{i,j,k} \underline{\End}_{D \otimes_{F} C_{i,j,k}}(V_{i,j,k}) & \cong & \prod_{i,j,k} \underline{\End}_{C}(C), \\
	\prod_{i,j} \underline{\End}_{D \otimes_{F} C_{i,j}} \left( \bigoplus_{k} V_{i,j,k} \right) & \cong & \prod_{i,j} \underline{\End}_{C}\left( C^{\oplus |I_{3}|} \right), \\
	\prod_{i} \underline{\End}_{D \otimes_{F} C_{i}} \left( \bigoplus_{j,k} V_{i,j,k} \right) & \cong & \prod_{i} \underline{\End}_{C} \left( C^{ \oplus (|I_{2}| \times |I_{3}|) } \right), \\
	\underline{\End}_{D \otimes_{F} C} \left( \bigoplus_{i,j,k} V_{i,j,k} \right) & \cong & \underline{\End}_{C} \left( C^{\oplus (|I_{1}| \times |I_{2}| \times |I_{3}|}) \right).  
\end{eqnarray*}
By taking the multiplicative group, we obtain
\[
\begin{array}{rcccl}
	T \times_{F} C & \cong & \prod_{i,j,k} \underline{\Aut}_{D \otimes_{F} C_{i,j,k}}(V_{i,j,k}) & \cong & \mathbb{G}_{m}{}^{\times (|I_{1}| \times |I_{2}| \times |I_{3}|)}, \\
	H' \times_{F} C & \cong & \prod_{i,j} \underline{\Aut}_{D \otimes_{F} C_{i,j}} \left( \bigoplus_{k} V_{i,j,k} \right) & \cong &  \GL_{|I_{3}|}{}^{\times (|I_{1}| \times |I_{2}|)}, \\
	H \times_{F} C & \cong & \prod_{i} \underline{\Aut}_{D \otimes_{F} C_{i}} \left( \bigoplus_{j,k} V_{i,j,k} \right) & \cong &  \GL_{|I_{2}| \times |I_{3}|}{}^{\times |I_{1}| }, \\
	G \times_{F} C & \cong & \underline{\Aut}_{D \otimes{F} C} \left( \bigoplus_{i,j,k} V_{i,j,k} \right) & \cong & \GL_{|I_{1}| \times |I_{2}| \times |I_{3}|}.  
\end{array}
\]
Therefore $H' \times_{F} C$ and $H \times_{F} C$ are Levi subgroups in $G \times _{F} C$ with a split maximal torus $T \times_{F} C$.  
Since $C$ is a finite, tamely ramified Galois extension of $F$, the sequence $(H', H, G)$ is a tame twisted Levi sequence.  

Moreover, we have $\left( Z(H')/Z(G) \right) (F)=E'^{\times}/F^{\times}$, which is compact.  
Then $Z(H')/Z(G)$ is anisotropic.  
\end{prf}

Let $C=\bar{F}$.  
For distinct elements $(i',j',k'), (i'', j'', k'') \in I_{1} \times I_{2} \times I_{3}$, we define the root $\alpha_{(i',j',k'),(i'',j'',k'')} \in \Phi(G, T; \bar{F})$ as
\[
	\alpha_{(i',j',k'),(i'',j'',k'')} : \prod_{i,j,k} \Aut_{D \otimes_{F} \bar{F}_{i,j,k} } (V_{i,j,k}) \to \bar{F}^{\times} ; (t_{i,j,k})_{i,j,k} \mapsto t_{i',j',k'} t_{i'',j'',k''}^{-1}.  
\]

Therefore we have
\begin{eqnarray*}
	\Phi(H, T; \bar{F}) & = & \left\{ \alpha_{(i',j',k'),(i'',j'',k'')} \in \Phi(G, T; \bar{F}) | i'=i'' \right\} \\
	\Phi(H', T; \bar{F}) & = & \left\{ \alpha_{(i',j',k'),(i'',j'',k'')} \in \Phi(G, T; \bar{F}) | i'=i'', j'=j'' \right\}, 
\end{eqnarray*}
and we obtain
\[
	\Phi(H, T; \bar{F}) \setminus \Phi(H', T; \bar{F}) = \left\{ \alpha_{(i',j',k'),(i'',j'',k'')} \in \Phi(G, T; \bar{F}) | i'=i'', j' \neq j'' \right\}.  
\]

Moreover, the coroot $\check{\alpha}_{(i',j',k'),(i'',j'',k'')}$ with respect to $\alpha_{(i',j',k'),(i'',j'',k'')}$ is as follows:  
\[
	\check{\alpha}_{(i',j',k'),(i'',j'',k'')} : \bar{F}^{\times} \to \prod_{i,j,k} \Aut_{D \otimes_{F} \bar{F}_{i,j,k} } (V_{i,j,k}) \cong \prod_{i,j,k} \bar{F}^{\times} ; t \mapsto (t_{i,j,k})_{i,j,k}, 
\]
where 
\[
	t_{i,j,k} = \begin{cases}
		t & \left( (i,j,k) = (i',j',k') \right), \\
		t^{-1} & \left( (i,j,k) = (i'',j'',k'') \right), \\
		1 & otherwise.  
		\end{cases}
\]

Then we have $\dr \check{\alpha}_{(i',j',k'),(i'',j'',k'')}(u)=(u_{i,j,k})_{i,j,k}$ where
\[
	u_{i,j,k} = \begin{cases}
		u & \left( (i,j,k) = (i',j',k') \right), \\
		-u & \left( (i,j,k) = (i'',j'',k'') \right), \\
		0 & otherwise.  
		\end{cases}
\]

Conversely, we determine the set of tame twisted Levi subgroup $G'$ in $G$ with $Z(G')/Z(G)$ anisotropic.  

\begin{lem}
\label{detofTLG}
Let $G'$ be a tame twisted Levi subgroup of $G=\underline{\Aut}_{D}(V)$.  
Suppose $Z(G')/Z(G)$ is anisotropic.  
Then there exists a finite, tamely ramified extension $E$ of $F$ such that $G' \cong \Res_{E/F} \underline{\Aut}_{D \otimes _{F} E} (V)$.  
\end{lem}

\begin{prf}
Let $F^{\mathrm{tr}}$ be the maximal tamely ramified extension of $F$.  
Since $G'$ is a tame twisted Levi subgroup in $G$, $G'_{F^{\mathrm{tr}}}$ is a Levi subgroup in $G_{F^{\mathrm{tr}}} \cong \underline{\Aut}_{D \otimes F^{\mathrm{tr}}}(V \otimes F^{\mathrm{tr}})$.  
There exists a one-to-one relationship between Levi subgroups in $G_{F^{\mathrm{tr}}}$ and direct decompositions of $V \otimes F^{\mathrm{tr}}$ as a right $D \otimes F^{\mathrm{tr}}$-module.  
Let $V \otimes F^{\mathrm{tr}} = \bigoplus_{k=1}^{j} V_{k}$ be the corresponding decomposition with $G'_{F^{\mathrm{tr}}}$.  
Then we have $G'_{F^{\mathrm{tr}}} = \prod_{k=1}^{j} \underline{\Aut}_{D \otimes F^{\mathrm{tr}}}(V_{k})$.  
We remark that the right-hand-side group is the multiplicative group of $\underline{\End}_{D \otimes F^{\mathrm{tr}}}(V_{k})$ with a $Gal(F^{\mathrm{tr}}/F)$-action defined by its $F$-structure.  
Let $Z_{k}$ be the center of $\underline{\End}_{D \otimes F^{\mathrm{tr}}}(V_{k})$, which is $F^{\mathrm{tr}}$-isomorphic to $\underline{\End}_{F^{\mathrm{tr}}}(F^{\mathrm{tr}})$.  
Then $Z(G')_{F^{\mathrm{tr}}}$ is the multiplicative group of $Z=\prod_{k=1}^{j} Z_{k}$, equipped with the same $Gal(F^{\mathrm{tr}}/F)$-action.  
Therefore, we consider the structure of $Z_{k}$.  
Let $\mathbf{1}_{k}$ be (the $F^{\mathrm{tr}}$-rational point corresponding to) the identity element in $Z_{k}$.  
Since the $\Gal(F^{\mathrm{tr}}/F)$-action to $Z$ preserves the $F$-algebra structure, the set $\{ \mathbf{1}_{k} \mid k=1, \ldots, j \}$ is $\Gal(F^{\mathrm{tr}}/F)$-invariant.  
Then by changing the indices if necessary, we may assume there exist integers $0 = n_{0} < n_{1} < \cdots < n_{l} = j$ such that $\Gal(F^{\mathrm{tr}}/F)$ acts the set $\{ \mathbf{1}_{n_{i-1}+1}, \ldots, \mathbf{1}_{n_{i}} \}$ transitively for $l=1, \ldots, i$.  
We put $Y_{i}=\prod_{k=n_{i-1}+1}^{n_{i}} Z_{k}$.  
Since $a \in F^{\mathrm{tr}}, b \in Z$ and $\gamma \in \Gal(F^{\mathrm{tr}}/F)$ we have $\gamma (ab) = \gamma(a) \gamma(b)$ and $\{ \mathbf{1}_{n_{i-1}+1}, \ldots, \mathbf{1}_{n_{i}} \}$ is $\Gal(F^{\mathrm{tr}}/F)$-invariant, $Y_{i}$  is also $\Gal(F^{\mathrm{tr}}/F)$-invariant.  
Then $Y_{i}$ is defined over $F$.  
Let $X_{i}$ be the Galois descent of $Y_{i}$ to $F$.  
Let $\Gal(F^{\mathrm{tr}}/F_{i})$ be the stabilizer of $\mathbf{1}_{n_{i}}$. 
The fields $F_{i}$ is tamely ramified, and finite-dimensional over $F$ since $\Gal(F^{\mathrm{tr}}/F)/\Gal(F^{\mathrm{tr}}/F_{i})$ is $\Gal(F^{\mathrm{tr}}/F)$-isomorphic to the finite set $\{ \mathbf{1}_{n_{i-1}+1}, \ldots, \mathbf{1}_{n_{i}} \}$.  

We show $X_{i}$ is isomorphic to $\Res_{F_{i}/F} \underline{\End}_{F_{i}}(F_{i})$.  
If this follows, then we can show the multiplicative group of $X_{i}$ is $\Res_{F_{i}/F} \mathbb{G}_{m}$ and $Z(G') = \prod_{i=1}^{l} \Res_{F_{i}/F} \mathbb{G}_{m}$.  

Any $F^{\mathrm{tr}}$-rational point of $Y_{i}$ is uniquely represented as the form $\sum_{k'=n_{i-1}+1}^{n_{i}} a_{k'}\mathbf{1}_{k'}$, where $a_{k'} \in F^{\mathrm{tr}}$.  
Suppose $z=\sum_{k'=n_{i-1}+1}^{n_{i}} a_{k'}\mathbf{1}_{k'}$ is stabilized by $\Gal(F^{\mathrm{tr}}/F)$.  
For any $\gamma \in \Gal(F^{\mathrm{tr}}/F_{i})$, we have $z=\gamma(z)=\sum_{k'=n_{i-1}+1}^{n_{i}-1} \gamma(a_{k'})\gamma(\mathbf{1}_{k'}) + \gamma(a_{n_{i}}) \mathbf{1}_{n_{i}}$.  
Then we have $\gamma(a_{n_{i}}) = a_{n_{i}}$, that is, $a_{n_{i}} \in F_{i}$.  
For $n_{i-1}<k'<n_{i}$, we pick $\gamma_{k'} \in \Gal(F^{\mathrm{tr}}/F)$ such that $\gamma_{k'}(\mathbf{1}_{n_{i}})=\mathbf{1}_{k'}$.  
Then we have
\[
z = \gamma_{k'}(z) = \sum_{k''=n_{i-1}+1}^{n_{i}-1} \gamma_{k''}(a_{k''}) \gamma_{k'}(\mathbf{1}_{k'}) + \gamma_{k'}(a_{n_i}) \mathbf{1}_{k'}, 
\]
whence $a_{k'} = \gamma_{k'}(a_{n_i})$.  
Therefore any $F$-rational point of $X_{i}$ is the form
\[
\sum_{k'=n_{i-1}+1}^{n_{i}-1} \gamma_{k'}(a_{n_{i}}) \mathbf{1}_{k'} + a_{n_{i}} \mathbf{1}_{n_{i}}, 
\]
where $a_{n_{i}} \in F_{i}$, and the ring structure of $X_{i}(F)$ is isomorphic to $F_{i}$.  
Since the ring structure of $X_{i}(C)$ is isomorphic to $X_{i}(F) \otimes C$ for any $F$-algebra $C$, we obtain $X_{i} \cong \Res_{F_{i}/F} \underline{\End}_{F_{i}}(F_{i})$.  

We have shown $Z(G') = \prod_{i=1}^{l} \Res_{F_{i}/F} \mathbb{G}_{m}$.  
Since $Z(G')/Z(G)$ is anisotropic and $Z(G)=\mathbb{G}_{m}$, we have $l=1$ and $Z(G') = \Res_{E/F} \mathbb{G}_{m}$, where we put $E=F_{1}$.  

The field $E$ can be regarded as a $F$-subfield in $A$ via $X \subset \underline{\End}_{D}(V)$.  
We put $H=\underline{\Aut}_{D \otimes E}(V)$.  
Then $H$ is a tame twisted Levi subgroup in $G$ and we have $Z(H)=Z(G')$.  
Since there exists a one-to-one relationship between subtori in $G$ defined over $F$ and Levi subgroups in $G$ defined over $F$, we obtain $G'=H$.  
\end{prf}

\section{Embeddings of buildings for Levi sequences of $G$}
\label{BroussousLemaire}

\subsection{Lattice functions in $V$}

First, we recall the lattice functions in $V$ and their properties from \cite{BL}.  

\begin{defn}
The map $\Lcal$ from $\R$ to the set of $\ofra_{D}$-lattices in $V$ is a lattice function in $V$ if
\begin{enumerate}
\item we have $\Lcal(r)\varpi_{D}=\Lcal \left( r+(1/d) \right)$ for some uniformizer $\varpi_{D}$ of $D$ and $r \in \R$, 
\item $\Lcal$ is decreasing, that is, $\Lcal(r) \supset \Lcal(r')$ if $r \leq r'$, and
\item $\Lcal$ is left-continuous, where the set of lattices in $V$ is equipped with the discrete topology.  
\end{enumerate}
\end{defn}

The set of lattice functions in $V$ is denoted by $\Latt^{1}(V)$.  
The groups $G$ and $\R$ act on $\Latt^{1}(V)$ by $(g \cdot \Lcal)(r) = g \cdot \left( \Lcal(r) \right)$ and $(r' \cdot \Lcal) (r)=\Lcal(r+r')$ for $g \in G$, $r, r' \in \R$ and $\Lcal \in \Latt^{1}(V)$.  
These actions are compatible, and then $\Latt(V) := \Latt^{1}(V)/\R$ is equipped with the canonical $G$-action.  
The $G$-sets $\Latt^{1}(V)$ and $\Latt(V)$ are also equipped with an affine structure.  
Then there exists a canonical $G$-equivariant, affine isomorphism $\Bscr^{E}(G,F) \to \Latt^{1}(V)$.  
This isomorphism induces a $G$-equivariant, affine isomorphism $\Bscr^{R}(G,F) \to \Latt(V)$.  

We construct lattice functions from $\ofra_{D}$-sequences.  
Let $c \in \R$ and let $(\Lcal_{i})_{i \in \Z}$ be an $\ofra_{D}$-sequence with period $e$.  
Then 
\[
	\Lcal(r) = \Lcal_{\lceil de(r - c) \rceil}, \, r \in \R
\]
is a lattice function in $V$.  

\begin{prop}
\label{lfconstbylc}
Let $\Lcal$ be a lattice function in $V$.  
The following assertions are equivalent:  
\begin{enumerate}
\item $\Lcal$ is constructed from an $\ofra_{D}$-chain.  
\item There exists $c \in \R$ and $e \in \Z_{>0}$ such that the set of discontinuous points of $\Lcal$ is equal to $c+(de)^{-1}\Z$.  
\end{enumerate}
Moreover, if (1) (and (2)) holds, $e$ is equal to the period of some $\ofra_{D}$-chain which $\Lcal$ constructed from.  
\end{prop}

\begin{prf}
First, suppose $\Lcal$ is constructed from an $\ofra_{D}$-chain.  
Then there exists $c \in \R$ and an $\ofra_{D}$-chain $(\Lcal_{i})_{i \in \Z}$ with period $e$ such that $\Lcal(r) = \Lcal_{\lceil de(r-c) \rceil}$ for $r \in \R$.  
Since $(\Lcal_{i})$ is an $\ofra_{D}$-chain, the set of discontinuous points of $\Lcal$ is equal to $c+(de)^{-1}\Z$, whence (2) holds.  

Conversely, suppose (2) holds.  
For $i \in \Z$, we put $\Lcal_{i}=\Lcal(c+(de)^{-1}i)$.  
Since $\Lcal$ is not right-continuous at $r=c+(de)^{-1}(i+1)$, we have 
\[
	\Lcal_{i} = \Lcal(c+(de)^{-1}i) \supsetneq \Lcal(c+(de)^{-1}(i+1)) = \Lcal_{i+1}.  
\]
Moreover, we also have
\[
	\Lcal_{i+e} = \Lcal(c+(de)^{-1}(i+e)) = \Lcal(c+(de)^{-1}i+d^{-1}) = \Lcal(c+(de)^{-1}i)\varpi_{D} = \Lcal_{i}\varpi_{D}.  
\]
Then $(\Lcal_{i})_{i \in \Z}$ is an $\ofra_{D}$-chain with period $e$.  

Let $\Lcal'$ be the lattice function constructed from $c \in \R$ and the $\ofra_{D}$-chain $(\Lcal_{i})$.  
We show $\Lcal = \Lcal'$.  
For $i \in \Z$, we have $\Lcal'(c+(de)^{-1}i) = \Lcal_{i} = \Lcal(c+(de)^{-1}i)$ and $\Lcal = \Lcal'$ on $c+(de)^{-1} \Z$.  
For $r \in \R$, there exists $i \in \Z$ such that $r \in ( c+(de)^{-1}(i-1), c+(de)^{-1}i ]$.  
Since the set of discontinuous points of $\Lcal$ is $c+(de)^{-1}\Z$, then $\Lcal|_{( c+(de)^{-1}(i-1), c+(de)^{-1}i ]}$ is continuous and 
\[
	\Lcal(r) = \Lcal(c+(de)^{-1}i) = \Lcal_{i} = \Lcal_{\lceil de \left( r-c \right) \rceil} = \Lcal'(r).  
\]
Therefore $\Lcal=\Lcal'$ is the lattice function constructed from the $\ofra_{D}$-chain $(\Lcal_{i})$ of period $e$.  
The last assertion follows from the above argument.  
\end{prf}

Conversely, for any lattice function $\Lcal$ there exists an $\ofra_{D}$-chain $(\Lcal_{i})_{i \in \Z}$ such that $\{ \Lcal(r) \mid r \in \R \} = \{ \Lcal_{i} \mid i \in \Z \}$, unique up to translation.  
Since $\Lcal(r+(1/d)) = \Lcal(r) \varpi_{D}$ for $r \in \R$, the period of $(\Lcal_{i})$ is equal to the number of discontinuous points of $\Lcal$ in $[0, 1/d)$.  

\subsection{Comparison of filtrations:  hereditary orders and Moy--Prasad filtration}

Let $x$ be an element in $\Bscr^{E}(G, F)$, corresponding to a lattice function $\Lcal$ via $\Bscr^{E}(G, F) \cong \Latt^{1}(V)$.  
We can define a filtration $(\mathfrak{a}_{x,r})_{r \in \R}$ in $A$ associated with $x$ as 
\[
	\mathfrak{a}_{x,r} = \mathfrak{a}_{\Lcal,r} = \{ a \in A \mid a\Lcal(r') \subset \Lcal(r+r'), \, r' \in \R \}
\]
for $r \in \R$.  
We also put $\mathfrak{a}_{x,r+} = \bigcup _{r<r'} \mathfrak{a}_{x,r'}$.  
Then we can define a hereditary $\ofra_{F}$-order $\Afra = \mathfrak{a}_{x,0}$ associated with $x$.  
The radical of $\Afra$ is equal to $\Pfra = \mathfrak{a}_{x,0+}$.  
We also put $\U_{0}(x) = \Afra^{\times}$, and  $\U_{r}(x) = 1 + \mathfrak{a}_{x,r}$ for $r \in \R_{> 0}$.  

\begin{prop}[{{\cite[Appendix A]{BL}}}]
\label{compoffiltA}
Let $x \in \Bscr^{E}(G,F)$.  
\begin{enumerate}
\item When we identify $A$ with the Lie algebra $\gfra(F)$ of $G$, we have $\mathfrak{a}_{x,r} = \gfra(F)_{x,r}$ for $r \in \R$.  
\item For $r \geq 0$, we have $\U_{r}(x) = G(F)_{x,r}$.  
\end{enumerate}
\end{prop}

Suppose $\Lcal$ is constructed from an $\ofra_{D}$-chain.  
Then there exist $c \in \R$ and an $\ofra_{D}$-chain $(\Lcal_{i})_{i \in \Z}$ with period $e$ such that $\Lcal(r) = \Lcal_{ \lceil de(r-c) \rceil}$.  
Since $\Lcal_{i+e} = \Lcal_{i}\varpi_{D}$ for $i \in \Z$, we have $\Lcal_{i+de} = \Lcal_{i}\varpi_{F}$, and then $e(\Afra|\ofra_{F}) = de$.  

\begin{prop}
\label{compoffiltG}
Let $x, \Lcal$ be as above, and let $r \in \R$.  
\begin{enumerate}
\item We have $\Pfra^{\lceil r \rceil} = \gfra(F)_{x, r/e(\Afra|\ofra_{F})}$.  
\item Suppose $r \geq 0$.  Then $\U^{\lceil r \rceil}(\Afra) = G(F)_{x,r/e(\Afra|\ofra_{F})}$.  
\item We have $\Kfra(\Afra) = G(F)_{[x]}$.  
\end{enumerate}
\end{prop}

\begin{prf}
We show (1).  
By Proposition \ref{compoffiltA} (1), it suffices to show $\Pfra^{\lceil r \rceil} = \mathfrak{a}_{\Lcal, r/e(\Afra|\ofra_{F})}$.  
We put $n= \lceil r \rceil$.  
Suppose $a \in \mathfrak{a}_{x, r/e(\Afra|\ofra_{F})}$.  
For $n' \in \Z$, we put $r'_{n'} = c + e(\Afra|\ofra_{F})^{-1}n'$.  
Then we have $\Lcal(r'_{n'}) = \Lcal_{\lceil de(r'_{n'} -c) \rceil} = \Lcal_{n'}$, and $\Lcal\left( e(\Afra|\ofra_{F})^{-1}r+r_{n'} \right) = \Lcal_{\lceil de\left( e(\Afra|\ofra_{F})^{-1}r+r'_{n'} -c \right) \rceil} = \Lcal_{n'+\lceil n \rceil}$.  
Since $a \in \mathfrak{a}_{n, r/e(\Afra|\ofra_{F})}$, in particular
\[
	a\Lcal_{n'} = a\Lcal(r') \subset \Lcal(e(\Afra)|\ofra_{F})^{-1}r+r') = \Lcal_{n+n'}
\]
for $n' \in \Z$.  
Since $\{ a \in A \mid a\Lcal_{n'} \subset \Lcal_{n+n'}, \, n' \in \Z \} = \Pfra^{n}$, we have $a \in \Pfra^{n}$.  

Conversely, suppose $a \in \Pfra^{n}$.  
For $r' \in \R$, we have $\Lcal(r') = \Lcal_{\lceil de(r'-c) \rceil}$ and $\Lcal(e(\Afra|\ofra_{F})^{-1}r+r') = \Lcal _{\lceil r+de(r'-c) \rceil}$.  
Since $\lceil r+de(r'-c) \rceil < r+de(r'-c)+1$ and $\lceil de(r'-c) \rceil \geq de(r'-c)$, we have 
\[
	\lceil r+de(r'-c) \rceil - \lceil de(r'-c) \rceil <r+de(r'-c)+1 - de(r'-c) = r+1.  
\]
Since $\lceil r+de(r'-c) \rceil - \lceil de(r'-c) \rceil \in \Z$, we also have $\lceil r+de(r'-c) \rceil - \lceil de(r'-c) \rceil \leq \lceil r \rceil$.  
When we put $n' = \lceil de(r'-c) \rceil$, we have $n + n' \geq \lceil r+de(r'-c) \rceil$.  
Therefore, 
\[
	a \Lcal(r') = a \Lcal_{ \lceil de(r'-c) \rceil} = a \Lcal_{n'} \subset \Lcal_{n+n'} \subset \Lcal_{\lceil r+de(r'-c) \rceil} = \Lcal(e(\Afra|\ofra_{F})^{-1}r+r') 
\]
for $r' \in \R$, which implies $a \in \mathfrak{a}_{\Lcal, r/e(\Afra|\ofra_{F})}$.  
Thus (1) holds.  

To show (2), it is enough to show $\U^{\lceil r \rceil}(\Afra) = \U_{r/e(\Afra|\ofra_{F})}(x)$ by Proposition \ref{compoffiltA} (2).  
Therefore (2) follows from (1).  

(3) is a corollary of \cite[I Lemma 7.3]{BL}, as $\Lcal$ is constructed from an $\ofra_{D}$-chain.  
\end{prf}

\begin{prop}
\label{compoffiltS}
Let $x \in \Bscr^{E}(G, F)$ correspond with a lattice function constructed from an $\ofra_{D}$-chain, and let $n \in \Z$.  
\begin{enumerate}
\item \begin{enumerate}
	\item $\Pfra^{n} = \gfra(F)_{x, n/e(\Afra|\ofra_{F})}$, 
	\item $\Pfra^{n+1} = \gfra(F)_{x, n/e(\Afra|\ofra_{F})+}$, 
	\item $\Pfra^{\lfloor (n+1)/2 \rfloor} = \gfra(F)_{x, n/2e(\Afra|\ofra_{F})}$, 
	\item $\Pfra^{\lfloor n/2 \rfloor +1} = \gfra(F)_{x, n/2e(\Afra|\ofra_{F})+}$.  
	\end{enumerate}
\item Suppose $n \geq 0$.  
	Then we have 
	\begin{enumerate}
	\item $\U^{n}(\Afra) = G(F)_{x, n/e(\Afra|\ofra_{F})}$, 
	\item $\U^{n+1}(\Afra) = G(F)_{x, n/e(\Afra|\ofra_{F})+}$, 
	\item $\U^{\lfloor (n+1)/2 \rfloor}(\Afra) = G(F)_{x, n/2e(\Afra|\ofra_{F})}$, 
	\item $\U^{\lfloor n/2 \rfloor +1}(\Afra) = G(F)_{x, n/2e(\Afra|\ofra_{F})+}$.  
	\end{enumerate}
\end{enumerate}
\end{prop}

\begin{prf}
We show (1), and (2) can be shown in the same way as (1).  

(a) follows from Proposition \ref{compoffiltG} (1).  
(c) also follows from Proposition \ref{compoffiltG} (1) and the fact $\lfloor (n+1)/2 \rfloor = \lceil n/2 \rceil$ for $n \in \Z$.  

We show (b).  
For $r \in (n, n+1]$, we have $\lceil r \rceil = n+1$.  
Then we have 
\begin{eqnarray*}
	\gfra(F)_{x,n/e(\Afra|\ofra_{F})+} & = & \bigcup_{n/e(\Afra|\ofra_{F})<r'} \gfra(F)_{x,r'} \\
	& = & \bigcup_{n/e(\Afra|\ofra_{F})<r' \leq (n+1)/e(\Afra|\ofra_{F})} \Pfra^{\lceil r'e(\Afra|\ofra_{F}) \rceil} \\
	& = & \Pfra^{n+1}.  
\end{eqnarray*}

To show (d), we consider two cases.  
First, suppose $n \in 2\Z$.  
Then we have $\Pfra^{\lfloor n/2 \rfloor + 1} = \Pfra^{(n/2)+1} = \gfra(F)_{x, ((n/2)+1)/e(\Afra|\ofra_{F})}$ by (a).  
Since $n/2 \in \Z$, for any $r \in ( n/2, (n/2)+1 ]$ we have $\lceil r \rceil = (n/2)+1$ and $\gfra(F)_{x,r/e(\Afra|\ofra_{F})} = \Pfra{\lceil r \rceil} = \Pfra^{(n/2)+1}$.  
Therefore 
\begin{eqnarray*}
	\gfra(F)_{x,n/2e(\Afra|\ofra_{F})+} & = & \bigcup_{n/2e(\Afra|\ofra_{F}) < r'}\gfra(F)_{x, r'} \\
	& = & \bigcup_{n/2e(\Afra|\ofra_{F}) < r' \leq ((n/2)+1)/e(\Afra|\ofra_{F})}\Pfra^{\lceil r'e(\Afra|\ofra_{F}) \rceil} \\
	& = & \Pfra^{(n/2)+1} = \Pfra^{\lfloor n/2 \rfloor + 1}.  
\end{eqnarray*}

Next, suppose $n \in \Z \setminus 2\Z$.  
Then we have $\lfloor n/2 \rfloor + 1 = (n+1)/2 = \lceil (n+1)/2 \rceil$ and $\Pfra^{\lfloor n/2 \rfloor + 1} = \gfra(F)_{x,n/2e(\Afra|\ofra_{F})}$ by (b).  
Since $\lceil r \rceil = (n+1)/2$ for $r \in (n/2, (n+1)/2]$, we obtain
\begin{eqnarray*}
	\gfra(F)_{x, n/2e(\Afra|\ofra_{F})+} & = & \bigcup_{n/2e(\Afra|\ofra_{F}) < r'}\gfra(F)_{x, r'} \\
	& = & \bigcup_{n/2e(\Afra|\ofra_{F}) < r' \leq (n+1)/2e(\Afra|\ofra_{F})}\Pfra^{\lceil r'e(\Afra|\ofra_{F}) \rceil} \\
	& = & \Pfra^{(n+1)/2} = \Pfra^{\lceil n/2 \rceil +1}.  
\end{eqnarray*}
\end{prf}

Let $(H', H, G)$ be a tame twisted Levi sequence.  
Then there exists a tower $E'/E/F$ of tamely ramified extensions in $A$ such that $H' = \Res_{E/F}\underline{\Aut}_{D \otimes_{F} E'}(V)$ and $H = \Res_{E'/F}\underline{\Aut}_{D \otimes_{F} E}(V)$.  
We put $B=\Cent_{A}(E)$ and $B'=\Cent_{A}(E')$.  
There exist a division $E$-algebra $D_{E}$ and a right $D_{E}$-module $W$ such that $B \cong \End_{D_{E}}(W)$.  
Similarly, there exist a division $E'$-algebra $D_{E'}$ and a right $D_{E'}$-module $W'$ such that $B' \cong \End_{D_{E'}}(W')$.  
Since $E'/E/F$ is a tower of tamely ramified extensions, we have canonical identifications 
\begin{eqnarray*}
	\Bscr^{E}(H, F) \cong & \Bscr^{E}( \underline{\Aut}_{D \otimes E}(V), E) & \cong \Bscr^{E}(\underline{\Aut}_{D_{E}}(W), E), \\
	\Bscr^{E}(H', F) \cong & \Bscr^{E}( \underline{\Aut}_{D \otimes E'}(V), E') & \cong \Bscr^{E}(\underline{\Aut}_{D_{E'}}(W'), E').  
\end{eqnarray*}

Let $x \in \Bscr^{E}(H', F) \cong \Bscr^{E}(\underline{\Aut}_{D_{E'}}(W'), E')$, and let $\Lcal$ be the corresponding lattice function in $W'$ with $x$.  

\begin{prop}
\label{beingvertex}
The following assertions are equivalent.  
\begin{enumerate}
\item $[x]$ is a vertex in $\Bscr^{R}(H', F)$.  
\item The hereditary $\ofra_{E'}$-order $\Bfra'$ associated with $x$ is maximal.  
\item $\Lcal$ is constructed from an $\ofra_{D_{E'}}$-chain of period 1.  
\end{enumerate}
\end{prop}

\begin{prf}
The element $[x]$ is a vertex if and only if the stabilizer $\Stab_{H'(F)}(x)$ of $x$ in $H'(F)$ is a maximal compact subgroup in $H'(F)$.  
Since $\Lcal$ is identified with $x$ via the $H'(F)$-isomorphism $\Latt^{1}(W') \cong \Bscr^{E}(H', F)$, we have $\Stab_{H'(F)}(x) = \Stab_{\Aut_{D_{E'}}(W')}(\Lcal) = \U(\Bfra')$.  
The group $\U(\Bfra')$ is a maximal compact subgroup in $H'(F)$ if and only if $\Bfra'$ is maximal, which implies the equivalence of (1) and (2).  

To show the equivalence of (2) and (3), let $(\Lcal_{i})$ be an $\ofra_{D_{E'}}$-chain in $W'$ such that $\{ \Lcal(r) \mid r \in \R \} = \{ \Lcal_{i} \mid i \in \Z \}$.  
Since $\Bfra' = \{ b' \in B' \mid b' \Lcal(r) \subset \Lcal(r), \, r \in \R \} = \{ b' \in B' \mid b' \Lcal_{i} \subset \Lcal_{i}, \, i \in \Z \}$, the hereditary $\ofra_{E'}$-order $\Bfra'$ is maximal if and only if the period of $(\Lcal_{i})$ is equal to 1.  
Since the period of $(\Lcal_{i})$ is also equal to the number of discontinuous points of $\Lcal$ in $[0, 1/d_{E'})$, where $d_{E'} = (\dim_{E'}D_{E'})^{1/2}$, (2) holds if and only if there exists a unique discontinuous point $c$ in $[0,1/d_{E'})$.  
Here, since $\Lcal(r+(1/d_{E'}))=\Lcal(r)\varpi_{D_{E'}}$, $\Lcal$ is discontinuous at $c \in \R$ if and only if $\Lcal$ is discontinuous at the unique element $c'$ in $(c+d_{E'}^{-1} \Z) \cap [0, 1/d_{E'})$.  
Therefore (2) holds if and only if the discontinuous points of $\Lcal$ is equal to $c+d_{E'}^{-1} \Z$ for some $c \in \R$, which is also equivalent to (3) by Proposition \ref{lfconstbylc}.  
\end{prf}

We fix an $H'(F)$-equivalent, affine embedding $\iota_{H/H'}:\Bscr^{E}(H', F) \hookrightarrow \Bscr^{E}(H, F)$ and an $H(F)$-equivalent, affine embedding $\iota_{G/H}:\Bscr^{E}(H, F) \hookrightarrow \Bscr^{E}(G, F)$.  
We also put $\iota_{G/H'} = \iota_{G/H} \circ \iota_{H/H'}$.  

\begin{prop}
\label{resultofBL}
Let $x \in \Bscr^{E}(H, F)$.  
\begin{enumerate}
\item The canonical identification $\Bscr^{E}(H, F) \cong \Bscr^{E}(\underline{\Aut}_{D_{E}}(W), E)$ and $\iota_{G/H}$ induce 
\[
	j:\Bscr^{R}(\underline{\Aut}_{D_{E}}(W), E) \hookrightarrow \Bscr^{R}(G, F), 
\]
which is equal to $j_{E}^{-1}$ in \cite[II-Theorem 1.1]{BL}.  
\item Let $(\mathfrak{a}_{\iota_{G/H}(x),r})_{r \in \R}$ be the filtration in $A$ associated with $\iota_{G/H}(x)$, and let $(\mathfrak{b}_{x,r})_{r \in \R}$ be the filtration in $B$ associated with $x$.  
Then 
\[
	\mathfrak{b}_{x, r} = B \cap \mathfrak{a}_{\iota_{G/H}(x), r/e(E/F)}.  
\]
\item The hereditary $\ofra_{F}$-order $\mathfrak{a}_{\iota_{G/H}(x), 0}$ is $E$-pure.  
\end{enumerate}
\end{prop}

\begin{prf}
Since $\Bscr^{E}(H, F) \cong \Bscr^{E}(\underline{\Aut}_{D_{E}}(W), E)$ and $\iota_{G/H}$ are $H(F)$-equivalent and affine, they induce the $H(F)$-equivalent, affine embedding
\[
	j:\Bscr^{R}(\underline{\Aut}_{D_{E}}(W), E) \cong \Bscr^{R}(H', F) \hookrightarrow \Bscr^{R}(G, F).  
\]
However, $H(F)$-equivalent, affine embedding $\Bscr^{R}(\underline{\Aut}_{D_{E}}(W), E) \hookrightarrow \Bscr^{R}(G, F)$ is unique.  
Since $j$ and $j_{E}^{-1}$ are $H(F)$-equivalent and affine, we obtain $j = j_{E}^{-1}$.  
The remainder assertions are results from \cite[II-Theorem 1.1]{BL}.  
\end{prf}

\begin{prop}
\label{beingprincipal}
Let $x \in \Bscr^{E}(H', F)$ such that $[x]$ is a vertex.  
\begin{enumerate}
\item The corresponding lattice function $\Lcal$ in $W$ with $\iota_{H/H'}(x)$ is constructed from a uniform $\ofra_{D_{E}}$-chain.  
In particular, the hereditary $\ofra_{E}$-order $\Bfra$ in $B$ associated with $\Lcal$ is principal.  
\item Let $\Bfra'$ be the hereditary $\ofra_{E'}$-order in $B'$ associated with $x$.  
Then $\Bfra$ is the unique $E'$-pure hereditary $\ofra_{E}$-order in $B$ such that $\Bfra' = B' \cap \Bfra$.  
\end{enumerate}
\end{prop}

\begin{prf}
By Proposition \ref{beingvertex}, the corresponding lattice function in $W'$ with $x$ is constructed from an $\ofra_{D_{E'}}$-chain with period 1.  
Since an $\ofra_{D_{E'}}$-chain with period 1 is uniform, (1) follows from Proposition \ref{resultofBL} and \cite[II-Proposition 5.4]{BL}.  
The claim (2) follows from Proposition \ref{resultofBL} and \cite[Lemme 1.6]{S2}.  
\end{prf}

We regard $\Bscr^{E}(H', F)$ as a subset in $\Bscr^{E}(H, F)$ via $\iota_{H/H'}$, and $\Bscr^{E}(H, F)$ as a subset in $\Bscr^{E}(G, F)$ via $\iota_{G/H}$.  

\begin{prop}
\label{compoffiltC}
Let $x \in \Bscr^{E}(H', F)$ such that $[x]$ is a vertex.  
Let $\Afra$ be the hereditary $\ofra_{F}$-order in $A$ associated with $x \in \Bscr^{E}(G, F)$, and let $\Pfra$ be the radical of $\Afra$.  
We put $\mathfrak{h}(F) = \Lie(H) = B$.  
\begin{enumerate}
\item Let $n \in \Z$.  
	\begin{enumerate}
	\item $B \cap \Pfra^{n} = \mathfrak{h}(F)_{x, n/e(\Afra|\ofra_{F})}$, 
	\item $B \cap \Pfra^{n+1} = \mathfrak{h}(F)_{x, n/e(\Afra|\ofra_{F})+}$, 
	\item $B \cap \Pfra^{ \lfloor (n+1)/2 \rfloor } = \mathfrak{h}(F)_{x, n/2e(\Afra|\ofra_{F})}$, 
	\item $B \cap \Pfra^{ \lfloor n/2 \rfloor +1 } = \mathfrak{h}(F)_{x, n/2e(\Afra|\ofra_{F})+}$.  
	\end{enumerate}
\item Let $n \in \Z_{\geq 0}$.  
	\begin{enumerate}
	\item $B^{\times} \cap \U^{n}(\Afra) = H(F)_{x, n/e(\Afra|\ofra_{F})}$, 
	\item $B^{\times} \cap \U^{n+1}(\Afra) = H(F)_{x, n/e(\Afra|\ofra_{F})+}$, 
	\item $B^{\times} \cap \U^{ \lfloor (n+1)/2 \rfloor }(\Afra) = H(F)_{x, n/2e(\Afra|\ofra_{F})}$, 
	\item $B^{\times} \cap \U^{ \lfloor n/2 \rfloor +1}(\Afra) = H(F)_{x, n/2e(\Afra|\ofra_{F})+}$.  
	\end{enumerate}
\end{enumerate}
\end{prop}

\begin{prf}
By \cite[Proposition 1.9.1]{Ad}, we have $B \cap \gfra_{x, r} = \mathfrak{h}(F) \cap \gfra_{x, r} = \mathfrak{h}_{x, r}$ for $r \in \tilde{\R}$ and $B^{\times} \cap G(F)_{x, r} = H(F) \cap G(F)_{x, r} = H(F)_{x, r}$ for $r \in \tilde{\R}_{\geq 0}$.  
On the other hand, $x \in \Bscr^{E}(G, F)$ is constructed from an $\ofra_{D}$-chain by Proposition \ref{beingprincipal}.  
Then we can apply Proposition \ref{compoffiltS} and assertions follow.  
\end{prf}

\begin{lem}
\label{lemforfindc}
Let $x \in \Bscr^{E}(H', F)$ such that $[x]$ is a vertex.  
Let $\Bfra$ be the hereditary $\ofra_{E}$-order in $B$ with $x \in \Bscr^{E}(H, F)$, and let $\Qfra$ be the radical of $\Bfra$.  
\begin{enumerate}
\item For $r \in \R$, we have $\Qfra^{\lceil r \rceil} =\mathfrak{h}(F)_{x, r/e(\Bfra|\ofra_{E})e(E/F)}$.  
\item For $r \in \R_{\geq 0}$, we have $\U^{\lceil r \rceil}(\Bfra) = H(F)_{x, r/e(\Bfra|\ofra_{E})e(E/F)}$.  
\item Let $r \in \R_{\geq 0}$.  
If $H(F)_{x, r} \neq H(F)_{x, r+}$, then $n=re(\Bfra|\ofra_{E})e(E/F)$ is an integer, and we have 
	\begin{enumerate}
	\item $H(F)_{x, r} = \U^{n}(\Bfra)$, 
	\item $H(F)_{x, r+} = \U^{n+1}(\Bfra)$, and
	\item $H(F)_{x, r/2+} = \U^{\lfloor n/2 \rfloor +1}(\Bfra)$.  
	\end{enumerate}
\end{enumerate}
\end{lem}

\begin{prf}
We show (1), and (2) follows from (1).  
Let $(\mathfrak{a}_{x,r})$ be the filtration in $A$ with $x$, and let $(\mathfrak{b}_{x,r})$ be the filtration in $B$ with $x$.  
Since $[x] \in \Bscr^{R}(H', F)$ is a vertex, by Proposition \ref{beingprincipal} (1) $x \in \Bscr^{E}(H, F)$ is constructed from an $\ofra_{D_{E}}$-chain.  
Then by Proposition \ref{compoffiltA} and Proposition \ref{compoffiltG} we have $\Qfra^{\lceil r \rceil} = \mathfrak{b}_{x, r/e(\Bfra|\ofra_{E})}$.  
On the other hand, by Proposition \ref{resultofBL} (2), we also have $\mathfrak{b}_{x, r/e(\Bfra|\ofra_{E})}=B \cap \mathfrak{a}_{x, r/e(\Bfra|\ofra_{E})e(E/F)}$.  
Since $\mathfrak{a}_{x,r/e(\Bfra|\ofra_{E})e(E/F)} = \gfra(F)_{x, r/e(\Bfra|\ofra_{E})e(E/F)}$ by Proposition \ref{compoffiltA}, we obtain $\Qfra^{\lceil r \rceil} = \mathfrak{h}(F) \cap \gfra(F)_{x, r/e(\Bfra|\ofra_{E})e(E/F)} = \mathfrak{h}(F)_{x, r/e(\Bfra|\ofra_{E})e(E/F)}$, where the last equality follows from \cite[Proposition 1.9.1]{Ad}.  

To show (3), let $r \in \R_{\geq 0}$ and suppose $H(F)_{x, r} \neq H(F)_{x, r+}$.  
If $re(\Bfra|\ofra_{E})e(E/F) \notin \Z$, then $( re(\Bfra|\ofra_{E})e(E/F) , \lceil re(\Bfra|\ofra_{E})e(E/F) \rceil ]$ is nonempty and
\begin{eqnarray*}
	H(F)_{x, r+} & = & \bigcup_{r < r'} H(F)_{x, r'} \\
	& = & \bigcup_{re(\Bfra|\ofra_{E})e(E/F) < r' \leq \lceil re(\Bfra|\ofra_{E})e(E/F) \rceil} H(F)_{x, r'/e(\Bfra|\ofra_{E})e(E/F)} \\
	& = & \U^{\lceil r' \rceil}(\Bfra) = \U^{\lceil re(\Bfra|\ofra_{E})e(E/F) \rceil}(\Bfra) \\
	& = & H(F)_{x,r}, 
\end{eqnarray*}
which is a contradiction.  Therefore we have $n=re(\Bfra|\ofra_{E})e(E/F) \in \Z$.  
We put $G' = \underline{\Aut}_{D \otimes_{F} E}(V)$.  
Then we can regard $x$ as an element in $\Bscr^{E}(G', E)$, and for any $r' \in \R_{\geq 0}$ we have $\U^{\lceil r' \rceil}(\Bfra) = G'(E)_{x, r'/e(\Bfra|\ofra_{E})}$ by Proposition \ref{compoffiltG} (2).   
Therefore we obtain $H(F)_{x,r'} = \U^{\lceil r'e(\Bfra|\ofra_{E})e(E/F) \rceil}(\Bfra) = G'(E)_{x, r'e(E/F)}$ and $H(F)_{x,r'+} = G'(E)_{x, r'e(E/F)+}$ for $r' \in \R$.  
Then by Proposition \ref{compoffiltS} (2)-(a), (b) and (d), 
\begin{eqnarray*}
	\U^{n}(\Bfra) = & G'(E)_{x,re(E/F)} & = H(F)_{x,r}, \\
	\U^{n+1}(\Bfra) = & G'(E)_{x, re(E/F)+} & = H(F)_{x,r+} \\
	\U^{\lfloor n/2 \rfloor +1}(\Bfra) = & G'(E)_{x, re(E/F)/2+} & = H(F)_{x, r/2+}, 
\end{eqnarray*}
which completes the proof of (3).  
\end{prf}

\section{Generic elements and generic characters of $G$}
\label{Generic}

In this section, we discuss generic elements and generic characters, using descriptions of tame twisted Levi subgroups in $G$, given in \S \ref{TTLS}.  

\subsection{Standard representatives}

We fix a uniformizer $\varpi_{F}$ of $F$.  
Let $E$ be a finite, tamely ramified extension of $F$.  
Then we can consider the subgroup $C_{E}$ of ``standard representatives" in $E^{\times}$ defined in \cite[\S 5]{May}.  

We recall the construction of $C_{E}$.  
Since $E/F$ is tamely ramified, there exist a uniformizer $\varpi_{E}$ of $E$ and a root $\zeta$ of unity in $E$ with order prime to the residual characteristic $p$ of $E$, such that $\varpi_{E}{}^{v_{E}(\varpi_{F})} \zeta = \varpi_{F}$.  

\begin{defn}[{{\cite[Definition 5.3]{May}}}]
Let $E/F$ be a finite, tamely ramified extension, and let $\varpi_{E} \in E$ be as above.  
We denote by $C_{E}$ the subgroup in $E^{\times}$ which is generated by $\varpi_{E}$ and roots of unity in $E$ with order prime to $p$.  
\end{defn}

By \cite[Proposition 5.4]{May}, this definition of $C_{E}$ is independent of the choice of $\varpi_{E}$.  
Then $C_{E}$ depends only the choice of $\varpi_{F}$, which we already fixed.  

We recall properties of $C_{E}$.  

\begin{prop}
\label{propforsr}
Let $E/F$ be a finite, tamely ramified extension.  
\begin{enumerate}
\item Let $c \in E^{\times}$.  
Then there exists a unique $\sr(c) \in C_{E}$, called the standard representative of $c$, such that $\sr(c) \in c(1+\pfra_{E})$.  
\item For any $c \in E^{\times}$, the standard representative $\sr(c)$ is the unique element in $C_{E}$ such that $\ord \left( \sr(c)-c \right) >\ord(c)$.  
\item Let $E'/E$ be also a finite, tamely ramified extension.  
Then we have an inclusion $C_{E} \subset C_{E'}$ as groups.  
\item Let $s \in C_{E}$.  
Let $\sigma_{1}, \sigma_{2} \in \Hom_{F}(E, \bar{F})$ such that $\sigma_{1}(s) \neq \sigma_{2}(s)$.  
Then we have $\ord \left( \sigma_{1}(s) - \sigma_{2}(s) \right) = \ord (s)$.  
\end{enumerate}
\end{prop}

\begin{prf}
The claims (1), (2) and (3) are results by Mayeux \cite[Proposition 5.5, 5.6(i)]{May}.  
The claim (4) is proved in \cite[Proposition 5.6 (ii)]{May} when $E/F$ is Galois.  
For general $E$, since $E/F$ is tamely ramified, then $E/F$ is separable and we can take the Galois closure $\tilde{E}$ of $E$ in $\bar{F}$.  
Then $\tilde{E}/F$ is a finite Galois extension.  
Let $\tilde{\sigma}_{1}, \tilde{\sigma}_{2} \in \Hom_{F}(\tilde{E}, \bar{F})$ be a extension of $\sigma_{1}, \sigma_{2}$, respectively.  
By (3), $s \in C_{E} \subset C_{\tilde{E}}$.  
We also have $\tilde{\sigma}_{1}(s) = \sigma_{1}(s) \neq \sigma_{2}(s) = \tilde{\sigma}_{2}(s)$.  
Therefore by applying (4) for $\tilde{E}/F$ we have
\[
	\ord (s) = \left( \tilde{\sigma}_{1}(s) -\tilde{\sigma}_{2}(s) \right) = \ord \left( \sigma_{1}(s) - \sigma_{2}(s) \right), 
\]
which is what we wanted.  
\end{prf}

By using standard representatives, we can judge whether some element in $E$ is minimal or not.  

\begin{prop}[{{\cite[Proposition 5.9]{May}}}]
\label{srformin}
Let $E/F$ be a finite, tamely ramified extension, and let $c \in E$ such that $E=F[c]$.  
Then the following assertions are equivalent.  
\begin{enumerate}
\item $c$ is minimal over $F$.  
\item $E=F[\sr(c)]$.  
\end{enumerate}
\end{prop}

\begin{lem}
\label{preservemin}
Let $E/F$ be a finite, tamely ramified extension, and let $c, c' \in E$ such that $c^{-1}c' \in 1+\pfra_{E}$.  
Then $c$ is minimal relative to $E/F$ if and only if $c'$ is minimal relative to $E/F$.  
\end{lem}

\begin{prf}
It suffices to show that if $c$ is minimal relative to $E/F$, then $c'$ is also minimal relative to $E/F$.  
Suppose $c$ is minimal relative to $E/F$.  
In particular, $E$ is generated by $c$ over $F$.  
Then we have $E=F[\sr(c)]$ by Proposition \ref{srformin}.  
Since $\sr(c) \in c(1+\pfra_{E}) = c'(1+\pfra_{E})$, we have $\sr(c') = \sr(c)$ by Proposition \label{propforsr} (1).  
If $E$ is also generated by $c'$ over $F$, then we can apply Proposition \ref{srformin} and $c'$ is minimal relative to $E/F$.  
Thus it is enough to show $E=F[c']$.  

We put $\Hom_{F}(E, \bar{F}) = \{ \tau_{1}, \ldots, \tau_{[E:F]} \}$.  
We have $\tau_{i} \neq \tau_{j}$ for distinct $i,j \in \{ 1, \ldots, [E:F] \}$ as $E/F$ is separable.  
Since $E=F[\sr(c)]$, if $i \neq j$ we have $\tau_{i}(\sr(c)) \neq \tau_{j}(\sr(c))$ and $\ord \left( \tau_{i}(\sr(c'))-\tau_{j}(\sr(c')) \right) = \ord(c')$ by Proposition \ref{propforsr} (4).  
On the other hand, since $\ord( \sr(c') -c' ) > \ord(c')$ by Proposition \ref{propforsr} we have
\[
	\ord \left( \tau_{i}(\sr(c')-c') \right) = \ord \left( \sr(c')-c' \right) > \ord(c').  
\]
For $i \neq j$, we obtain
\begin{eqnarray*}
	& & \ord ( \tau_{i}(c') - \tau_{j}(c') ) \\
	& = & \ord \Bigl( \bigl( \tau_{i}(\sr(c')) - \tau_{j}(\sr(c')) \bigr) - \bigl( \tau_{i}(\sr(c')-c') \bigr) + \bigl( \tau_{j}(\sr(c')-c') \bigr) \Bigr),  
\end{eqnarray*}
and then 
\[
	\ord ( \tau_{i}(c') - \tau_{j}(c') ) = \ord \left( \tau_{i}(\sr(c'))-\tau_{j}(\sr(c')) \right) = \ord(c') \in \R.  
\]
In particular, we have $\tau_{i}(c') \neq \tau_{j}(c')$.  
Since $\Hom_{F}(E, \bar{F}) = \{ \tau_{1}, \ldots, \tau_{[E:F]} \}$, the element $c'$ generates $E$ over $F$.  
\end{prf}

\subsection{Concrete description of \textbf{GE1} for $G$}
\label{defofXc}

Let $E'/E/F$ be a tamely ramified field extension in $A$.  
We put 
\[
	H=\Res_{E/F} \underline{\Aut}_{D \otimes _{F} E}(V), \, H'=\Res_{E'/F} \underline{\Aut}_{D \otimes _{F} E'}(V).  
\]
Then $(H', H, G)$ is a tame twisted Levi sequence by Corollary \ref{findTTLS}.  
And also, we have a natural isomorphism $\Lie(H') \cong \End_{D \otimes _{F} E'}(V)$.  
For $c \in \End_{D \otimes _{F} E'}(V)$, we can define $X_{c}^{*} \in \Lie^{*}(H')$ as
\[
X_{c}^{*}(z)= \Tr_{E'/F} \circ \Trd_{\End_{D \otimes_{F} E'}(V)/E'}(cz), 
\]
for $z \in \Lie(H') \cong \End_{D \otimes_{F} E'} (V)$, where $cz$ is a product of $c$ and $z$ as elements in $\End_{D \otimes _{F} E'}(V)$.  
Since $E'/F$ is separable, $Tr_{E'/F}$ is surjective and there exists $e' \in E'$ such that $\Tr_{E'/F}(e') \neq 0$.  
Here, suppose $c \neq 0$.  
Since the map $(c, z) \mapsto \Trd_{\End_{D \otimes_{F} E'}(V)/E'}(cz)$ is a non-degenerate bilinear form on $\End_{D \otimes _{F} E'}(V)$, there exists $z \in \End_{D \otimes _{F} E'}(V)$ such that $\Trd_{\End_{D \otimes_{F} E'}(V)/E'}(cz) = e'$.  
In this case, we have $X_{c}^{*}(z) \neq 0$.  
Then, the map $c \mapsto X_{c}^{*}$ gives an isomorphism
\[
\Lie(\Res_{E/F} \underline{\Aut}_{D \otimes _{F} E}(V)) \cong \Lie^{*}(\Res_{E/F} \underline{\Aut}_{D \otimes_{F} E}(V)).  
\]
Since $\Trd_{A/F} |_{\End_{D \otimes_{F} E'}(V)} = \Tr_{E'/F} \circ \Trd_{\End_{D \otimes _{F} E'}(V)/E'}$, we also have
\[
X_{c}^{*}(z) = \Trd_{A/F}(cz).  
\]
For any $h \in H'(F)$ and $z \in \Lie(H')$, we have
\[
X_{c}^{*}(hzh^{-1})=\Trd_{A/F}(chzh^{-1})=\Trd_{A/F}(h^{-1}chz)=X_{h^{-1}ch}^{*}(z).  
\]
Then the linear form $X_{c}^{*}$ is invariant under $H'(F)$-conjugation if and only if $c = h^{-1}ch$ for any $h \in H'(F)=\Aut_{D \otimes_{F}E'}(V)$, that is, $c \in \Cent \left( \End_{D \otimes_{F} E}'(V) \right) = E'$.  

Let $c \in E'^{\times}$.  
We denote by $X_{c, \bar{F}}^{*}$ the image of $X_{c}^{*}$ in $\Lie^{*} \left( Z(H') \times_{F} \bar{F} \right)$.  
In other words, we put $X_{c, \bar{F}}^{*} = X_{c}^{*} \otimes_{F} \id_{\bar{F}} \in \Lie^{*} \left( Z(H') \right) \otimes_{F} \bar{F} \cong \Lie^{*} \left( Z(H') \times_{F} \bar{F} \right)$.  

To describe $X_{c,\bar{F}}^{*}(H_{\alpha})$ concretely, we use the notations in \S \ref{TTLS}.  

\begin{prop}
\label{concreteGE1}
Let $c \in E'^{\times}$ and $\alpha=\alpha_{(i',j',k'),(i'',j'',k'')} \in \Phi(G, T; \bar{F})$.  
Then we have $X_{c,\bar{F}}^{*}(H_{\alpha})=\sigma_{i',j'}(c)-\sigma_{i'',j''}(c)$.  
\end{prop}

\begin{prf}
Let $z=\sum_{i}z_{i} \otimes_{F} a_{i} \in \Lie^{*} (G) \otimes_{F} \bar{F} \cong \Lie \left( G \times_{F} \bar{F} \right)$.  
Then we have 
\begin{eqnarray*}
	X_{c,\bar{F}}^{*}(z) & = & \sum_{i} \Trd_{A/F}(cz_{i}) \otimes_{F} a_{i} = \sum_{i} \Tr_{A \otimes_{F} \bar{F} / \bar{F}} (cz_{i} \otimes_{F} a_{i}) \\
	& = & \Tr_{A \otimes_{F} \bar{F}/\bar{F}} \left( (c \otimes_{F} 1) \sum_{i}z_{i} \otimes_{F} a_{i} \right) \\
	& = & \Tr_{. \left( \End_{D \otimes \bar{F}}(\bigoplus_{i,j,k}V_{i,j,k}) \right) /\bar{F}} (m_{c,\bar{F}} z),  
\end{eqnarray*}
where $\End_{D \otimes \bar{F}}(\bigoplus_{i,j,k}V_{i,j,k}) \cong \M_{|I_{1}| \times |I_{2}| \times |I_{3}|}\left( \End_{D \otimes \bar{F}} (V \otimes_{L} \bar{F}) \right) \cong \M_{[L:F]}(\bar{F})$.  
Then, to calculate $\Tr_{A \otimes \bar{F}/\bar{F}}(m_{c,\bar{F}}H_{\alpha})$ we consider the value $m_{c,\bar{F}} \circ H_{\alpha} (v_{i,j,k})$ for some $v_{i,j,k} \in V_{i,j,k} \setminus \{ 0 \}$.  
By construction of $H_{\alpha}$ and Proposition \ref{Levidecom} (iii), we obtain
\[
	m_{c,\bar{F}} \circ H_{\alpha} (v_{i,j,k}) = \begin{cases}
		\sigma_{i',j',k'}(c)v_{i',j',k'} & \left( (i,j,k) = (i',j',k') \right), \\
		-\sigma_{i'',j'',k''}(c)v_{i'',j'',k''} & \left( (i,j,k) = (i'',j'',k'') \right), \\
		0 & otherwise.  
		\end{cases}
\]
Then we have $X_{c}^{*}(H_{\alpha})=\Tr_{A \otimes \bar{F}/\bar{F}}(m_{c,\bar{F}}H_{\alpha})=\sigma_{i',j',k'}(c)-\sigma_{i'',j'',k''}(c)$.  
Since $c \in E'$, we have $\sigma_{i',j',k'}(c)=\sigma_{i',j'}(c)$ and $\sigma_{i'',j'',k''}(c) = \sigma_{i'',j''}(c)$, which complete the proof.  
\end{prf}

\subsection{General elements of $G$}

\begin{prop}
\label{genelem}
Let $c \in E'^{\times}$.  
We put $r = -\ord(c)$.  
\begin{enumerate}
\item $X_{c}^{*} \in \Lie^{*}(Z(H'))_{-r}$.  
\item $X_{c}^{*}$ is $H$-generic of depth $r$ if and only if $c$ is minimal relative to $E'/E$.  
\end{enumerate}
\end{prop}

\begin{prf}
We show (1).  
By definition of the filtration on $\Lie^{*} \left( Z(H') \right)$, the lattice $\Lie^{*} \left( Z(H') \right) _{-r}$ is the set of $F$-linear forms $f$ of $\Lie \left( Z(H') \right)$ such that $f \left( \Lie \left(Z(H') \right)_{r+} \right) \subset \pfra_{F}$.  
Let $z \in \Lie \left( Z(H') \right) \cong E'$.  
Then $z \in \Lie \left( Z(H') \right)_{r+}$ if and only if $\ord(z) > r$.  
Therefore $c\Lie \left( Z(H') \right)_{r+} \subset \pfra_{E}$.  
Since $\Trd_{\End_{D \otimes_{F} E}(V)/E}(\pfra_{E}) \subset \pfra_{E}$ and $\Tr_{E/F}(\pfra_{E}) \subset \pfra_{F}$, we have
\[
	X_{c}^{*} \left( \Lie \left( Z(H') \right)_{r+} \right) \subset \Tr_{E/F} \circ \Trd_{\End_{D \otimes_{F} E}(V)/E} (\pfra_{E}) \subset \pfra_{F}.  
\]

To show (2), first suppose $X_{c}^{*}$ is $H$--generic of depth $r$.  

We will show $E'=E[c]$.  
We fix an embedding $\sigma_{i}:E \to \bar{F}$.  
Then we have $\Hom_{E}(E', \bar{F}) = \{ \sigma_{i,j} \mid j \in I_{2} \}$.  
Since $E/F$ is separable, to show $E'=E[c]$ it suffices to show $\sigma_{i,j}(c) \neq \sigma_{i,j'}(c)$ for any distinct $j,j' \in I_{2}$.  
We fix $k \in I_{3}$ and we put $\alpha = \alpha_{(i,j,k),(i,j',k)} \in \Phi(G, T; \bar{F})$.  
Then $\alpha \in \Phi(H, T; \bar{F}) \setminus \Phi(H', T; \bar{F})$.  
Since $X_{c}^{*}$ is $H$-generic of depth $r$, we have $-r=\ord \left( X_{c,\bar{F}}^{*}(H_{\alpha}) \right)=\ord \left( \sigma_{i,j}(c)-\sigma_{i,j'}(c) \right)$, where the last equality follows from Proposition \ref{concreteGE1}.  
In particular, we have $\ord \left( \sigma_{i,j}(c)-\sigma_{i,j'}(c) \right) \in \R$.  
Then $\sigma_{i,j}(c)-\sigma_{i,j'}(c) \neq 0$, that is, $\sigma_{i,j}(c) \neq \sigma_{i,j'}(c)$.  

Therefore, we can apply Proposition \ref{srformin} to show $c$ is minimal, and it is enough to show $\sigma_{i,j}(\sr(c)) \neq \sigma_{i,j'}(\sr(c))$ for any distinct $j,j' \in I_{2}$.  
We already have $-r=\ord \left( \sigma_{i,j}(c)-\sigma_{i,j'}(c) \right)$.  
On the other hand, we also have
\begin{eqnarray*}
	& & \ord \left( \sigma_{i,j}(\sr(c)) - \sigma_{i,j'}(\sr(c)) \right) \\
	& = & \ord \Bigl( \bigl( \sigma_{i,j}(c) - \sigma_{i,j'}(c) \bigr) + \bigl( \sigma_{i,j}(\sr(c)-c)-\sigma_{i,j'}(\sr(c)-c) \bigr) \Bigr).  
\end{eqnarray*}
Since $\ord ( \sr(c)-c ) > \ord(c)=-r$ by Proposition \ref{propforsr} (2), we have 
\[
	 \ord \bigl( \sigma_{i,j}(\sr(c)-c)-\sigma_{i,j'}(\sr(c)-c) \bigr) > -r = \ord \left( \sigma_{i,j}(c)-\sigma_{i,j'}(c) \right).  
\]
Then we obtain $\ord \left( \sigma_{i,j}(\sr(c)) - \sigma_{i,j'}(\sr(c)) \right) = \ord \left( \sigma_{i,j}(c)-\sigma_{i,j'}(c) \right) = -r \in \R$, and $\sigma_{i,j}(\sr(c)) \neq \sigma_{i,j'}(\sr(c))$.  

Conversely, suppose $c$ is minimal relative to $E'/E$.  
In particular, we have $E'=E[c]$.  
By Corollary \ref{GE1toGE2}, to show that $X_{c}^{*}$ is $H$-generic it suffices to check $X_{c}^{*}$ satisfies \textbf{GE1}.  
Let $\alpha = \alpha_{(i,j,k), (i',j',k')} \in \Phi(H, T; F) \setminus \Phi(H', T; F)$.  
Then we have $i=i'$ and $j \neq j'$.  
We equip $\bar{F}$ with $E$-structure via $\sigma_{i}$.  
Then we have $\Hom_{E}(E', \bar{F}) = \{ \sigma_{i,j} \mid j \in I_{2} \}$.  
Since $c$ is minimal over $E$ and $E'=E[c]$, the element $\sr(c)$ also generate $E'$ over $E$ by Proposition \ref{srformin}, and $\sigma_{i,j}(\sr(c)) \neq \sigma_{i,j'}(\sr(c))$.  
Then we have $\ord \bigl( \sigma_{i,j}(\sr(c)) - \sigma_{i,j'}(\sr(c)) \bigr) = \ord (\sr(c))$ by Proposition \ref{propforsr} (4).  
Moreover, $\ord (\sr(c) - c) > \ord(c)$ by Proposition \ref{propforsr}, and we have $\ord(\sr(c))=\ord(c)=-r$.  
On the other hand, we also have $\ord \bigl( \sigma_{i,j}(\sr(c)-c)-\sigma_{i,j'}(\sr(c)-c) \bigr) > -r$ as in the same way as above.    
Therefore
\begin{eqnarray*}
	X_{c}^{*}(H_{\alpha}) & = & \ord \left( \sigma_{i,j}(c) - \sigma_{i,j'}(c) \right) \\
	& = & \ord \Bigl( \bigl( \sigma_{i,j}(\sr(c)) - \sigma_{i,j'}(\sr(c)) \bigr) - \bigl( \sigma_{i,j}(\sr(c)-c)-\sigma_{i,j'}(\sr(c)-c) \bigr) \Bigr) \\
	& = & -r, 
\end{eqnarray*}
which implies $X_{c}^{*}$ is $H$-generic of depth $r$.  
\end{prf}

\subsection{General characters of $G$}

In this subsection, we discuss smooth characters of $G$.  
The goal of this subsection is to prove the following proposition.  

\begin{prop}
\label{propforgench}
Let $\chi$ be a smooth character of $G$.  
Let $\Afra$ be a principal hereditary $\ofra_F$--order.  
Suppose $\chi$ is trivial on $\U^{n+1}(\Afra)$, but not trivial on $\U^{n}(\Afra)$ for some $n \in \Z_{ \geq 0}$.  
Then there exists $c \in F$ such that $v_{F}(c)=-n/e(\Afra|\ofra_{F})$ and
\[
	\chi |_{\U^{\lfloor n/2 \rfloor +1}(\Afra)} (1+y) = \psi \circ \Trd_{A/F}(cy)
\]
for $y \in \Pfra^{\lfloor n/2 \rfloor +1}$.  
\end{prop}

To prove Proposition \ref{propforgench}, we need some preliminary.  
We put $e=e(\Afra|\ofra_{F})$.  
If Proposition \ref{propforgench} holds for some $\chi, n$ and $\Afra$, then it also holds for any $G$-conjugation of $\Afra$ and the same $\chi, n$ as above.  
Therefore we may assume
\[
	\Afra = \left(
	\begin{array}{ccc}
		\M_{md/e}(\ofra_{D}) & \cdots & \M_{md/e}(\ofra_{D}) \\
		\vdots & \ddots & \vdots \\
		\M_{md/e}(\pfra_{D}) & \cdots & \M_{md/e}(\ofra_{D})
	\end{array}
	\right).  	
\]

\begin{lem}
\label{exttrivpart}
Suppose $\chi$ is trivial on $\U^{e(n+1)}(\Afra)$.  
Then $\chi$ is also trivial on $\U^{en+1}(\Afra)$.  
\end{lem}

\begin{prf}
Since $\chi$ factors through $\Nrd_{A/F}$, it is enough to show that
\[
\Nrd_{A/F} \left( \U^{e(n+1)}(\Afra) \right) = \Nrd_{A/F} \left( \U^{en+1}(\Afra) \right).  
\]
We can deduce it from the following lemma.  
\end{prf}

\begin{lem}
We have $\Nrd_{A/F}(1+\Pfra^{n})=1+\pfra_{F}^{\lceil n/e \rceil}$.  
\end{lem}

\begin{prf}
First we show the lemma when $A$ is split.  
In this case, we have 
\[
	1+\Pfra^{n} = \left(
	\begin{array}{ccc}
		1+\M_{N/e} \left( \pfra_{F}^{\lceil n/e \rceil} \right) & & * \\
		 & \ddots & \\
		** & & 1+\M_{N/e} \left( \pfra_{F}^{\lceil n/e \rceil} \right)
	\end{array}
	\right), 
\]
where each block in $**$ is equal to $\M_{N/e} \left( \pfra_{F}^{\lceil n/e \rceil} \right)$ or $\M_{N/e} \left( \pfra_{F}^{\lceil n/e \rceil +1} \right)$.  
Then any element $a$ in $1+\Pfra^{n}$ are upper triangular modulo $\pfra_{F}^{\lceil n/e \rceil}$, and $\det_{A/F}(a)$ is 1 modulo $\pfra_{F}^{\lceil n/e \rceil}$, whence $\det_{A/F} \left( 1+\Pfra^{n} \right) \subset 1 + \pfra_{F}^{\lceil n/e \rceil}$.  
To obtain the converse inclusion, let $1+b \in 1+\pfra_{F}^{\lceil n/e \rceil}$.  
Let $a$ be an element in $A$ with the $(1,1)$-entry $b$, and other entries 0.  
Then $1+a \in 1+\Pfra^{n}$ and $\det_{A/F} (1+a) = 1+b$.  

In general case, we take a maximal unramified extension $E/F$ in $D$.  
Then $A \otimes_{F} E$ is split, and the subring $\Afra_{E}:=\Afra \otimes_{\ofra_{F}} \ofra_{E}$ in $A \otimes_{F} E$ is a hereditary $\ofra_{E}$-order with $e(\Afra_{E}|\ofra_{E}) = e(\Afra|\ofra_{F})=e$.  
Let $\Pfra_{E}$ be the radical of $\Afra_{E}$.  
Then ${\Pfra_{E}}^{n} = \Pfra^{n} \otimes_{\ofra_{F}} \ofra_{E}$ and $\det_{A \otimes_{F}E/E}(1+{\Pfra_{E}}^{n})=1+\pfra_{E}^{\lceil n/e \rceil}$ by the split case.  
Since $\Nrd_{A/F}(A^{\times}) = \det_{A \otimes_{F} E/E} \left((A \otimes_{F} 1)^{\times} \right) = F^{\times}$, we have 
\[
	\Nrd_{A/F}(1+\Pfra^{n}) \subset \left( 1+\pfra_{E}^{\lceil n/e \rceil} \right) \cap F^{\times} = 1+\pfra_{F}^{\lceil n/e \rceil}, 
\]
where the last equality follows from the assumption $E/F$ is unramified.  

To obtain the converse inclusion, let $1+b \in 1+\pfra_{F}^{\lceil n/e \rceil}$.  
Since $E/F$ is unramified, we have $\Nrm_{E/F} \left(1+\pfra_{E}^{\lceil n/e \rceil} \right) = 1+\pfra_{F}^{\lceil n/e \rceil}$, and there exists $b' \in \pfra_{E}^{\lceil n/e \rceil}$ such that $\Nrm_{E/F}(1+b') = 1+b$.  
Let $a$ be an element in $A \cong \M_{m}(D)$ with the $(1,1)$-entry $b'$, and other entries 0.  
Then $1+a \in \Cent_{A}(E) \cong \M_{m}(E)$, and
\[
	\Nrd_{A/F}(1+a) = \Nrm_{E/F} \circ \det{}_{\Cent_{A}(E)/E}(1+a) = \Nrm_{E/F}(1+b') = 1+b.  
\]
Therefore it suffice to check $a \in \Pfra^{n}$.  
We have 
\[
	\Pfra^{n} = \left(
	\begin{array}{ccc}
		\M_{md/e} \left( \pfra_{D}^{\lceil nd/e \rceil} \right) & & * \\
		 & \ddots & \\
		** & & \M_{N/e} \left( \pfra_{D}^{\lceil nd/e \rceil} \right)
	\end{array}
	\right), 
\]
Then $a \in \Pfra^{n}$ if and only if $b' \in \pfra_{D}^{\lceil nd/e \rceil}$.  
However, $b' \in \pfra_{E}^{\lceil n/e \rceil} \subset \pfra_{D}^{\lceil n/e \rceil d} \subset \pfra_{D}^{\lceil nd/e \rceil}$, where $\lceil n/e \rceil d \geq \lceil nd/e \rceil$ since $nd/e \leq \lceil n/e \rceil d \in \Z$.  
\end{prf}

\begin{prop}
\label{correctbyF}
Suppose $n>0$.  
Furthermore, assume $\chi$ is trivial on $\U^{en+1}(\Afra)$, but not on $\U^{en}(\Afra)$.  
Then there exists $c \in F$ with $v_{F}(c) = -n$ such that
\[
	\chi |_{\U^{en}(\Afra)}(1+y) = \psi \circ \Trd_{A/F}(cy)
\]
for $y \in \Pfra^{en}$.  
\end{prop}

\begin{prf}
We have $\U^{en}(\Afra)/\U^{en+1}(\Afra) \cong \Pfra^{en}/\Pfra^{en+1}$, we can regard any smooth character $\U^{en}(\Afra)/\U^{en+1}(\Afra)$ as a smooth character of $\Pfra^{en}/\Pfra^{en+1}$.  
For any smooth character $\phi$ of $\Pfra^{en}/\Pfra^{en+1}$, there exists $c_{0} \in \Pfra^{-en}$, unique up to modulo $\Pfra^{-en+1}$, such that $\phi(y) = \psi \circ \Trd_{A/F}(cy)$ for any $y \in \Pfra^{en}$.  
Since $\chi$ is not trivial on $\U^{en}(\Afra)$, we have $c_{0} \notin \Pfra^{-en+1}$.  
By the uniqueness of $c_{0}$, it suffices to show that $c_{0} + \Pfra^{-en+1}$ contains some element $c$ in $F$ with $v_{F}(c)=-n$.  
In particular, there exists $c_{0} \in \Pfra^{-en}$ such that $\chi(1+y) = \psi \circ \Trd_{A/F}(c_{0}y)$ for any $y \in \Pfra^{en}$.  
Here, let $g \in \Kfra(\Afra)$ and $y \in \Pfra^{-en}$.  
Since $\chi$ is a character of $G$, we have $\chi(1+y) = \chi \left( g (1+y) g^{-1} \right)$.  
However, we have $g (1+y) g^{-1} = 1+gyg^{-1}$ and $gyg^{-1} \in \Pfra^{en}$ since $g \in \Kfra(\Afra)$.  
Then we obtain
\begin{eqnarray*}
	\chi \left( g (1+y) g^{-1} \right) & = & \chi \left( 1+gyg^{-1} \right) \\
	& = & \psi \circ \Trd_{A/F}(c_{0}gyg^{-1}) \\
	& = & \psi \circ \Trd_{A/F}(g^{-1}c_{0}gy).  
\end{eqnarray*}
Since $g^{-1}c_{0}g \in \Pfra^{-en}$, we have $c_{0} + \Pfra^{-en+1} = g^{-1}c_{0}g + \Pfra^{-en+1}$ by the uniqueness of $c_{0}$.  
We take $t \in F^{\times}$ such that $v_{F}(t) = -n$.  
Then we have 
\[
	tc_{0} + \Pfra = t(c_{0} + \Pfra^{-en+1}) = tg^{-1}c_{0}g+t\Pfra^{-en+1} = g^{-1}(tc_{0})g + \Pfra
\]
for $g \in \Kfra(\Afra)$.  
If we put $c' = tc_{0}$, then $c', g^{-1}c'g \in \Afra$ and $c' + \Pfra = g^{-1}c'g + \Pfra$.  
Since $c_{0} \in \Pfra^{-en} \setminus \Pfra^{-en+1}$, we have $c' \in t \left( \Pfra^{^en} \setminus \Pfra^{-en+1} \right) = \Afra \setminus \Pfra$.  
Therefore we obtain $\overline{c'} = \overline{g^{-1}c'g}$ for $g \in \Kfra(\Afra)$, where for $a \in \Afra$ we denote by $\overline{a}$ the image of $a$ in $\Afra/\Pfra$.  
By the form of $\Afra$, we have an isomorphism $\Afra/\Pfra \cong \M_{md/e}(k_{D})$ as
\begin{eqnarray*}
	\Afra/\Pfra \cong \left(
	\begin{array}{ccc}
		\M_{md/e}(k_{D}) & & \\
		 & \ddots & \\
		 & & \M_{md/e}(k_{D})
	\end{array}
	\right)
	& \ni & \left(
	\begin{array}{ccc}
		b_{1}  & \\
		 & \ddots & \\
		 & & b_{e/d}
	\end{array}
	\right) \\
	 & \mapsto & (b_{1}, \ldots, b_{e/d}) \in \prod_{i=1}^{e/d} \M_{md/e}(k_{D}).  
\end{eqnarray*}
Here, let $g \in \U(\Afra)$.  
Then $g \in \Afra$ and we have $\overline{c'} = \overline{g}^{-1} \cdot \overline{c'} \cdot \overline{g}$.  
Since $\U(\Afra) \to (\Afra/\Pfra)^{\times}$ is surjective, $\overline{c'} \in Z(\Afra/\Pfra) \cong Z \left( \prod_{i=1}^{e/d} \M_{md/e}(k_{D}) \right) = \prod_{i=1}^{e/d} k_{D}$.  
Let $(b_1, \ldots, b_{e/d})$ be the image of $\overline{c'}$ in $\prod_{i=1}^{e/d} k_{D}$.  

We take $g \in \Kfra(\Afra)$ with $v_{\Afra}(g) = -1$.  
Then $\overline{g^{-1}c'g} = (b_{2}, \ldots, b_{e/d}, \tau(b_{1}))$, where $\tau \in \Gal(k_{D}/k_{F})$ is a generator.  
Since $\overline{c'}=\overline{g^{-1}c'g}$, we have $b_{1} = b_{2} = \cdots = b_{e/d} = \tau(b_{1})$.  
Since $\tau$ is a generator of $\Gal(k_{D}/k_{F})$, the element $b_{1}$ is stabilized by $\Gal(k_{D}/k_{F})$, that is, $b_{1} \in k_{F}$.  
Therefore $\overline{c'} \in k_{F} \subset \prod_{i=1}^{e/d} k_{D}$.  
We take a lift $a$ of $b_{1}$ to $\ofra_{F}$.  
Since $\overline{c'} \neq 0$, we have $b_{1} \neq 0$ and then $a \in \ofra_{F}^{\times}$.  
Therefore $c=t^{-1}a$ satisfies the desired condition.  
\end{prf}

\begin{lem}
\label{existchar}
Let $c \in F^{\times}$ such that $v_{F}(c)=-n<0$.  
Then there exists a smooth character $\theta$ of $\Afra$ such that 
\[
	\theta | _{\U^{\lfloor en/2 \rfloor +1}(\Afra)}(1+y) = \psi \circ \Trd_{A/F}(cy) 
\]
for $y \in \Pfra^{en}$.  
\end{lem}

\begin{prf}
Since $v_{\Afra}(c)=-en$, the 4-tuple $[\Afra, en, 0, c]$ is a simple stratum.  
Then we can take an element $\theta$ in $\mathscr{C}(c, 0, \Afra)$, which is nonempty by Remark \ref{existsimpch}.  
Since $\theta$ is simple, $\theta |_{\Cent_{A}(F[c])^{\times} \cap H^{1}(c, \Afra)}$ can be extended to a character of $\Cent_{A}(F[c])^{\times}$.  
However, we have $F[c]=F$ and then $\Cent_{A}(F[c]) = A$.  
Therefore, $\theta$ can be extended to a character of $A^{\times}$.  
Since $\theta$ is simple and $c \in F$ is minimal over $F$, we have 
\[
	\theta |_{\U^{\lfloor en/2 \rfloor +1}(\Afra)}(1+y) = \psi_{c}(1+y) = \psi \circ \Trd_{A/F}(cy)
\]
for $y \in \Pfra^{\lfloor en/2 \rfloor}+1$.  
\end{prf}

Let us start the proof of Proposition \ref{propforgench}.  

\begin{prf}
First, if $n=0$, then $c=1$ satisfies the condition.  
Then we may assume $n>0$.  

If $n \notin e\Z$ and $\chi$ is trivial on $\U^{n+1}(\Afra)$, then $\chi$ is also trivial in $\U^{n}(\Afra)$ by Lemma \ref{exttrivpart}, which is a contradiction.  
Then $n \in e\Z$.  
Let $i_{0}$ be the smallest integer satisfying $\lfloor n/2 \rfloor +1 \leq ei_{0}$.  
Since $n \geq 1$, we have $i_{0} \geq 1$.  
For $i = i_{0}, \ldots, n/e$, we construct $c_{i} \in F$ and a character $\theta_{i}$ of $F^{\times}$ such that $\theta_{i} |_{\U^{\lfloor ei/2 \rfloor +1}(\Afra)} = \psi_{c_{i}}$ and $\chi \cdot \left( \prod_{j=i}^{n/e} \theta_{j} \right) ^{-1}$ is trivial on $\U^{ei}(\Afra)$, by downward induction.  

Let $i=n/e$.  
Since $\chi$ is not trivial on $\U^{n}(\Afra)$, then there exists $c_{n/e} \in F$ such that $v_{F}(c_{n/e}) = -n$ and $\chi$ is equal to $\psi_{c_{n/e}}$ by Proposition \ref{correctbyF}.  
Then we take a character $\theta_{n/e}$ of $F^{\times}$ as Lemma \ref{existchar} for $c_{i}$, and $\chi \cdot \theta_{n/e}{}^{-1}$ is trivial on $\U^{n}(\Afra)=\U^{ei}(\Afra)$.  

Let $i_{0} \leq i<n/e$, and suppose we construct $c_{j}$ and $\theta_{i}$ for $i<j \leq n/e$.  
Since $\chi \cdot \left( \prod_{j=i+1}^{n/e} \theta_{j} \right) ^{-1}$ is trivial on $\U^{e(i+1)}(\Afra)$ by induction hypothesis, it is also trivial on $\U^{ei+1}(\Afra)$ by Lemma \ref{exttrivpart}.  
If $\chi \cdot \left( \prod_{j=i+1}^{n/e} \theta_{j} \right) ^{-1}$ is also trivial on $\U^{ei}(\Afra)$, then we put $c_{i}=0$ and $\theta_{i}=1$, whence $c_{i}$ and $\theta_{i}$ satisfy the condition.  
Otherwise, there exists $c_{i} \in F$ such that $v_{F}(c_{i})=-i$ and $\chi \cdot \left( \prod_{j=i+1}^{n/e} \theta_{j} \right) ^{-1} $ is equal to $\psi_{c_{i}}$ on $\U^{ei}(\Afra)$ by Proposition \ref{correctbyF}.  
Then we take a character $\theta_{i}$ of $F^{\times}$ as Lemma \ref{existchar} for $c_{i}$, and $\chi \cdot \left( \prod_{j=i}^{n/e} \theta_{j} \right) ^{-1}$ is trivial on $\U^{ei}(\Afra)$.  

Therefore $\chi \cdot \left( \prod_{i=i_{0}}^{n/e} \theta_{i} \right) ^{-1}$ is trivial on $\U^{ei_{0}}(\Afra)$.  
By Lemma \ref{exttrivpart}, it is also trivial on $\U^{e(i_{0}-1)+1}(\Afra)$.  
Since $i_{0}$ is the smallest integer satisfying $\lfloor n/2 \rfloor +1 \leq ei_{0}$, we have $e(i_{0}-1) < \lfloor n/2 \rfloor +1$, that is, $e(i_{0}-1)+1 \leq \lfloor n/2 \rfloor +1$.  
Then $\U^{\lfloor n/2 \rfloor +1}(\Afra) \subset \U^{e(i_{0}-1)+1}(\Afra)$, whence $\chi \cdot \left( \prod_{i=i_{0}}^{n/e} \theta_{i} \right) ^{-1}$ is trivial on $\U^{\lfloor n/2 \rfloor +1}(\Afra)$.  
This implies $\chi$ is equal to $\prod_{i=i_{0}}^{n/e} \theta_{i}$ on $\U^{ \lfloor n/2 \rfloor +1}(\Afra)$.  
For $i=i_{0}, \ldots, n/e$, we have $\lfloor ei/2 \rfloor +1 \leq \lfloor e(n/e)/2 \rfloor + 1 = \lfloor n/2 \rfloor +1$.  
By construction of $\theta_{i}$, the restriction of $\theta$ to $\U^{.\lfloor n/2 \rfloor +1 }(\Afra) \subset \U^{\lfloor ei/2 \rfloor +1}$ is equal to $\psi_{c_{i}}$.  
Then $\chi$ is equal to $\prod_{i=i_{0}}^{n/e} \psi_{c_{i}} = \psi_{\left( \sum_{i=i_{0}}^{n/e}c_{i} \right)}$ on $\U^{.\lfloor n/2 \rfloor +1 }(\Afra)$.  
We put $c=\sum_{i=i_{0}}^{n/e}c_{i}$.  
Since $v_{F}(c_{n/e})=-n$ and $v_{F}(c_{i}) \geq -i > -n$ for $i=i_{0}, \ldots, (n/e)-1$, we have $v_{F}(c)=-n$, which completes the proof.  
\end{prf}

\section{Some lemmas on depth-zero simple types}
\label{Depth0}

In this section, we show some lemmas which are used when we take the ``depth-zero" part of S\'echerre--Stevens's datum or Yu's datum.  

\begin{lem}
\label{lem1forfindingdepth0}
Let $\Lambda, \Lambda'$ be extensions of a maximal simple type $(J, \lambda)$ to $\tilde{J}=\tilde{J}(\lambda)$.  
Then there exists a character $\chi$ of $\tilde{J}(\lambda)/J$ such that $\Lambda' \cong \chi \otimes \Lambda$.  
\end{lem}

\begin{prf}
Since $\Lambda|_{J} = \lambda = \Lambda'|_{J}$ is irreducible, we have $\Hom_{J}(\Lambda, \Lambda') \cong \C$.  
The group $\tilde{J}$ acts $\Hom_{J}(\Lambda, \Lambda') \cong \C$ as the character $\chi$ of $\tilde{J}
$ by
\[
	g \cdot f := \Lambda'(g) \circ f \circ \Lambda(g^{-1}) = \chi(g)f
\]
for $g \in \tilde{J}$ and $f \in \Hom_{J}(\Lambda, \Lambda')$.  
Since $f$ is a $J$--homomorphism, $\chi$ is trivial on $J$.  
We take a nonzero element $f$ in $\Hom_{J}(\Lambda, \Lambda')$.  
Then for $g \in \tilde{J}$ we have 
\[
	\Lambda'(g) \circ f = f \circ \left( \chi(g)\Lambda(g) \right) = f \circ \left( \Lambda \otimes \chi (g) \right)
\]
and an $\tilde{J}$--isomorphism $\Lambda' \cong \Lambda \otimes \chi$.  
\end{prf}

If a maximal simple type $(J, \lambda)$ is associated with a simple stratum $[\Afra, n, 0, \beta]$, we put $\hat{J} = \hat{J}(\beta, \Afra)$ as in Definition \ref{defofJhat}.  

\begin{lem}
\label{lem2forfindingdepth0}
Let $(J=\U(\Afra), \lambda)$ be a simple type of depth zero, where $\Afra$ is a maximal hereditary $\ofra_{F}$-order in $A$, and let $(\tilde{J}, \Lambda)$ be a maximal extension of $(J, \lambda)$.  
We put $\rho=\Ind_{\tilde{J}}^{\hat{J}}\Lambda$.  
\begin{enumerate}
\item $\cInd_{\hat{J}}^{G} \rho$ is irreducible and supercuspidal.  
\item $\rho$ is irreducible.  
\item $\rho$ is trivial on $\U^{1}(\Afra)$.  
\end{enumerate}
\end{lem}

\begin{prf}
Since $(\tilde{J}, \Lambda)$ is a maximal extension of a simple type of depth zero, $\cInd_{\tilde{J}}^{G} \Lambda$ is irreducible and supercuspidal.  
However, by the transitivity of compact induction, we also have $\cInd_{\tilde{J}}^{G} \Lambda = \cInd_{\hat{J}}^{G} \cInd_{\tilde{J}}^{\hat{J}} \Lambda = \cInd_{\hat{J}}^{G} \rho$, which implies (1).  

Since $\cInd_{\hat{J}}^{G} \rho$ is irreducible, $\rho$ is also irreducible, that is, (2) holds.  

To show (3), we consider the Mackey decomposition of $\Res_{J}^{\hat{J}} \Ind_{\tilde{J}}^{\hat{J}} \Lambda$.  
We have 
\[
	\Res_{J}^{\hat{J}} \Ind_{\tilde{J}}^{\hat{J}} \Lambda = \bigoplus_{g \in J \backslash \hat{J} / \tilde{J}} \Ind_{J \cap {}^{g}\tilde{J} }^{J} \Res_{J \cap {}^{g}\tilde{J}}^{{}^{g} \tilde{J}} {}^{g} \Lambda = \bigoplus_{i=0}^{i-1} {}^{h^i} \lambda, 
\]
where $l = ( \hat{J} : \tilde{J} )$ and $h \in \hat{J}$ such that the image of $h$ in $\hat{J} / J \cong \Z$ is 1.  
Since $h\U^{1}(\Afra)h^{-1}=\U^{1}(\Afra)$, the representation ${}^{h^i} \lambda$ is trivial on $\U^{1}(\Afra)$ for $i=0, \ldots, l-1$.  
Therefore $\rho$ is also trivial on $\U^{1}(\Afra)$.  
\end{prf}

\begin{lem}
\label{inf_and_ind}
Let $[\Afra, n, 0, \beta]$ be a simple stratum with $\Bfra$ maximal.  
Let $\sigma^{0}$ be an irreducible cuspidal representation of $\U(\Bfra)/\U^{1}(\Bfra)$, and let $(\tilde{J}^{0}, \tilde{\sigma}^{0})$ be a maximal extension of $(\U(\Bfra), \sigma^{0})$ in $\Kfra(\Bfra)$.  
We put $\rho = \cInd_{\tilde{J}^{0}}^{\Kfra(\Bfra)} \tilde{\sigma}^{0}$.  
We denote by $\tilde{\sigma}$ the representation $\tilde{\sigma}^{0}$ as a representation of $\tilde{J}=\tilde{J}^{0}J^{1}(\beta, \Afra)$ via the isomorphism $\tilde{J}^{0}/\U^{1}(\Bfra) \cong \tilde{J}/J^{1}(\beta, \Afra)$.  
Then $\cInd_{\tilde{J}}^{\hat{J}(\beta, \Afra)} \tilde{\sigma}$ is the representation $\rho$ regarded as a representation of $\hat{J}=\hat{J}(\beta, \Afra)$ via $\Kfra(\Bfra)/\U^{1}(\Bfra) \cong \hat{J}(\beta, \Afra)/J^{1}(\beta, \Afra)$.  
\end{lem}

\begin{prf}
Since $(\U(\Bfra), \sigma^{0})$ is a simple type of $B^{\times}$ of depth zero, $\rho$ is trivial on $\U^{1}(\Bfra)$ by Lemma \ref{lem2forfindingdepth0} (3).  
Then we can regard $\rho$ as a $\hat{J}(\beta, \Afra)$-representation.  

Since $\rho=\cInd_{\tilde{J}^{0}}^{\Kfra(\Bfra)} \tilde{\sigma}^{0}$, the dimension of $\rho$ is equal to $(\Kfra(\Bfra) : \tilde{J}^{0}) \dim \tilde{\sigma}^{0}$.  
On the other hand, the dimension of $\cInd_{\tilde{J}}^{\hat{J}} \tilde{\sigma}$ is equal to $(\hat{J} : \tilde{J}) \dim \tilde{\sigma}$.  
Since $\tilde{\sigma}$ is an extension of $\tilde{\sigma}^{0}$, we have $\dim \tilde{\sigma}^{0} = \dim \tilde{\sigma}$.  
Moreover, we also have $\Kfra(\Bfra)/\tilde{J}^{0} \cong \hat{J}/\tilde{J}$ and $(\Kfra(\Bfra) : \tilde{J}^{0})=(\hat{J} : \tilde{J})$ as $\hat{J}=\Kfra(\Bfra)J(\beta, \Afra)=\Kfra(\Bfra)\tilde{J}$ and $\Kfra(\Bfra) \cap \tilde{J} = \tilde{J}^{0}$.  
Since $\rho$ is irreducible by \ref{lem2forfindingdepth0} (2), it is enough to show that there exists a nonzero $\hat{J}$--homomorphism $\rho \to \cInd_{\tilde{J}}^{\hat{J}} \tilde{\sigma}$.  

First, since $\hat{J}$ is compact modulo center in $G$ and $\tilde{J}^{0}$ contains the center of $G$, for any subgroups $J' \subset J''$ between $\hat{J}$ and $\tilde{J}^{0}$ we have $\Ind_{J'}^{J''} = \cInd_{J'}^{J''}$.   
By the Frobenius reciprocity, $\Hom_{\tilde{J}} \left( \Ind_{\tilde{J}}^{\hat{J}} \tilde{\sigma}, \tilde{\sigma} \right) \neq 0$.  
Restricting these representations to $\tilde{J}^{0}$, we have $\Hom_{\tilde{J}^{0}} \left( \Ind_{\tilde{J}}^{\hat{J}} \tilde{\sigma}, \tilde{\sigma}^{0} \right) \neq 0$.  
Using the Frobenius reciprocity, we have $\Hom_{\Kfra(\Bfra)} \left( \Ind_{\tilde{J}}^{\hat{J}} \tilde{\sigma}, \Ind_{\tilde{J}^{0}}^{\Kfra(\Bfra)} \tilde{\sigma}^{0} \right) \neq 0$.  
Since $\Kfra(\Bfra)$ is compact modulo center, $\Kfra(\Bfra)$-representations are semisimple and $\Hom_{\Kfra(\Bfra)} \left( \Ind_{\tilde{J}^{0}}^{\Kfra(\Bfra)} \tilde{\sigma}^{0}, \Ind_{\tilde{J}}^{\hat{J}} \tilde{\sigma} \right) \neq 0$.  

Here, since $J^{1}(\beta, \Afra)$ is normal in $\hat{J}$ and $\tilde{\sigma}$ is trivial on $J^{1}(\beta, \Afra)$, the restriction of $\Ind_{\tilde{J}}^{\hat{J}} \tilde{\sigma}$ to $J^{1}(\beta, \Afra)$ is also trivial.  
Then, if we extend $\Ind_{\tilde{J}^{0}}^{\Kfra(\Bfra)} \tilde{\sigma}^{0} = \rho$ to $\hat{J}=\Kfra(\Bfra)J^1(\beta, \Afra)$ as trivial on $J^{1}(\beta, \Afra)$, there exists a nonzero $\hat{J}$-homomorphism $\rho \to \Ind_{\tilde{J}}^{\hat{J}} \tilde{\sigma}$.  
\end{prf}

The following lemma guarantees the existence of extensions of $\beta$-extensions for simple characters.  

\begin{lem}
\label{exist_ext_of_b-ext}
Let $[\Afra, n, 0, \beta]$ be a simple stratum of $A$ with $\Bfra$ maximal.  
Let $\theta \in \mathscr{C}(\beta, 0, \Afra)$, and let $\kappa$ be a $\beta$-extension of the Heisenberg representation $\eta_{\theta}$ of $\theta$ to $J(\beta, \Afra)$
\begin{enumerate}
\item There exists an extension $\hat{\kappa}$ of $\kappa$ to $\hat{J}(\beta, \Afra)$.  
\item Let $\hat{\kappa}'$ be another extension of $\eta_{\theta}$ to $\hat{J}(\beta, \Afra)$.  
Then there exists a character $\chi$ of $\hat{J}(\beta, \Afra)/J^{1}(\beta, \Afra)$ such that $\hat{\kappa}' \cong \hat{\kappa} \otimes \chi$.  
\end{enumerate}
\end{lem}

\begin{prf}
We fix $g \in \Kfra(\Bfra)$ with $v_{\Bfra}(g) = 1$.  
Since $\Kfra(\Bfra) \subset B^{\times} \subset I_{G}(\kappa)$ and $\Kfra(\Bfra)$ normalizes $J(\beta, \Afra)$, we can take a $J(\beta, \Afra)$-isomorphism $f:{}^{g}kappa \to \kappa$.  
The group $\hat{J}(\beta, \Afra)/J(\beta, \Afra)$ is a cyclic group generated by the image of $g$, and then we can define $\hat{\kappa}$ as
\[
	\hat{\kappa}(g^{l}u) = f^{l} \circ \kappa(u)
\]
for $l \in \Z$ and $u \in J(\beta, \Afra)$.  
It is enough to show $\hat{kappa}$ is a group homomorphism.  
Let $g_{1}, g_{2} \in \hat{J}(\beta, \Afra)$.  
Then there exist $l_{1}, l_{2} \in \Z$ and $u_{1}, u_{2} \in J(\beta, \Afra)$ such that $g_{i} = g^{l_{i}}u_{i}$ for $i=1,2$.  
We have $g_{1}g_{2} = g^{l_{1}+l_{2}}(g^{-l_{2}}u_{1}g^{l_{2}})u_{2}$ with $g^{-l_{2}}u_{1}g^{l_{2}} \in J(\beta, \Afra)$.  
Therefore we obtain
\begin{eqnarray*}
	\hat{\kappa}(g_{1}g_{2}) & = & f^{l_{1}+l_{2}} \circ \kappa(g^{-l_{2}}u_{1}g^{l_{2}}) \circ \kappa(u_{2}) \\
	& = & f^{l_{1}} \circ \kappa(u_{1}) \circ f^{l_{2}} \circ \kappa(u_{2}) = \hat{\kappa}(g_{1}) \circ \hat{\kappa}(g_{2}), 
\end{eqnarray*}
whence (1) holds.  

Let $\hat{\kappa}'$ be another extension of $\eta_{\theta}$ to $\hat{J}(\beta, \Afra)$.  
Then we have $\Hom_{J^{1}(\beta, \Afra)} (\hat{\kappa}, \hat{\kappa}') = \Hom_{J^{1}(\beta, \Afra)} (\eta_{\theta}, \eta_{\theta}) \cong \C$.  
The group $\hat{J}(\beta, \Afra)$ acts on $\Hom_{J^{1}(\beta, \Afra)} (\hat{\kappa}, \hat{\kappa}') \cong \C$.  
Then as in the proof of Lemma \ref{lem1forfindingdepth0} we obtain $\chi$ and (2) also holds.  
\end{prf}

The following proposition is one of the key points to construct a Yu datum from a S\'echerre--Stevens datum.  

\begin{prop}
\label{findingdepth0}
Let $(J, \lambda)$ be a maximal simple type associated to a simple stratum $[\Afra, n, 0, \beta]$.  
Let $\theta \in \mathscr{C}(\beta, 0, \Afra)$ be a subrepresentation in $\lambda$, and let $\eta_{\theta}$ be the Heisenberg representation of $\theta$.  
For any extension $\Lambda$ of $\lambda$ to $\tilde{J}$ and any extension $\hat{\kappa'}$ of $\eta_{\theta}$ to $\hat{J}$, there exists an irreducible $\Kfra(\Bfra)$--representation $\rho$ such that
\begin{enumerate}
\item $\rho |_{\U(\Bfra)}$ is trivial on $\U^1(\Bfra)$ and cuspidal as a representation of $\U(\Bfra)/\U^1(\Bfra)$, 
\item $\cInd _{\Kfra(\Bfra)}^{B^{\times}} \rho$ is irreducible and supercuspidal, and
\item regarding $\rho$ as a $\hat{J}$--representation via the isomorphism $\Kfra(\Bfra)/\U^{1}(\Bfra) \cong \hat{J}/J^1$, the representation $\hat{\kappa'} \otimes \rho$ is isomorphic to $\cInd_{\tilde{J}}^{\hat{J}} \Lambda$.  
\end{enumerate}
\end{prop}

\begin{prf}
Let $\lambda = \kappa \otimes \sigma$ be a decomposition as in Definition \ref{defofsimpletype}.  
We take an extension $\hat{\kappa}$ of $\kappa$ to $\hat{J}$, which exists by Lemma \ref{exist_ext_of_b-ext} (1).  
Then there exists a character $\chi_{1}$ of $\hat{J}/J^{1}(\beta, \Afra)$ such that $\hat{\kappa} \cong \hat{\kappa}' \otimes \chi_{1}$ by Lemma \ref{exist_ext_of_b-ext}.  
Let $\tilde{\sigma}$ be an extension of $\sigma$ to $\tilde{J}$.  
Then the $\tilde{J}$-representations $\Lambda$ and $\hat{\kappa}' \otimes \chi_{1} \otimes \tilde{\sigma}$ are extensions of $\lambda$.  
By Lemma \ref{lem1forfindingdepth0}, there exists a character $\chi_2$ of $\tilde{J}$ such that $\Lambda \cong \hat{\kappa}' \otimes \chi_1 \otimes \tilde{\sigma} \otimes \chi_2$.  
Since $\chi_2$ is trivial on $J$ and $\hat{J}/J \cong \Z$, we can extend $\chi_2$ to $\hat{J}$.  
Let $J'$ be a subgroup in $\Kfra(\Bfra)$ corresponding to $\tilde{J}$ via the isomorphism $\Kfra(\Bfra)/\U(\Bfra) \cong \hat{J}/J$.  
Then $(J', \tilde{\sigma} \otimes \chi_1 \chi_2)$ is a maximal extension of the depth zero simple type $(\U(\Bfra), \sigma)$.  
Therefore we obtain a $\Kfra(\Bfra)$--representation $\rho = \cInd_{J'}^{\Kfra(\Bfra)} (\tilde{\sigma} \otimes \chi_1 \chi_2)$.  
Regarding $\rho$ as a $\hat{J}$--representation, $\rho$ is equal to $\cInd_{\tilde{J}}^{\hat{J}} (\tilde{\sigma} \otimes \chi_1 \chi_2)$ by Lemma \ref{inf_and_ind}.  
Then we have 
\[
	\hat{\kappa'} \otimes \rho = \hat{\kappa'} \otimes \cInd_{\tilde{J}}^{\hat{J}} (\tilde{\sigma} \otimes \chi_1 \chi_2) \cong \cInd_{\tilde{J}}^{\hat{J}} (\hat{\kappa'} \otimes \tilde{\sigma} \otimes \chi_1 \chi_2) \cong \cInd_{\tilde{J}}^{\hat{J}} \Lambda.  
\]
Therefore $\rho$ satisfies the desired conditions by Lemma \ref{lem2forfindingdepth0}.  
\end{prf}

Conversely, the following proposition is used to construct S\'echerre--Stevens data from Yu data.  

\begin{prop}
\label{depth0ofYu}
Let $(x, (G^{i}), (\rbf_{i}), (\Phibf_{i}), \rho)$ be a Yu datum of $G \cong \GL_{m}(D)$.  
\begin{enumerate}
\item $[x]$ is a vertex in $\Bscr^{R}(G^{0}, F)$.  
\item There exists a simple type $(G^{0}(F)_{x}, \sigma)$ of depth zero and a maximal extension $(\tilde{J}, \tilde{\sigma})$ of $(G^{0}(F)_{x}, \sigma)$ such that $\rho \cong \Ind_{\tilde{J}}^{G^{0}(F)_{[x]}} \tilde{\sigma}$.  
\end{enumerate}
\end{prop}

\begin{prf}
In the beginning, $G^{0}$ is a tame twisted Levi subgroup in $G$ with $Z(G^{0})/Z(G)$ anisotropic.  
Then there exists a tamely ramified field extension $E_{0}/F$ in $A \cong \M_{m}(D)$ such that $G^{0}(F)$ is the multiplicative group of $\Cent_{A}(E_{0})$.  
Since $\Cent_{A}(E_{0})$ is a central simple $E_{0}$-algebra, there exists $m_{E_{0}} \in \Z_{>0}$ and a division $E_{0}$-algebra $D_{E_{0}}$ such that $\Cent_{A}(E_{0}) \cong \M_{m_{E_{0}}}(D_{E_{0}})$.  

By our assumption, $\pi := \cInd_{G^{0}(F)_{[x]}}^{G} \rho$ is an irreducible and supercuspidal representation of depth zero.  
Then there exists $y \in \Bscr^{E}(G^{0}, F)$ and an irreducible $G^{0}(F)_{y}$--representation $\sigma$ such that $[y]$ is a vertex and $(G^{0}(F)_{y}, \sigma)$ is a $[G, \pi]_{G}$--type.  
Since vertices in $\Bscr^{R}(G^{0}, F)$ are permuted transitively by the action of $G^{0}(F)$, we may assume $G^{0}(F)_{y} \supset G^{0}(F)_{x}$.  

We show that $\Ind_{G^{0}(F)_{x}}^{G^{0}(F)_{y}} \Res_{G^{0}(F)_{x}}^{G^{0}(F)_{[x]}} \rho$ has a non-trivial $G^{0}(F)_{y,0+}$--fixed part.  
Since $G^{0}(F)_{x} \cap G^{0}(F)_{y,0+} \subset G^{0}(F)_{x,0+}$, the representation $\rho$ is trivial on $G^{0}(F)_{x} \cap G^{0}(F)_{y,0+}$.  
Then $\Ind_{G^{0}(F)_{x} \cap G^{0}(F)_{y,0+}} ^{G^{0}(F)_{y,0+}} \Res_{G^{0}(F)_{x} \cap G^{0}(F)_{y,0+}} ^{G^{0}(F)_{[x]}} \rho$ has a non-trivial $G^{0}(F)_{y,0+}$--fixed part by the Frobenius reciprocity.   
However, 
\[
	\Ind_{G^{0}(F)_{x} \cap G^{0}(F)_{y,0+}} ^{G^{0}(F)_{y,0+}} \Res_{G^{0}(F)_{x} \cap G^{0}(F)_{y,0+}} ^{G^{0}(F)_{[x]}} \rho \subset \Res_{G^{0}(F)_{y,0+}}^{G^{0}(F)_{y}} \Ind_{G^{0}(F)_{x}}^{G^{0}(F)_{y}} \Res_{G^{0}(F)_{x}}^{G^{0}(F)_{[x]}} \rho
\]
by the Mackey decomposition.  
Therefore $\Ind_{G^{0}(F)_{x}}^{G^{0}(F)_{y}} \Res_{G^{0}(F)_{x}}^{G^{0}(F)_{[x]}} \rho$ has a non-trivial $G^{0}(F)_{y,0+}$--fixed part.  

Since $\Ind_{G^{0}(F)_{x}}^{G^{0}(F)_{y}} \Res_{G^{0}(F)_{x}}^{G^{0}(F)_{[x]}} \rho \subset \Res_{G^{0}(F)_{y}}^{G} \cInd_{G^{0}(F)_{[x]}}^{G} \rho = \Res_{G^{0}(F)_{y}}^{G} \pi$ by the Mackey decomposition, then we also may assume $\sigma \subset \Ind_{G^{0}(F)_{x}}^{G^{0}(F)_{y}} \Res_{G^{0}(F)_{x}}^{G^{0}(F)_{[x]}} \rho$ by $G^{0}(F)_{[x]}$--conjugation if necessary by \cite[Theorem 5.5(ii)]{GSZ}.  
By the Frobenius reciprocity, $\Res_{G^{0}(F)_{x}}^{G^{0}(F)_{y}} \sigma$ is a subrepresentation of $\Res_{G^{0}(F)_{x}}^{G^{0}(F)_{[x]}} \rho$, which is trivial on $G^{0}(F)_{x,0+}$.  
Therefore, $\sigma$ has a non-trivial $G^{0}(F)_{x,0}G^{0}(F)_{y,0+}$--fixed part.  
Since the image of $G^{0}(F)_{x,0}$ in $G^{0}(F)_{y}/G^{0}(F)_{y,0+}$ is a parabolic subgroup of $G^{0}(F)_{y}/G^{0}(F)_{y,0+}$ and $\sigma$ is cuspidal when we regard $\sigma$ as a $G^{0}(F)_{y}/G^{0}(F)_{y,0+}$--representation, we have $G^{0}(F)_{x,0}G^{0}(F)_{y,0+} = G^{0}(F)_{y}$, which implies $[x]=[y]$, that is, (1) holds.  

To show (2), let $(\tilde{J}, \tilde{\sigma})$ be the unique extension of $(G^{0}(F)_{x}, \sigma)$ such that $\pi \cong \cInd_{\tilde{J}}^{G} \tilde{\sigma}$.  

We show the $G^{0}(F)_{x,0+}$--fixed part in $\pi$ is contained in $\Ind_{\tilde{J}}^{G^{0}(F)_{[x]}} \tilde{\sigma}$.  
By the Mackey decomposition, we have 
\[
	\pi \cong \bigoplus_{g \in G^{0}(F)_{x,0+} \backslash G / \tilde{J} } \Ind_{G^{0}(F)_{x,0+} \cap {}^{g}\tilde{J}}^{G^{0}(F)_{x,0+}} \Res_{G^{0}(F)_{x,0+} \cap {}^{g}\tilde{J}}^{{}^{g}\tilde{J}} {}^{g} \tilde{\sigma}.  
\]
We put $\tau(g)=\Ind_{G^{0}(F)_{x,0+} \cap {}^{g}\tilde{J}}^{G^{0}(F)_{x,0+}} \Res_{G^{0}(F)_{x,0+} \cap {}^{g}\tilde{J}}^{{}^{g}\tilde{J}} {}^{g} \tilde{\sigma}$.  
Suppose $\tau(g)$ has a non-trivial $G^{0}(F)_{x,0+}$--fixed part.  
Then $\Hom_{G^{0}(F)_{x,0+} \cap {}^{g}\tilde{J}}(\mathbf{1}, {}^{g}\tilde{\sigma}) \neq 0$ by the Frobenius reciprocity.  
Here, since $[x]$ is a vertex, $G^{0}(F)_{x}$ is a maximal compact open subgroup.  
Therefore we may assume $G^{0}(F)_{x}=\GL_{m_{E_{0}}}(\ofra_{D_{E_{0}}})$ by $G$--conjugation if necessary.  
Then there exist $k,k' \in G^{0}(F)_{x}$ and a diagonal matrix $g'$ such that the $(i,i)$--coefficient of $g'$ is $\varpi_{D_{E_{0}}}^{a_i}$ with $a_1 \geq a_2 \geq \cdots \geq a_{m}$, and such that $g=kg'k'$.  
Since $G^{0}(F)_{x,0+}$ is normal in $G^{0}(F)_{x}$ and $G^{0}(F)_{x} \subset \tilde{J}$, the condition $\Hom_{G^{0}(F)_{x,0+} \cap {}^{g}\tilde{J}}(\mathbf{1}, {}^{g}\tilde{\sigma}) \neq 0$ holds if and only if $\Hom_{{}^{{(g')}^{-1}}G^{0}(F)_{x,0+} \cap \tilde{J}}(\mathbf{1}, \tilde{\sigma}) \neq 0$.  
Therefore $\sigma$ has a non-trivial $G^{0}(F)_{x,0+}\left( {}^{{(g')}^{-1}}G^{0}(F)_{x,0+} \cap \tilde{J} \right)$--fixed part.  
If $a_i > a_{i+1}$ for some $i$, the image of ${}^{{(g')}^{-1}}G^{0}(F)_{x,0+} \cap \tilde{J}$ in $G^{0}(F)_{x}/G^{0}(F)_{x,0+}$ is a proper parabolic subgroup, which is a contradiction since $\sigma$ is cuspidal.  
Then $g' \in D_{E_{0}}^{\times} \subset G^{0}(F)_{[x]}$ and $g=kg'k' \in G^{0}(F)_{[x]}$.  
Therefore the $G^{0}_{x,0+}$--fixed part in $\pi$ is contained in $\bigoplus_{g \in G^{0}(F)_{x,0+} \backslash G^{0}(F)_{[x]} / \tilde{J} } \tau(g) = \cInd_{\tilde{J}}^{G^{0}(F)_{[x]}} \tilde{\sigma}$.  

Then we have $\rho \subset \Ind_{\tilde{J}}^{G^{0}(F)_{[x]}} \tilde{\sigma}$.  
Since these representations are irreducible, we obtain $\rho = \Ind_{\tilde{J}}^{G^{0}(F)_{[x]}} \tilde{\sigma}$.  
\end{prf}

\section{Factorization of tame simple characters}

\label{Factorization}
Let $[\Afra, n, 0, \beta]$ be a tame simple stratum of $A$.  
If $n=0$, suppose $\beta \in \ofra_{F}^{\times}$.  
By Proposition \ref{appfortame}, there exists a defining sequence $\left( [\Afra, n, r_i, \beta_i] \right)_{i=0, 1, \ldots, s}$ of $[\Afra, n, 0, \beta]$ such that
\begin{enumerate}
\item $F[\beta_{i}] \supsetneq F[\beta_{i+1}]$, 
\item $\beta_{i}-\beta_{i+1}$ is minimal over $F[\beta_{i+1}]$ and
\item $v_{\Afra}(\beta_{i}-\beta_{i+1})=k_{0}(\beta_{i}, \Afra)=r_{i+1}$
\end{enumerate}
for $i=0, 1, \ldots, s-1$.  

We put $E_i=F[\beta_i]$.  
Let $B_{i}$ be the centralizer of $E_i$ in $A$.  
Let $c_{i}=\beta_{i}-\beta_{i+1}$ for $i=0, \ldots, s-1$ and let $c_{s}=\beta_{s}$.  

\begin{prop}
\label{factchar}
Let $0 \leq t < -k_{0}(\beta, \Afra)$.  
Let $\theta \in \mathscr{C}(\beta, t, \Afra)$.  
Then for $i=0, 1, \ldots, s$ there exists a smooth character $\phi_{i}$ of $E_{i}^{\times}$ such that we have $\theta = \prod_{i=0}^{s} \theta^{i}$, where the characters $\theta^{i}$ of $H^{t+1}(\beta, \Afra)$ are defined as in the following:  
\begin{enumerate}
\item $\theta^{i} | _{B_{i}^{\times} \cap H^{t+1}(\beta, \Afra)} = \phi_i \circ \Nrd_{B_i/E_i}$, and
\item $\theta^{i} | _{H^{t_{i}+1}(\beta, \Afra)} = \psi_{c_i}$, where $t_{i}=\max \left\{ t, \lfloor -v_{\Afra}(c_{i})/2 \rfloor \right\}$.  
\end{enumerate}
\end{prop}

\begin{prf}
We show this proposition by induction on the length $s$ of the defining sequence.  

First, suppose that $s=0$, that is, $\beta$ is minimal over $F$.  
We have $\theta = \theta^{0}$.  
Then it is enough to show that $\theta$ satisfies (1) and (2).  
Since $\theta$ is simple, $\theta | B_{0}^{\times} \cap H^{t+1}(\beta, \Afra)$ factors through $\Nrd_{B_{0}/E_{0}}$.  
Then there exists a character $\phi_{0}$ of $E_{0}^{\times}$ such that $\theta = \phi_{0} \circ \Nrd_{B_{0}/E_{0}}$, whence (1) holds.  
We have $v_{\Afra}(c_0)=v_{\Afra}(\beta)=-n$ and $t_{0}=\max \{ t, \lfloor -v_{\Afra}(c_{0})/2 \rfloor \} \geq \lfloor n/2 \rfloor$.  
Then we have $H^{t_{0}+1}(\beta, \Afra) \subset H^{\lfloor n/2 \rfloor +1}(\beta, \Afra)$.  
Since $\theta$ is simple, we have $\theta|_{H^{t_{0}+1}}(\beta, \Afra)=\psi_{\beta}=\psi_{c_{0}}$, whence (2) also holds.  

Next, suppose that $s>0$, that is, $\beta$ is not minimal over $F$.  
We put $t'=\max \{ t, \lfloor -k_{0}(\beta, \Afra)/2 \rfloor \}$.  
Since $k_{0}(\beta, \Afra)=v_{\Afra}(c_{0})$, we have $t' = \max \{ t, \lfloor -v_{\Afra}(c_{0})/2 \rfloor \}=t_{0}$.  
Since $\theta$ is simple, there exists $\theta' \in \mathscr{C}(\beta, t', \Afra)$ such that $\theta |_{H^{t'+1}(\beta, \Afra)} =  \psi_{c_0}\theta'$.  
By induction hypothesis, for $i=1, \ldots, s$, there exist a smooth character $\phi_{i}$ of $E_{i}^{\times}$ and a smooth character $\theta'^{i}$ of $H^{t'+1}(\beta_1, \Afra)=H^{t'+1}(\beta, \Afra)$ such that $\theta'^{i}|_{B_{i}^{\times} \cap H^{t'+1}(\beta, \Afra)}=\phi_{i} \circ \Nrd_{B_{i}/E_{i}}$ and $\theta'^{i}|_{H^{t'_{i}+1}(\beta, \Afra)}=\psi_{c_i}$, where $t'_{i}= \max \{ t', \lfloor -v_{\Afra}(c_{i})/2 \rfloor \}$.  
Here we have $r_1 < \ldots < r_{s} < n$, whence we obtain $-v_{\Afra}(c_{0}) < \ldots < -v_{\Afra}(c_{s-1}) < -v_{\Afra}(c_{s})$.  
Since $t'=\max \{ t, \lfloor -v_{\Afra}(c_{0})/2 \rfloor \}$ and $-v_{\Afra}(c_{0}) < -v_{\Afra}(c_{i})$, we have $t'_{i} = t_{i}$.  
We want to extend $\theta'^{i}$ to a character $\theta^{i}$ of $H^{t+1}(\beta, \Afra)$ as $\theta^{i} |_{B_{i}^{\times} \cap H^{t+1}(\beta, \Afra)} = \phi_{i} \circ \Nrd_{B_{i}/E_{i}}$.  
Suppose we obtain $\theta^{i}$ in such a way.  
Then $\theta^{i}$ satisfies (1) in Proposition by construction of $\theta^{i}$, and (2) in Proposition since $H^{t_{i}+1}(\beta, \Afra)=H^{t'_{i}+1}(\beta, \Afra)$ and $\theta^{i}|_{H^{t_{i}+1}(\beta, \Afra)}=\theta'^{i}|_{H^{t_{i}+1}(\beta, \Afra)}=\psi_{c_{i}}$.  

We have $B_{i}^{\times} \cap H^{t+1}(\beta, \Afra) \cap H^{t'+1}(\beta, \Afra) = B_{i}^{\times} \cap H^{t'+1}(\beta, \Afra)$, whence restrictions of $\theta'^{i}$ and $\phi_{i} \circ \Nrd_{B_{i}/E_{i}}$ to $B_{i}^{\times} \cap H^{t+1}(\beta, \Afra) \cap H^{t_{i}+1}(\beta, \Afra)$ are equal.  
Let $b_{1}, b_{2} \in B_{i}^{\times} \cap H^{t+1}(\beta, \Afra)$ and $h'_{1}, h'_{2} \in H^{t'+1}(\beta, \Afra)$ with $b_{1}h'_{1} = b_{2}h'_{2}$.  
Then $b_{1}^{-1}b_{2} = h'_{1}(h'_{2})^{-1} \in B_{i}^{\times} \cap H^{t+1}(\beta, \Afra) \cap H^{t_{i}+1}(\beta, \Afra)$ and $\phi_{i} \circ \Nrd_{B_{i}/E_{i}}(b_{1}^{-1}b_{2}) = \theta'^{i}(h'_{1}(h'_{2})^{-1})$.  
Therefore we also have 
\[
	\phi_{i} \circ \Nrd_{B_{i}/E_{i}}(b_{1})\theta'^{i}(h'_{1}) = \phi_{i} \circ \Nrd_{B_{i}/E_{i}}(b_{2}) \theta'^{i}(h'_{2}).  
\]
Then $\theta^{i}$ is well-defined as a map from $H^{t+1}(\beta, \Afra)$ to $\C^{\times}$.  

We show $\psi_{c_{i}}(b^{-1}hb)=\psi_{c_{i}}(h)$ for $b \in B_{i}^{\times} \cap H^{t+1}(\beta, \Afra)$ and $h \in H^{t_{i}+1}(\beta, \Afra)$.  
By definition of $\psi_{c_{i}}$, we have
\[
	\psi_{c_{i}}(b^{-1}hb) = \Trd_{A/F}(c_{i}(b^{-1}hb-1)).  
\]
Since $c_{i} \in E_{i}$ and $b \in B_{i} = \Cent_{A}(E_{i})$, we have 
\[
	c_{i}(b^{-1}hb-1)=c_{i}b^{-1}(h-1)b=b^{-1}c_{i}(h-1)b.  
\]
Therefore, we obtain
\[
	\psi_{c_{i}}(b^{-1}hb) = \Trd_{A/F}(b^{-1}c_{i}(h-1)b) = \Trd_{A/F}(c_{i}(h-1)) = \psi_{c_{i}}(h).  
\]

To show $\theta^{i}$ is a character, let $h_{1}, h_{2} \in H^{t+1}(\beta, \Afra)$.  
Then there exist $b_{1}, b_{2} \in B_{i}^{\times} \cap H^{t+1}(\beta, \Afra)$ and $h'_{1}, h'_{2} \in H^{t_{i}+1}(\beta, \Afra)$ such that $h_{1}=b_{1}h'_{1}$ and $h_{2}=b_{2}h'_{2}$.  
Therefore we have 
\begin{eqnarray*}
	\theta^{i}(h_{1}h_{2}) & = & \theta^{i}(b_{1}h'_{1} b_{2}h'_{2}) \\
	& = & \theta^{i}\left( (b_{1}b_{2}) (b_{2}^{-1}h'_{1}b_{2}h'_{2} \right) \\
	& = & (\phi_{i} \circ \Nrd_{B_{i}/E_{i}}) (b_{1}) (\phi_{i} \circ \Nrd_{B_{i}/E_{i}})(b_{2}) \psi_{c_{i}}(b_{2}^{-1}h'_{1}b_{2}) \psi_{c_{i}}(h'_{2}) \\
	& = &  (\phi_{i} \circ \Nrd_{B_{i}/E_{i}}) (b_{1}) \psi_{c_{i}}(h'_{1}) (\phi_{i} \circ \Nrd_{B_{i}/E_{i}})(b_{2}) \psi_{c_{i}}(h'_{2}) \\
	& = & \theta^{i}(b_{1}h'_{1}) \theta^{i}(b_{2}h'_{2}) = \theta^{i}(h_{1}) \theta^{i}(h_{2}).  
\end{eqnarray*}

We put $\theta^{0} = \theta \prod_{i=1}^{s}(\theta^{i})^{-1}$.  
To complete the proof, it is enough to show $\theta^{0}$ satisfies (1) and (2).  

To see (1), we show the restrictions of $\theta$ and $\theta^{i}$ ($i=1, \ldots, s$) to $B_{0}^{\times} \cap H^{t+1}(\beta, \Afra)$ factor through $\Nrd_{B_{0}/E_{0}}$.  
Since $\theta$ is simple, $\theta | _{B_{0}^{\times} \cap H^{t+1}(\beta, \Afra)}$ factors through $\Nrd_{B_{0}/E_{0}}$.  
We already have $\theta | _{B_{i}^{\times} \cap H^{t+1}(\beta, \Afra)} = \phi_{i} \circ \Nrd_{B_{i}/E_{i}}$.  
Since $B_{0}^{\times} \subset B_{i}^{\times}$, we have $\theta_{B_{0}^{\times} \cap H^{t+1}(\beta, \Afra)} = \phi_{i} \circ (\Nrd_{B_{i}/E_{i}} | _{B_{0}^{\times}})$.  
However, the equation $\Nrd_{B_{i}/E_{i}} | _{B_{0}^{\times}} = \mathrm{N}_{E_{0}/E_{i}} \circ \Nrd_{B_{0}/E_{0}}$ holds.  
Then $\theta^{i} |_{B_{0}^{\times} \cap H^{t+1}(\beta, \Afra)}$ factors through $\Nrd_{B_{0}/E_{0}}$.  
Therefore $\theta^{0} |_{B_{0}^{\times} \cap H^{t+1}(\beta, \Afra)}$ also factors through $\Nrd_{B_{0}/E_{0}}$, and there exists a character $\phi_{0}$ of $E_{0}^{\times}$ such that $\theta^{0} |_{B_{0}^{\times} \cap H^{t+1}(\beta, \Afra)} = \phi_{0} \circ \Nrd_{B_{0}/E_{0}}$.  

By restricting $\theta^{0} = \theta \prod_{i=1}^{s} (\theta^{i})^{-1}$ to $H^{t_{0}+1}(\beta, \Afra)=H^{t'+1}(\beta, \Afra)$, we have
\begin{eqnarray*}
	\theta^{0}|_{H^{t_{0}+1}(\beta, \Afra)} & = & (\theta|_{H^{t'+1}(\beta, \Afra)} \prod_{i=1}^{s} (\theta^{i}|_{H^{t'+1}(\beta, \Afra)})^{-1} \\
	& = & \psi_{c_{0}} \theta' \prod_{i=1}^{s} (\theta'^{i})^{-1} = \psi_{c_{0}} \theta' \theta'^{-1} = \psi_{c_{0}}.  
\end{eqnarray*}
Therefore, (2) also holds and complete the proof.  
\end{prf}

\section{Construction of a Yu datum from a S\'echerre--Stevens datum}
\label{SStoYu}

Let $[\Afra, n, 0, \beta]$ be a tame simple stratum, $(J(\beta, \Afra), \lambda)$ be a maximal simple type with $[\Afra, n, 0, \beta]$, and let $(\tilde{J}(\lambda), \Lambda)$ be a maximal extension of $(J(\beta, \Afra), \lambda)$.  
We construct a Yu datum from the data of $[\Afra, n, 0, \beta]$ and $\lambda$.  

We put $\Pfra = \Pfra(\Afra)$.  
Let $([\Afra ,n, r_{i}, \beta_{i}])_{i=0}^{s}$, $E_{i}$, $B_{i}$ and $c_{i}$ be as in \S \ref{Factorization}.  
For $i=0, 1, \ldots, s$, we put $G^{i}=\Res_{E_{i}/F} \underline{\Aut}_{D \otimes E_{i}}(V)$ and $\rbf_{i}=-\ord(c_{i})$.  
If $\beta_s \in F$, we put $d=s$.  
If $\beta_s \notin F$, we put $d=s+1$, $G^{d}=G$ and $\rbf_{d}=\rbf_{s}$.  
Then $(G^{0}, \ldots, G^{d})$ is a tame twisted Levi sequence by Corollary \ref{findTTLS}.  
We also put $\rbf_{-1}=0$.  
For $i=-1, 0, 1, \ldots, d$, we put $\mathbf{s}_{i}=\rbf_{i}/2$.  

\begin{prop}
\label{compofTLS}
We fix a $G^{i-1}(F)$-equivalent and affine embedding 
\[
	\iota_{i}:\Bscr^{E}(G^{i-1}, F) \hookrightarrow \Bscr^{E}(G^{i}, F)
\]
for $i=1, \ldots, d$ and we put $\tilde{\iota}_{i}=\iota_{i} \circ \cdots \circ \iota_{1}$.  
We also put $\tilde{\iota}_{0}=\id_{\Bscr^{E}(G^{0}, F)}$.  
\begin{enumerate}
\item There exists $x \in \Bscr^{E}(G^{0}, F)$ such that $[x]$ is a vertex and
\begin{enumerate}
\item $G^{0}(F)_{[x]} = \Kfra(\Bfra_{0})$, 
\item $G^{0}(F)_{x} = B_{0}^{\times} \cap \U(\Afra) = \U(\Bfra_{0})$, 
\item $G^{0}(F)_{x,0+} = B_{0}^{\times} \cap \U^{1}(\Afra)$, 
\item $\gfra^{0}(F)_{x} = B_{0} \cap \Afra = \Bfra_{0}$, and
\item $\gfra^{0}(F)_{x,0+} = B_{0} \cap \Pfra$.  
\end{enumerate}
\item For $i=1, \ldots, d$, we have
\begin{enumerate}
\item $G^{i}(F)_{\tilde{\iota}_{i}(x), \mathbf{s}_{i-1} } = B_{i}^{\times} \cap \U^{\lfloor (-v_{\Afra}(c_{i-1})+1)/2 \rfloor}(\Afra)$, 
\item $G^{i}(F)_{\tilde{\iota}_{i}(x), \mathbf{s}_{i-1}+ } = B_{i}^{\times} \cap \U^{\lfloor -v_{\Afra}(c_{i-1})/2 \rfloor +1}(\Afra)$, 
\item $G^{i}(F)_{\tilde{\iota}_{i}(x), \rbf_{i-1} } = B_{i}^{\times} \cap \U^{-v_{\Afra}(c_{i-1})}(\Afra)$, 
\item $G^{i}(F)_{\tilde{\iota}_{i}(x), \rbf_{i-1}+ } = B_{i}^{\times} \cap \U^{-v_{\Afra}(c_{i-1})+1}(\Afra)$, 
\item $\gfra^{i}(F)_{\tilde{\iota}_{i}(x), \mathbf{s}_{i-1} } = B_{i} \cap \Pfra^{\lfloor (-v_{\Afra}(c_{i-1})+1)/2 \rfloor}$, 
\item $\gfra^{i}(F)_{\tilde{\iota}_{i}(x), \mathbf{s}_{i-1}+ } = B_{i} \cap \Pfra^{\lfloor -v_{\Afra}(c_{i-1})/2 \rfloor +1}$, 
\item $\gfra^{i}(F)_{\tilde{\iota}_{i}(x), \rbf_{i-1} } = B_{i} \cap \Pfra^{-v_{\Afra}(c_{i-1})}$, and
\item $\gfra^{i}(F)_{\tilde{\iota}_{i}(x), \rbf_{i-1}+ } = B_{i} \cap \Pfra^{-v_{\Afra}(c_{i-1})+1}$.  
\end{enumerate}
\item For $i=0, \ldots, s$, we have
\begin{enumerate}
\item $G^{i}(F)_{\tilde{\iota}_{i}(x), \mathbf{s}_{i}+ } = B_{i}^{\times} \cap \U^{ \lfloor -v_{\Afra}(c_{i})/2 \rfloor +1}(\Afra)$, 
\item $G^{i}(F)_{\tilde{\iota}_{i}(x), \rbf_{i} }  = B_{i}^{\times} \cap \U^{-v_{\Afra}(c_{i})}(\Afra)$, 
\item $G^{i}(F)_{\tilde{\iota}_{i}(x), \rbf_{i}+ } = B_{i}^{\times} \cap \U^{-v_{\Afra}(c_{i})+1} (\Afra)$, 
\item $\gfra^{i}(F)_{\tilde{\iota}_{i}(x), \rbf_{i} } = B_{i} \cap \Pfra^{-v_{\Afra}(c_{i})}(\Afra)$, and
\item $\gfra^{i}(F)_{\tilde{\iota}_{i}(x), \rbf_{i}+ } = B_{i} \cap \Pfra^{-v_{\Afra}(c_{i})}(\Afra)$.  
\end{enumerate}
\end{enumerate}
\end{prop}

\begin{prf}
We find $x \in \Bscr^{E}(G^{0}, F)$.  
Since $B$ is a central simple $E_{0}$-algebra, there exists a division $E_{0}$-algebra $D_{E_{0}}$ and a right $D_{E_{0}}$-module $W_{0}$ such that $B \cong \End_{D_{E_{0}}}(W_{0})$.  
Since $\Bfra_{0}$ is a maximal hereditary $\ofra_{E_{0}}$-order in $B_{0}$, there exists an $\ofra_{D_{E_{0}}}$-chain $(\Lcal_{i})_{i \in \Z}$ in $W_{0}$ of period 1 such that $\Bfra_{0}$ is the hereditary $\ofra_{E_{0}}$-order associated with $(\Lcal_{i})$.  
Let $x \in \Bscr^{E}(G^{0}, F) \cong \Bscr^{E}(\underline{\Aut}_{D_{E_{0}}}(W_{0}), E_{0})$ be an element which corresponds to a lattice function constructed from $(\Lcal_{i})_{i \in \Z}$.  
Then by Proposition \ref{beingvertex} $[x]$ is a vertex in $\Bscr^{E}(G^{0}, F)$.  
Therefore by Proposition \ref{compoffiltG} (3) we have (1)-(a).  

To show the remainder assertion, we show $\Afra$ is the hereditary $\ofra_{F}$-order in $A$ associated with $\tilde{\iota}_{d}(x)$.  
Since $[\Afra, n, 0, \beta]$ is a stratum, $\Afra$ is $E=F[\beta]$-pure.  
Moreover, we have $\Afra \cap B_{0} = \Bfra_{0}$ by definition of $\Bfra_{0}$.  
Therefore by Proposition \ref{beingprincipal} (3) $\Afra$ is associated with $\tilde{\iota}_{d}(x)$.  
Since $v_{\Afra}(c_{i}) \in \Z_{\geq 0}$ and $v_{\Afra}(c_{i}) = \ord(c_{i})e(\Afra|\ofra_{F})$, the remainder assertions in Proposition follow from Proposition \ref{compoffiltC}.  
\end{prf}

In the following, we regard $\Bscr^{E}(G^{0}, F), \ldots, \Bscr^{E}(G^{d-1}, F)$ are subsets in $\Bscr^{E}(G, F)$ via $\tilde{\iota}_{1}, \ldots, \tilde{\iota}_{d}$.  

\begin{prop}
\label{compofHJ}
\begin{enumerate}
\item $H^{1}(\beta, \Afra)=K_{+}^{d}$, 
\item $J(\beta, \Afra)={}^{\circ}K^{d}$, 
\item $\hat{J}(\beta, \Afra)=K^{d}$.  
\end{enumerate}
\end{prop}

\begin{prf}
We show (1).  
We have $r_{i}=-v_{\Afra}(c_{i-1})=-\ord(c_{i-1})e(\Afra|\ofra_{F})=-e(\Afra|\ofra_{F}) \rbf_{i-1}$ for $i=1, \ldots, s$ and $n = -v_{\Afra}(c_{s}) = -e(\Afra|\ofra_{F}) \rbf_{s}$.  
We have $G^{0}(F)_{x,0+} = B_{0}^{\times} \cap \U^{1}(\Afra)$ by Proposition \ref{compofTLS} (1)-(c).  
For $i=1, \ldots, s$ we have $B_{i}^{\times} \cap \U^{\lfloor r_{i}/2 \rfloor +1}(\Afra) = B_{i}^{\times} \cap \U^{\lfloor -v_{\Afra}(c_{i-1})/2 \rfloor +1}(\Afra) = G^{i}(F)_{x, \mathbf{s}_{i-1}+}$ by Proposition \ref{compofTLS} (2)-(b).  
We also have $B_{s}^{\times} \cap \U^{\lfloor n/2 \rfloor +1}(\Afra) = G^{s}(F)_{x, \mathbf{s}_{s}+}$.  
If $d=s+1$, by comparing Lemma \ref{presenofHJ} (1) and Definition \ref{defofKi} (1) of $K_{+}^{d}$ we have $H^{1}(\beta, \Afra) = K^{d}$.  
If $d=s$, we have $H^{1}(\beta, \Afra) = K^{d} \U^{\lfloor n/2 \rfloor +1}(\Afra)$ and it suffices to show $K^{d} \supset \U^{\lfloor n/2 \rfloor +1}(\Afra)$.  
However, since $\mathbf{s}_{s-1} < \mathbf{s}_{s}$ we have 
\[
	\U^{\lfloor n/2 \rfloor +1}(\Afra) = G^{s}(F)_{x, \mathbf{s}_{s}+} \subset G^{s}(F)_{x, \mathbf{s}_{s-1}+} \subset K^{d}.  
\]

(2) is similarly shown as (1), using Proposition \ref{compofTLS} (1)-(b), (2)-(a), Lemma \ref{presenofHJ} (2) and Definition \ref{defofKi} (2) dispite of Proposition \ref{compofTLS} (1)-(c), (2)-(b), Lemma \ref{presenofHJ} (1) and Definition \ref{defofKi} (1), respectively.  
Since $J(\beta, \Afra) = {}^{\circ}K^{d}$ and $\Kfra(\Bfra_{0}) = G^{0}(F)_{[x]}$ by Proposition \ref{compofTLS} (1)-(a), we obtain 
\[
	\hat{J}(\beta, \Afra) = \Kfra(\Bfra)J(\beta, \Afra) = G^{0}(F)_{[x]} {}^{\circ}K^{d} = K^{d}, 
\]
whence (3) holds.  
\end{prf}

Let $\theta \in \mathscr{C}(\beta, 0, \Afra)$ be the unique character of $H^{1}(\beta, \Afra)$ in $\lambda$.  
Then we can take characters $\phi_{i}$ of $E_{i}^{\times}$ for $i=0, 1, \ldots, s$ and define characters $\theta^{i}$ as in Proposition \ref{factchar}.  
We put $\Phibf_{i} = \phi_{i} \circ \Nrd_{B_{i}/E_{i}}$.  
If $d=s+1$, we put $\Phibf_{d}=1$.  

\begin{prop}
For $i=0, 1, \ldots, d-1$, the character $\Phibf_{i}$ is $G^{i+1}$-generic relative to $x$ of depth $\rbf_{i}$.  
If $s=d$, then $\Phibf_{d}$ is of depth $\rbf_{d}$.  
\end{prop}

\begin{prf}
First, we show the restriction of $\Phibf_{i}$ to $G^{i}(F)_{x, \mathbf{s}_{i}+}$ is equal to $\psi_{c_{i}}$ for $i=0, \ldots, s$.  
We have 
\begin{eqnarray*}
	B_{i}^{\times} \cap H^{1}(\beta, \Afra) & = & G^{i}(F) \cap K_{+}^{d} \\
	& = & G^{0}(F)_{x, 0+} G^{1}(F)_{x,\mathbf{s}_{0}+} \cdots G^{i}(F)_{x, \mathbf{s}_{i-1}+}, 
\end{eqnarray*}
and $G^{i}(F)_{x, \mathbf{s}_{i}+} \subset B_{i}^{\times} \cap H^{1}(\beta, \Afra)$, as we have $s_{i} > s_{i-1}$ and then
\[
	G^{i}(F)_{x, \mathbf{s}_{i}+} \subset G^{i}(F)_{x, \mathbf{s}_{i-1}+} \subset G^{0}(F)_{x, 0+} G^{1}(F)_{x,\mathbf{s}_{0}+} \cdots G^{i}(F)_{x, \mathbf{s}_{i-1}+}.  
\]
To show $G^{i}(F)_{x,\mathbf{s}_{i}+} \subset H^{t_{i}+1}(\beta, \Afra)$, where $t_{i} = \max \{ 0, \lfloor -v_{\Afra}(c_{i})/2 \rfloor \} = \lfloor -v_{\Afra}(c_{i})/2 \rfloor$, we consider two cases.  
If $i<d$, we have 
\begin{eqnarray*}
	H^{t_{i}+1}(\beta, \Afra) & = & H^{1}(\beta, \Afra) \cap \U^{t_{i}+1}(\Afra) \\
	& = & K_{+}^{d} \cap G(F)_{x, \mathbf{s}_{i}+ } \\
	& = & G^{i+1}(F)_{x, \mathbf{s}_{i}+} \cdots G^{d}(F)_{x, \mathbf{s}_{d-1}+}, 
\end{eqnarray*}
and $G^{i}(F)_{x, \mathbf{s}_{i}+} \subset H^{t_{i}+1}(\beta, \Afra)$ since 
\[
	G^{i}(F)_{x, \mathbf{s}_{i}+} \subset G^{i+1}(F)_{x, \mathbf{s}_{i}+} \subset G^{i+1}(F)_{x, \mathbf{s}_{i}+} \cdots G^{d}(F)_{x, \mathbf{s}_{d-1}+}.  
\]
Otherwise, that is, if $i=s=d$, we also have
\[
	H^{t_{d}+1}(\beta, \Afra) = K_{+}^{d} \cap G(F)_{x, \mathbf{s}_{d}+ } = G^{d}(F)_{x, \mathbf{s}_{d}+}.  
\]
Therefore $G^{i}(F)_{x, \mathbf{s}_{i}+} \subset \left( B_{i}^{\times} \cap H^{1}(\beta, \Afra) \right) \cap H^{t_{i}+1}(\beta, \Afra)$, and we obtain 
\[
	\Phibf_{i} |_{G^{i}(F)_{x, \mathbf{s}_{i}+}} = \theta^{i} |_{G^{i}(F)_{x, \mathbf{s}_{i}+}} \\ = \psi_{c_{i}} |_{G^{i}(F)_{x, \mathbf{s}_{i}+}}.  
\]
In particular, $\Phibf_{i}$ is trivial on 
\[
\U^{-v_{\Afra}(c_{i})+1}(\Afra) \cap G^{i}(F)_{x, \mathbf{s}_{i}+} = G(F)_{x, \rbf_{i}+} \cap G^{i}(F)_{x, \mathbf{s}_{i}+} = G^{i}(F)_{x, \rbf_{i}+}.  
\]

Next, we show $\Phibf_{i}$ is not trivial on $G^{i}(F)_{x,\rbf_{i}}$.  
We have $G^{i}(F)_{x,\rbf_{i}} =\U(\Bfra_{i}) \cap \U^{-v_{\Afra}(c_{i})}(\Afra) = B_{i} \cap \left( 1+\Pfra^{-v_{\Afra}(c_{i})} \right) = 1 + \left( B_{i} \cap \Pfra^{-v_{\Afra}(c_{i})} \right)$.  
Then 
\begin{eqnarray*}
	\Phibf_{i}(G^{i}(F)_{x,\rbf_{i}}) & = & \psi_{c_{i}} \Bigl( 1 + \bigl( B_{i} \cap \Pfra^{-v_{\Afra}(c_{i})} \bigr) \Bigr) = \psi \circ \Trd_{A/F} \Bigl( c_{i} \bigl( B_{i} \cap \Pfra^{-v_{\Afra}(c_{i})} \bigr) \Bigr) \\
	& = & \psi \circ \Trd_{A/F} \left( B_{i} \cap \Afra \right) = \psi \circ \Tr_{E_{i}/F} \circ \Trd_{B_{i}/E_{i}} (\Bfra_{i}).  
\end{eqnarray*}
Since $\Bfra_{i}$ is a hereditary $\ofra_{E_{i}}$-order in $B_{i}$, we have $\Trd_{B_{i}/E_{i}} (\Bfra_{i}) = \ofra_{E_{i}}$.  
Moreover, since $E_{i}/F$ is tamely ramified, $\Tr_{E_{i}/F} (\ofra_{E_{i}}) = \ofra_{F}$.  
Therefore $\Phibf_{i}$ is not trivial on $G^{i}(F)_{x,\rbf_{i}}$, as $\psi$ is not trivial on $\ofra_{F}$.  
In particular, we completed the proof when $i=s=d$ and we may assume $i<d$ in the following.  

Finally, let $X_{c_{i}}^{*} \in \Lie^{*}(Z(G^{i}))$ as \S \ref{defofXc}.  
Since $c_{i}$ is minimal relative to $E_{i}/E_{i+1}$, the element $X_{c_{i}}^{*}$ is in $\Lie^{*}(Z(G^{i}))_{-\rbf_{i}}$ and $G^{i+1}$-generic of depth $\rbf_{i}$ by Proposition \ref{genelem}.  
Then, to complete the proof it suffices to show that $\Phibf_{i} |_{G^{i}(F)_{x,\rbf_{i}:\rbf_{i}+}}$ is realized by $X_{c_{i}}^{*}$.  
The isomorphism $G^{i}(F)_{x,\rbf_{i}:\rbf_{i}+} \cong \gfra^{i}(F)_{x,\rbf_{i}:\rbf_{i}+}$ is induced from $1+y \mapsto y$.  
Therefore, when we regard $\psi \circ X_{c_{i}}^{*}$ as a character of $G^{i}(F)_{x,\rbf_{i}:\rbf_{i}+}$, for $1+y \in G^{i}(F)_{x, \rbf_{i}}$ we have
\[
	\left( \psi \circ X_{c_{i}}^{*} \right) (1+y) = \psi \circ X_{c_{i}}^{*}(y) = \psi \circ \Trd_{A/F}(c_{i}y) = \psi_{c_{i}}(1+y) = \Phibf_{i}(1+y).  
\]
\end{prf}

Then we have a 4-tuple $\left( x, (G^{i}), (\rbf_{i}), (\Phibf_{i}) \right)$.  
As in \S \ref{Yuconst}, we can define characters $\hat{\Phibf}_{i}$ of $K_{+}^{d}$.  

\begin{prop}
\label{phiistheta1}
For $i=0, 1, \ldots, s$, we have $\hat{\Phibf}_{i}=\theta^{i}$.  
\end{prop}

\begin{prf}
Recall the definition of $\hat{\Phibf}_{i}$.  The character $\hat{\Phibf}_{i}$ is defined as 
\begin{eqnarray*}
\hat{\Phibf}_{i}|_{K_{+}^{d} \cap G^{i}(F)} (g) & = & \Phibf_{i}(g), \\
\hat{\Phibf}_{i}|_{K_{+}^{d} \cap G(F)_{x, \mathbf{s}_{i}+ } } (1+y) & = & \Phibf_{i} (1+\pi_{i}(y) ).  
\end{eqnarray*}

Since $\left( K_{+}^{d} \cap G^{i}(F) \right) \left( K_{+}^{d} \cap G(F)_{x, \mathbf{s}_{i}+ } \right) =K^{d}_{+}$, it is enough to show that $\hat{\Phibf}_{i}$ is equal to $\theta^{i}$ on $K_{+}^{d} \cap G^{i}(F)$ and $K_{+}^{d} \cap G(F)_{x, \mathbf{s}_{i}+ }$.  

We have that $K^{d}_{+} \cap G^{i}(F) = B_{i}^{\times} \cap H^{1}(\beta, \Afra)$ and $K^{d}_{+} \cap G(F)_{r, \mathbf{s}_{i}+}=H^{t_{i}+1}(\beta, \Afra)$, where $t_{i} = \lfloor -v_{\Afra}(c_{i})/2 \rfloor$.  

If $g \in B_{i}^{\times} \cap H^{1}(\beta, \Afra)$, then $\hat{\Phibf}_{i}(g)=\Phibf_{i}(g)=\phi_{i} \circ \Nrd_{B_{i}/E_{i}}(g) = \theta^{i}(g)$.  

Suppose $1+y \in H^{t_{i}+1}(\beta, \Afra)$.  
Then $\pi_{i}(y) \in \mathfrak{g}^{i}(F)_{x, \mathbf{s}_{i}+}=B \cap \Pfra_{i}^{t_{i}+1}$ and $1+\pi_{i}(y) \in B_{i}^{\times} \cap H^{t_{i}+1}(\beta, \Afra)$.  
Therefore we have $\hat{\Phibf}_{i}(1+y)=\Phibf_{i}(1+\pi_{i}(y))=\theta^{i}(1+\pi_{i}(y))=\psi \circ \Trd_{A/F}(c_{i} \pi_{i}(y))$.  

Here, we show if $n \in \mathfrak{n}^{i}(F)$, then $\Trd_{A/F}(c_{i}n)=0$.  
Since $c_{i}$ is in the center of $B_{i}^{\times}$, the linear automorphism $z \mapsto c_{i}z$ of $A$ is also a $Z(G^{i})(F)$-automorphism.  
Then $c_{i}\gfra^{i}(F)$ is a trivial $Z(G^{i})(F)$-representation and $c_{i}\mathfrak{n}^{i}(F) \cong \mathfrak{n}^{i}(F)$ is a $Z(G^{i})(F)$-representation which does not contain any trivial representation.  
Therefore we have $c_{i}\gfra^{i}(F) \subset \gfra^{i}(F)$ and $c_{i} \mathfrak{n}^{i}(F) \subset \mathfrak{n}^{i}(F)$.  
On the other hand, $\Trd_{A/F}$ is a $Z(G^{i})(F)$-homomorphism from $\gfra(F)$ to the trivial representation $F$.  
Since $\mathfrak{n}^{i}(F)$ does not contain any trivial representations, $\Trd_{A/F}(\mathfrak{n}^{i}(F))=0$.  
In particular, $\Trd_{A/F}(c_{i}n)=0$ as $c_{i}n \in c_{i} \mathfrak{n}^{i}(F) \subset \mathfrak{n}^{i}(F)$.  

Since $\pi_{i}: \gfra^{i}(F) \oplus \mathfrak{n}^{i}(F) \to \gfra^{i}(F)$ is the projection, $y-\pi_{i}(y) \in \mathfrak{n}^{i}(F)$.  
Therefore we have $\Trd_{A/F}( c_{i}y)=\Trd_{A/F} \left( c_{i} \left( y-\pi_{i}(y) \right) \right) + \Trd_{A/F} \left( c_{i}\pi_{i}(y) \right) = \Trd_{A/F} \left( c_{i}\pi_{i}(y) \right)$ and
\[
	\hat{\Phibf}_{i}(1+y) = \psi \circ \Trd_{A/F}(c_{i} \pi_{i}(y)) = \psi \circ \Trd_{A/F}(c_{i}y) = \psi_{c_i}(1+y) = \theta^{i}(1+y).  
\]
\end{prf}

\begin{prop}
\label{betaext}
The representation $\kappa_{0} \otimes \cdots \otimes \kappa_{d}$ is an extension of $\eta_\theta$ to $K^{d}$.  
\end{prop}

\begin{prf}
We put $\hat{\kappa}'=\kappa_{0} \otimes \cdots \otimes \kappa_{d}$.  
By \cite[Lemma 3.27]{HM}, $\kappa_{i} |_{K_{+}^{d}}$ contains $\hat{\Phibf}_{i}$ for $i=0, \ldots, d$.  
If $d=s+1$, then $\hat{\Phibf}_{d}=\Phibf_{d}=1$ and 
\[
	\prod_{i=0}^{d} \hat{\Phibf}_{i} = \prod_{i=0}^{s} \hat{\Phibf}_{i} = \prod_{i=0}^{s} \theta^{i} = \theta
\]
by Proposition \ref{factchar} and Proposition \ref{phiistheta1}.  
Then $\hat{\kappa}'$ contains $\theta$ as a $K_{+}^{d}$-representation.  
Since $\eta_{\theta}$ is the unique irreducible $J^{1}(\beta, \Afra)$-representation which contains $\theta$, the $J^{1}(\beta, \Afra)$-representation $\hat{\kappa}'$ contains $\eta_{\theta}$.  

Then it suffices to show that the dimension of $\hat{\kappa}$ is equal to the dimension of $\eta_{\theta}$.  
The dimension of $\eta_{\theta}$ is $(J^{1}(\beta, \Afra) : H^{1}(\beta, \Afra))^{1/2}$.  
On the other hand, for $i=0, \ldots, d-1$ the dimension of $\kappa_{i}$ is $(J^{i} : J_{+}^{i})^{1/2}$, and the dimension of $\kappa_{d}$ is 1.  
Then the dimension of $\hat{\kappa}'$ is $\prod_{i=1}^{d} (J^{i} : J_{+}^{i})^{1/2}$, and it suffices to show that $(J^{1}(\beta, \Afra) : H^{1}(\beta, \Afra)) = \prod_{i=1}^{d} (J^{i} : J_{+}^{i})$.  
Here, $H^{1}(\beta, \Afra) = K_{+}^{d} = K_{+}^{0}J_{+}^{1} \cdots J_{+}^{d} = G^{0}(F)_{x,0+}J_{+}^{1} \cdots J_{+}^{d}$.  
Since $G^{i}(F)_{x,\mathbf{s}_{i}}J^{i+1}=G^{i+1}(F)_{x,\mathbf{s}_{i}}$ for $i=0, \ldots, d-1$, we also have 
\begin{eqnarray*}
	G^{0}(F)_{x,0+}J^{1} \cdots J^{d} & = & G^{0}(F)_{x,0+}G^{0}(F)_{x,\mathbf{s}_{0}}J^{1} \cdots J^{d} \\
	& = & G^{0}(F)_{x,0+}G^{1}(F)_{x,\mathbf{s}_{0}} J^{2} \cdots J^{d} = \cdots \\
	& = & G^{0}(F)_{x,0+}G^{1}(F)_{x, \mathbf{s}_{0}} \cdots G^{d}(F)_{x, \mathbf{s}_{d-1}} \\
	& = & G(F)_{x,0+} \cap K^{d}=\U^{1}(\Afra) \cap J(\beta, \Afra) = J^{1}(\beta, \Afra).  
\end{eqnarray*}
Since $G^{0}(F)_{x,0+} \cap (J^{1} \cdots J^{d}) = G^{0}(F)_{x,\rbf_{0}} \subset J_{+}^{1}$, we have $\left( G^{0}(F)_{x,0+} J_{+}^{1} \cdots J_{+}^{d} \right) \cap ( J^{1} \cdots J^{d} ) = J_{+}^{1} \cdots J_{+}^{d}$, and
\begin{eqnarray*}
	J^{1}(\beta, \Afra)/H^{1}(\beta, \Afra) & = & \left( G^{0}(F)_{x,0+}J^{1} \cdots J^{d} \right) / \left( G^{0}(F)_{x,0+} J_{+}^{1} \cdots J_{+}^{d} \right) \\
	& \cong & ( J^{1} \cdots J^{d} ) / ( J_{+}^{1} \cdots J_{+}^{d} ).  
\end{eqnarray*}
Then it is enough to show $\left( ( J^{1} \cdots J^{d} ) : ( J_{+}^{1} \cdots J_{+}^{d} ) \right) = \prod_{i=1}^{d} (J^{i} : J_{+}^{i})$, which is already proved in part (c) in the proof of \cite[Proposition 9.2]{May}.  
Therefore we obtain $\hat{\kappa}' |_{J^{1}(\beta, \Afra)} = \eta_{\theta}$.  
\end{prf}

\begin{thm}
\label{Main1}
Let $(J, \lambda)$ be a maximal simple type associated to a tame simple stratum $[\Afra, n, 0, \beta]$.  
Let $(\tilde{J}, \Lambda)$ be a maximal extension of $(J, \lambda)$.  
Then there exists a Yu datum $\left( x, (G^{i})_{i=0}^{d}, (\rbf_{i})_{i=0}^{d}, (\Phibf_{i})_{i=0}^{d}, \rho \right)$ such that
\begin{enumerate}
\item $\hat{J}(\beta, \Afra)=K^{d}$, and
\item $\rho_{d} \left( x, (G^{i}), (\rbf_{i}), (\Phibf_{i}), \rho \right) \cong \cInd_{\tilde{J}}^{\hat{J}(\beta, \Afra)} \Lambda$.  
\end{enumerate}
\end{thm}

\begin{prf}
In the above argument, we can take a 4-tuple $\left( x, (G^{i}), (\rbf_{i}), (\Phibf_{i}) \right)$ from a S\'echerre--Stevens datum.  
Therefore it is enough to show that we can take an irreducible $G^{0}(F)_{[x]}$--representation $\rho$ such that the Yu datum $\left( x, (G^{i}), (\rbf_{i}), (\Phibf_{i}), \rho \right)$ satisfies the desired conditions.  

Let $\eta$ be the unique $J^{1}(\beta, \Afra)$-subrepresentation in $\lambda|_{J^{1}(\beta, \Afra)}$.  
Then $\kappa_{0} \otimes \cdots \otimes \kappa_{d}$ is an extension of $\eta$ to $K^{d}=\hat{J}(\beta, \Afra)$ by Proposition \ref{betaext}.  
Therefore there exists an irreducible $\Kfra(\Bfra)$-representation $\rho$ such that $\rho$ is trivial on $\U^{1}(\Bfra)$ but not trivial on $\U(\Bfra)$, the representation $\cInd_{\Kfra(\Bfra)}^{B^{\times}}\rho$ is irreducible and supercuspidal, and $\cInd_{\tilde{J}}^{\hat{J}(\beta, \Afra)} \Lambda \cong \rho \otimes \kappa_{0} \otimes \cdots \otimes \kappa_{d}$ by Proposition \ref{findingdepth0}.  
Since we have equalities of groups $B^{\times}=G^{0}(F)$, $\Kfra(\Bfra)=G^{0}(F)_{[x]}$, $\U(\Bfra)=G^{0}(F)_{x}$ and $\U^{1}(\Bfra) = G^{0}(F)_{x+}$, then the 5-tuple $\left( x, (G^{i}), (\rbf_{i}), (\Phibf_{i}), \rho \right)$ is a Yu datum satisfying the condition in the theorem.  
\end{prf}

\section{Construction of a S\'echerre--Stevens datum from a Yu datum}
\label{YutoSS}

Let $(x, (G^{i})_{i=0}^{d}, (\rbf_{i})_{i=0}^{d}, (\Phibf_{i})_{i=0}^{d}, \rho)$ be a Yu datum.  

First, since $G^{i}$ are tame twisted Levi subgroups in $G$ with $Z(G^{i})/Z(G)$ anisotropic, there exist tamely ramified field extensions $E_{i}/F$ in $A$ such that 
\[
	G^{i} \cong \Res_{E_{i}/F} \underline{\Aut}_{D \otimes_{F} E_{i}} (V)
\]
by Lemma \ref{detofTLG}.  
Since $G^{0} \subsetneq \ldots \subsetneq G^{d}$, we have $E_{0} \supsetneq \ldots \supsetneq E_{d} = F$.  
We put $B_{i} = \Cent_{A}(E_{i})$.  

Since $\cInd_{G^{0}(F)_{[x]}}^{G^{0}(F)} \rho$ is supercuspidal, $[x]$ is a vertex in $\Bscr^{E}(G^{0},F)$ by Proposition \ref{depth0ofYu}.  
Let $\Bfra_{0}$ be the hereditary $\ofra_{E_{0}}$-order in $B_{0}$ associated with $x$.  
Then the hereditary $\ofra_{F}$-order $\Afra$ associated with $x \in \Bscr^{E}(G,F)$ is $E_{0}$-pure and principal, and $\Afra \cap B_{0} = \Bfra_{0}$ by Proposition \ref{beingprincipal}.  
We also put $\Pfra = \Pfra(\Afra)$.  

To obtain a simple stratum, we need an element $\beta \in E_{0}$.  
We will take $\beta$ by using information from characters $(\Phibf_{i})_{i}$.  
For $c_{i} \in E_{i} = \Lie(Z(G^{i}))$, let $X_{c_{i}}^{*} \in \Lie^{*} (Z(G^{i}))$ be as in \S \ref{defofXc}.  
We put $s=\sup \{ i \mid \Phibf_i \neq 1 \}$.  

\begin{prop}
\label{takec}
Suppose $s \geq 0$.  
\begin{enumerate}
\item For $i=0, \ldots, d$, the hereditary $\ofra_{E_{i}}$-order in $B_{i}$ associated with $x \in \Bscr^{E}(G^{i}, F)$ is equal to $\Bfra_{i} = B_{i} \cap \Afra$.  
\item There exists $c_{i} \in \Lie(Z(G^{i}))_{-\rbf_{i}}$ such that $\Phibf_{i} | _{G^{i}(F)_{x,\rbf_{i}/2+:\rbf_{i}+}}$ is realized by $X_{c_i}^{*}$ for $i=0, \ldots, d-1$.  
\item If $s=d$, then there also exists $c_{s} \in \Lie(Z(G))_{-\rbf_{s}}$ such that $\Phibf_{s} | _{G(F)_{x,\rbf_{s}/2+:\rbf_{s}+}}$ is realized by $X_{c_s}^{*}$.  
\item For $i=0, \ldots, s$, we have $\rbf_{i} = -\ord(c_{i})$.  
\item For $i=0, \ldots, d-1$, the element $c_{i}$ is minimal relative to $E_{i}/E_{i+1}$.  
In particular, we have $E_{i} = E_{i+1}[c_{i}]$.  
\end{enumerate}
\end{prop}

\begin{prf}
We show (1).  
First, we have $\Bfra_{i} \cap B_{0} = \Afra \cap B_{i} \cap B_{0} = \Afra \cap B_{0} = \Bfra_{0}$.  
Moreover, for $g \in E_{i}^{\times}$ we also have 
\[
	g \Bfra_{i} g^{-1} = g (\Afra \cap B_{i}) g^{-1} = g \Afra g^{-1} \cap g B_{i} g^{-1} = \Afra \cap B_{i} = \Bfra_{i}, 
\]
as $\Afra$ is $E_{0}$-pure and $E_{0} \subset B_{i}$.  
Therefore (1) holds by Proposition \ref{beingprincipal} (2).  

Next, we show (2), and (3) is similarly shown.  
Since $\Phibf_{i}$ is trivial on $G^{i}(F)_{x, \rbf_{i}+}$ but not on $G^{i}(F)_{x, \rbf_{i}}$, we have $G^{i}(F)_{x, \rbf_{i}} \neq G^{i}(F)_{x, \rbf_{i}+}$ in particular.  
Then $n_{i} = \rbf_{i}e(\Bfra_{i}|\ofra_{E_{i}})e(E_{i}/F)$ is a non-negative integer and we have $G^{i}(F)_{x, \rbf_{i}} = \U^{n}(\Bfra)$ and $G^{i}(F)_{x, \rbf_{i}+} = \U^{n+1}(\Bfra)$, by Lemma \ref{lemforfindc} (3).  
On the other hand, a character $\psi \circ \Tr_{E_{i}/F}$ of $E_{i}$ is with conductor $\pfra_{E_{i}}$ since $E_{i}/F$ is tamely ramified.  
Therefore, we can apply Proposition \ref{propforgench} for $\Bfra_{i}, n$ and $\psi \circ \Tr_{E_{i}/F}$ as $\Bfra_{i}$ is principal by (1) and Proposition \ref{beingprincipal} (1).  
Thus there exists $c_{i} \in E_{i}$ such that 
\[
	\Phibf_{i} (1+y) = (\psi \circ \Tr_{E_{i}/F}) \circ \Trd_{B_{i}/E_{i}}(c_{i}y) = \psi \circ (\Tr_{E_{i}/F} \circ \Trd_{B_{i}/E_{i}}) (c_{i}y) = \psi \circ X_{c_{i}}^{*}(y)
\]
for $1+y \in \U^{\lfloor n_{i}/2 \rfloor +1}(\Bfra_{i}) = G^{i}(F)_{x, \mathbf{r}_{i}/2+}$.  
Then (2) holds.  

We have $v_{E_{i}}(c_{i})=-n_{i}/e(\Bfra_{i}|\ofra_{E_{i}})=-\rbf_{i}e(E_{i}/F)$ by Proposition \ref{propforgench}, and
\[
\ord(c_{i})=v_{E_{i}}(c_{i})/e(E_{i}/F)=-\rbf_{i}, 
\]
whence (4) holds.  

To show (5), let $c'_{i} \in E_{i}^{\times}$ such that $X_{c'_{i}}^{*}$ is $G^{i+1}$-generic of $\rbf_{i}$ and $G^{i}(F)_{x, \rbf_{i}:\rbf_{i}+}$ is realized by $X_{c'_{i}}$.  
In particular, we have
\[
	(\psi \circ \Tr_{E_{i}/F}) \circ \Trd_{B_{i}/E_{i}}(c_{i}y) = \Phibf_{i}(1+y) = \psi \circ X_{c'_{i}}(y) = (\psi \circ \Tr_{E_{i}/F}) \circ \Trd_{B_{i}/E_{i}}(c'_{i}y)
\]
for $y \in \Qfra_{i}^{n_{i}}$, where $\Qfra_{i}$ is the radical of $\Bfra_{i}$.  
Then we have $c_{i} - c'_{i} \in \Qfra_{i}^{-n_{i}+1} \cap E_{i} \subset c_{i}(\Qfra_{i} \cap E_{i}) = c_{i}\pfra_{E_{i}}$ and $c_{i}^{-1}c'_{i} \in 1+\pfra_{E_{i}}$.  
Thus $(c'_{i})^{-1}c_{i} \in 1+\pfra_{E_{i}}$.  
On the other hand, $c'_{i}$ is minimal relative to $E_{i}/E_{i+1}$ by Proposition \ref{genelem}.  
Therefore, by Lemma \ref{preservemin} $c_{i}$ is also minimal relative to $E_{i}/E_{i+1}$.  
\end{prf}

Therefore if $s \geq 0$, we can take $c_i$ for $i=0,1, \ldots, s$.  
We put $\beta_{i} = \sum_{j=i}^{s} c_{j}$ for $i=0,1, \ldots, s$, $\beta = \beta_{0}$ and $n = -v_{\Afra}(\beta)$.  
Since 
\[
v_{\Afra}(c_i) = -e(\Afra|\ofra_{F}) \ord(c_{i}) = -e(\Afra|\ofra_{F})\rbf_{i} < -e(\Afra|\ofra_{F})\rbf_{j} = -e(\Afra|\ofra_{F}) \ord(c_{j}) = v_{\Afra}(c_j)
\]
for $i,j=0, 1, \ldots, s$ with $i>j$, we have $n=-v_{\Afra}(\beta_{i})$ for any $i=0, 1, \ldots, s$.  
We also put $r_{i} = -v_{\Afra}(c_{i-1})$ for $i=1, \ldots, s$ and $r_{0}=0$.  

\begin{prop}
Suppose $s \geq 0$.  
\begin{enumerate}
\item $E_{i}=F[\beta_{i}]$ for $i=0, 1, \ldots, s$.  
	In particular, $[\Afra, n, 0, \beta]$ is a simple stratum.  
\item $\left( [\Afra, n, r_{i}, \beta_{i}] \right) _{i=0}^{s}$ is a defining sequence of $[\Afra, n, 0, \beta]$.  
\end{enumerate}
\end{prop}

\begin{prf}
If $s=-\infty$, nothing has to be shown.  
Then we may assume $s \in \Z$.  

First, suppose $\Afra = \Afra(E_0)$.  
We will show this proposition by downward induction on $i$.  

If $i=s$, then $\beta_{s}=c_{s}$ is minimal over $F$.  
Therefore for any $r' \in \{ 0, 1, \ldots, n-1 \}$, the stratum $[\Afra, n, r', \beta_{s}]$ is simple.  
The equation $E_{s}=F[\beta_{s}]$ trivially holds.  
If $s=0$, then $\left( [\Afra, n, r_{i}, \beta_{i}] \right) _{i=0}^{0}$ is a defining sequence of $[\Afra, n, 0, \beta]$ and this proposition holds.  
If $s>0$, we have $r_{s} = -v_{\Afra}(c_{s-1}) < -v_{\Afra}(c_{s})$.  
Then $[\Afra, n, r_{s}, \beta_{s}]$ is simple and $\left( [\Afra, n, r_{i+s}, \beta_{i+s}] \right)_{i=0}^{0}$ is a defining sequence of $[\Afra, n, r_{s}, \beta_{s}]$.  

Let $i_{0} \in \{ 0, 1, \ldots, s-1 \} $ and suppose that $E_{i}=F[\beta_{i}]$ and that $\left( [\Afra, n, r_{j+i}, \beta_{j+i}] \right)_{j=0}^{s-i}$ is a defining sequence of a simple stratum $[\Afra, n, r_{i}, \beta_{i}]$ for any integer $i$ with $i_0 < i  \leq s$.  
The element $c_{i_{0}}$ is minimal over $E_{i_{0}+1}$.  
Since $r_{i_{0}+1}=-v_{\Afra}(c_{i_0})$, a 4-tuple $[\Bfra_{\beta_{i_0}+1}, r_{i_{0}+1}, r_{i_{0}+1}-1, c_{i_0}]$ is a simple stratum, where $\Bfra_{\beta_{i_{0}+1}} = \Afra \cap \Cent_{A(E_0)}(\beta_{i_{0}+1})$.  
Moreover, $c_{i_0} \notin E_{i_{0}+1}=F[\beta_{i_{0}+1}]$.  
Therefore, by Proposition \ref{Mforcons}, we have $F[\beta_{i_0}] = F[\beta_{i_{0}+1}, c_{i_0}] = E_{i_{0}+1}[c_{i_{0}}]$ and $[\Afra, n, r_{i_{0}+1}, \beta_{i_0}]$ is a pure stratum with $k_{0}(\beta_{i_0}, \Afra)=-r_{i_{0}+1}$, where $F[\beta_{i_{0}+1}, c_{i_0}] = E_{i_{0}+1}[c_{i_{0}}]$ follows from our induction hypothesis.  
If $i_0>0$, we have $r_{i_0}=-v_{\Afra}(c_{i_{0}-1}) < -v_{\Afra}(c_{i_{0}})=r_{i_{0}+1}$ and $[\Afra, n, r_{i_{0}}, \beta_{i_{0}}]$ is a simple stratum.  
Since $\left( [\Afra, n, r_{j+i_{0}+1}, \beta_{j+i_{0}+1}] \right)_{j=0}^{s-i_{0}-1}$ is a defining sequence of a simple stratum $[\Afra, n, r_{i_{0}+1}, \beta_{i_{0}+1}]$ by our induction hypothesis, $\left( [\Afra, n, r_{j+i_{0}}, \beta_{j+i_{0}}] \right) _{j=0}^{s-i_{0}}$ is also a defining sequence of a simple stratum $[\Afra, n, r_{i_{0}}, \beta_{i_{0}}]$.  
If $i_{0}=0$, then $[\Afra, n, 0, \beta]$ is simple and we can show $\left( [\Afra, n, r_{i}, \beta_{i}] \right) _{i=0}^{s}$ is also a defining sequence of a simple stratum $[\Afra, n, 0, \beta]$ in the same way as above.  
Then the proposition for $\Afra=\Afra(E_{0})$ case holds.  

We will show the proposition in general case.  
Since $\beta_{i} \in E_{i} \subset E_{0}$ for $i=0, \ldots, s$, we can regard $\beta_{i}$ as in $A(E_0)$.  
Then (1) follows from the proposition for $\Afra=\Afra(E_0)$ case.  
Moreover, if we put $n'=-v_{\Afra(E_0)}(\beta)$, $r'_{0}=0$ and $r'_{i}=-v_{\Afra(E_0)}(c_{i-1})$ for $i=1, \ldots, s$, then $\left( [\Afra(E_0), n', r'_{i}, \beta_{i}] \right) _{i=0}^{s}$ is a defining sequence of a simple type $[\Afra(E_0), n', 0, \beta]$ by the proposition for $\Afra=\Afra(E_0)$ case.  
Since for $c \in E_0$ we have $v_{\Afra}(c)=e(\Afra|\ofra_{F})e(E_{0}/F)^{-1}v_{\Afra(E_0)}(c)$, we also have
\[
	n=-v_{\Afra}(\beta)=-e(\Afra|\ofra_{F})e(E_{0}/F)^{-1}v_{\Afra(E_0)}(\beta)=e(\Afra|\ofra_{F})e(E_{0}/F)^{-1}n'
\]
and
\[
r_{i}=-v_{\Afra}(c_{i-1})=-e(\Afra|\ofra_{F})e(E_{0}/F)^{-1}v_{\Afra(E_0)}(c_{i-1})=e(\Afra|\ofra_{F})e(E_{0}/F)^{-1}r'_{i}
\]
for $i=1, \ldots, s$.  
Since $\left( [\Afra(E_0), n', r'_{i}, \beta_{i}] \right) _{i=0}^{s}$ is a defining sequence of a simple type $[\Afra(E_0), n', 0, \beta]$, we have $r'_{i}=-k_{0}(\beta_{i-1}, \Afra(E_0))$ for $i=1, \ldots, s$.  
We also have $k_{0}(c, \Afra)=e(\Afra|\ofra_{F})e(E_{0}/F)^{-1}k_{0}(c, \Afra(E_0))$ by Lemma \ref{compofk0}, whence
\[
r_{i}=e(\Afra|\ofra_{F})e(E_{0}/F)^{-1}r'_{i}=-e(\Afra|\ofra_{F})e(E_{0}/F)^{-1}k_{0}(\beta_{i-1}, \Afra(E_0))=-k_{0}(\beta_{i-1}, \Afra)
\]
for $i=1, \ldots, s$.  Then by Proposition \ref{BHforapp} strata $[\Afra, n, r_{i}, \beta_{i}]$ are simple and equivalent to $[\Afra, n, r_{i}, \beta_{i-1}]$ for $i=1, \ldots, s$.  
Therefore (2) holds.  
\end{prf}

Then we have a simple stratum $[\Afra, n, 0, \beta]$ with a defining sequence $([\Afra, n, r_{i}, \beta_{i}])_{i=0}^{s}$ if $s \geq 0$.  
If $s = - \infty$, we take a simple stratum $[\Afra, 0, 0, \beta]$ with $\Afra$ maximal and $c_{0} = \beta_{0} = \beta \in \ofra_{F}^{\times}$, and then we can define subgroups $H^{1}(\beta, \Afra)$ and $J(\beta, \Afra)$ in $G$ for any case.  
Moreover, since $\Bfra_{0}$ is maximal, we also can define $\hat{J}(\beta, \Afra) = \Kfra(\Bfra_{0})J(\beta, \Afra)$.  

\begin{prop}
\begin{enumerate}
\item We have
\begin{enumerate}
\item $G^{0}(F)_{[x]} = \Kfra(\Bfra_{0})$, 
\item $G^{0}(F)_{x} = B_{0}^{\times} \cap \U(\Afra) = \U(\Bfra_{0})$, 
\item $G^{0}(F)_{x,0+} = B_{0}^{\times} \cap \U^{1}(\Afra)$, 
\item $\gfra^{0}(F)_{x} = B_{0} \cap \Afra = \Bfra_{0}$, and
\item $\gfra^{0}(F)_{x,0+} = B_{0} \cap \Pfra$.  
\end{enumerate}
\item For $i=1, \ldots, d$, we have
\begin{enumerate}
\item $G^{i}(F)_{x, \mathbf{s}_{i-1} } = B_{i}^{\times} \cap \U^{\lfloor (-v_{\Afra}(c_{i-1})+1)/2 \rfloor}(\Afra)$, 
\item $G^{i}(F)_{x, \mathbf{s}_{i-1}+ } = B_{i}^{\times} \cap \U^{\lfloor -v_{\Afra}(c_{i-1})/2 \rfloor +1}(\Afra)$, 
\item $G^{i}(F)_{x, \rbf_{i-1} } = B_{i}^{\times} \cap \U^{-v_{\Afra}(c_{i-1})}(\Afra)$, 
\item $G^{i}(F)_{x, \rbf_{i-1}+ } = B_{i}^{\times} \cap \U^{-v_{\Afra}(c_{i-1})+1}(\Afra)$, 
\item $\gfra^{i}(F)_{x, \mathbf{s}_{i-1} } = B_{i} \cap \Pfra^{\lfloor (-v_{\Afra}(c_{i-1})+1)/2 \rfloor}$, 
\item $\gfra^{i}(F)_{x, \mathbf{s}_{i-1}+ } = B_{i} \cap \Pfra^{\lfloor -v_{\Afra}(c_{i-1})/2 \rfloor +1}$, 
\item $\gfra^{i}(F)_{x, \rbf_{i-1} } = B_{i} \cap \Pfra^{-v_{\Afra}(c_{i-1})}$, and
\item $\gfra^{i}(F)_{x, \rbf_{i-1}+ } = B_{i} \cap \Pfra^{-v_{\Afra}(c_{i-1})+1}$.  
\end{enumerate}
\item For $i=0, \ldots, s$, we have
\begin{enumerate}
\item $G^{i}(F)_{x, \mathbf{s}_{i}+ } = B_{i}^{\times} \cap \U^{\lfloor -v_{\Afra}(c_{i})/2 \rfloor +1}(\Afra)$, 
\item $G^{i}(F)_{x, \rbf_{i} }  = B_{i}^{\times} \cap \U^{-v_{\Afra}(c_{i})}(\Afra)$,
\item $G^{i}(F)_{x, \rbf_{i}+ } = B_{i}^{\times} \cap \U^{-v_{\Afra}(c_{i})+1} (\Afra)$, 
\item $\gfra^{i}(F)_{x, \rbf_{i} } = B \cap \Pfra^{-v_{\Afra}(c_{i})}$, and
\item $\gfra^{i}(F)_{x, \rbf_{i}+ } = B \cap \Pfra^{-v_{\Afra}(c_{i})}$.  
\end{enumerate}
\end{enumerate}
\end{prop}

\begin{prf}
Similar to the proof of Proposition \ref{compofTLS}.  
\end{prf}

\begin{prop}
\begin{enumerate}
\item $K^{d}_{+} = H^{1}(\beta, \Afra)$, 
\item ${}^{\circ}K^{d} = J(\beta, \Afra)$, and
\item $K^{d} = \hat{J}(\beta, \Afra)$.  
\end{enumerate}
\end{prop}

\begin{prf}
Similar to the proof of Proposition \ref{compofHJ}.  
\end{prf}

Next, we construct a simple character in $\mathscr{C}(\beta, 0, \Afra)$ from $(\Phibf_{i})_{i}$.  

\begin{lem}
\label{lemforbeingsimp}
Suppose $s \geq 0$.  
For $i=0, 1, \ldots, s$, the following assertions hold.  
\begin{enumerate}
\item $\hat{\Phibf}_{i} |_{B_{i}^{\times} \cap H^{1}(\beta, \Afra)}$ factors through $\Nrd_{B_{i}/E_{i}}$.  
\item $\hat{\Phibf}_{i} |_{H^{t_{i}+1}(\beta, \Afra)} = \psi_{c_{i}}$, where $t_{i} = \lfloor -v_{\Afra}(c_{i})/2 \rfloor$.  
\item $H^{t_{i}+1}(\beta, \Afra) = H^{t_{i}+1}(\beta_{i}, \Afra)$ is normalized by $B_{i}^{\times} \cap \Kfra(\Afra)$.  
\item For any $g \in B_{i}^{\times} \cap \Kfra(\Afra)$ and $h \in H^{1}(\beta, \Afra) \cap {}^{g}H^{1}(\beta, \Afra)$, we have $\hat{\Phibf}_{i}(g^{-1}hg) = \hat{\Phibf}_{i}(h)$.  
\end{enumerate}
\end{lem}

\begin{prf}
We have $B_{i}^{\times} \cap H^{1}(\beta, \Afra) = G^{i}(F) \cap K^{d}$.  
By construction of $\hat{\Phibf}_{i}$ we have $\hat{\Phibf}_{i}|_{B_{i}^{\times} \cap H^{1}(\beta, \Afra)} = \hat{\Phibf}_{i}|_{G^{i}(F) \cap K^{d}} = \Phibf_{i}$.  
The map $\Phibf_{i}$ is a character of $G^{i}(F)$, and then $\Phibf_{i}$ factors through $\Nrd_{B_{i}/E_{i}}$ and (1) holds.  

We also have $H^{t_{i}+1}(\beta, \Afra) = K^{d} \cap G(F)_{x, \mathbf{s}_{i}+}$.  
Since $\Phibf_{i}|_{G^{i}(F)_{x,\mathbf{s}_{i}+:\rbf_{i}+}}$ is realized by $X_{c_{i}}^{*}$ by Proposition \ref{takec} (2) or (3), we have 
\[
	\Phibf_{i}(1+y) = \psi \circ \Tr_{E_{i}/F} \circ \Trd_{B_{i}/E_{i}} (c_{i}y) = \psi \circ \Trd_{A/F}(c_{i}y)
\]
for $y \in B_{i} \cap \Pfra^{t_{i}+1}=\gfra^{i}(F)_{x, \mathbf{s}_{i}+}$.  
We recall that $\pi_{i}:\gfra(F) = \gfra^{i}(F) \oplus \mathfrak{n}^{i}(F) \to \gfra^{i}(F)$ is the projection and
\[
	\hat{\Phibf}_{i}(1+y) = \Phibf_{i}(1+\pi_{i}(y)) = \psi \circ \Trd_{A/F}(c_{i}\pi_{i}(y))
\]
for $1+y \in K^{d} \cap G(F)_{x, \mathbf{s}_{i}+} = H^{t_{i}+1}(\beta, \Afra)$.  
However, we also can show $\Trd_{A/F}(c_{i}\pi_{i}(y)) = \Trd_{A/F}(c_{i}y)$ as in the proof of Proposition \ref{phiistheta1}.  
In conclusion, for $1+y \in H^{t_{i}+1}(\beta, \Afra)$ we obtain $\hat{\Phibf}_{i}(1+y) = \psi \circ \Trd_{A/F}(c_{i}y) = \psi_{c_{i}}(y)$ and (2) holds.  

Let $g \in B_{i}^{\times} \cap \Kfra(\Afra)$.  
We check that $g$ normalizes $H^{t_{i}+1}(\beta, \Afra)$.  
We consider two cases.  
First, suppose $i < d$.  
Then we have $H^{t_{i}+1}(\beta, \Afra) = G^{i+1}(F)_{x, \mathbf{s}_{i}+} \cdots G^{d}(F)_{x, \mathbf{s}_{d-1}+}$.  
Thus it suffices to show $g$ normalizes $G^{j}(F)_{x, \mathbf{s}_{j-1}+}$ for $i=j+1, \ldots, d$.  
However, we have
\[
	gG^{j}(F)_{x, \mathbf{s}_{j-1}+}g^{-1} = g \left( B_{j}^{\times} \cap \U^{t_{j-1}+1}(\Afra) \right) g^{-1} = ( gB_{j}^{\times} g^{-1} ) \cap (g \U^{t_{j-1}+1}(\Afra) g^{-1}).  
\]
Since $g \in B_{i}^{\times} \subset B_{j}^{\times}$ we have $gB_{j}^{\times} g^{-1} = B_{j}^{\times}$.  
Moreover, we also have $g \U^{t_{j-1}+1}(\Afra) g^{-1} = \U^{t_{j-1}+1}(\Afra)$ as $g \in \Kfra(\Afra)$.  
Therefore we obtain $gG^{j}(F)_{x, \mathbf{s}_{j}+}g^{-1} = B_{j}^{\times} \cap \U^{t_{j-1}+1}(\Afra) = G^{i}(F)_{x, \mathbf{s}_{j}+}$.  
Next, suppose $i = d = s$.  
Then we have $H^{t_{s}+1}(\beta, \Afra) = G^{d}(F)_{x, \mathbf{s}_{s}+} = \U^{t_{s}+1}(\Afra)$.  
Since $g \in \Kfra(\Afra)$, we obtain 
\[
	g H^{t_{s}+1}(\beta, \Afra) g^{-1} = g \U^{t_{s}+1}(\Afra) g^{-1} = \U^{t_{s}+1}(\Afra) = H^{t_{s}+1}(\beta, \Afra).  
\]
Therefore we obtain (3).  

Here, let $g$ be as above and $h \in H^{1}(\beta, \Afra)$.  
Since 
\[
H^{1}(\beta, \Afra) = \left(B_{i}^{\times} \cap H^{1}(\beta, \Afra) \right) H^{t_{i}+1}(\beta, \Afra), 
\]
we have $h=bh'$ for some $b \in B_{i}^{\times} \cap H^{1}(\beta, \Afra)$ and $h' \in H^{t_{i}+1}(\beta, \Afra)$.  
By the above argument, we have $h' \in H^{t_{i}+1}(\beta, \Afra) = g H^{t_{i}+1}(\beta, \Afra) g^{-1}$ and $h'$ is an element in $H^{1}(\beta, \Afra) \cap gH^{1}(\beta, \Afra) g^{-1}$.  
Then, $h \in H^{1}(\beta, \Afra) \cap g H^{1}(\beta, \Afra) g^{-1}$ if and only if $b \in H^{1}(\beta, \Afra) \cap g H^{1}(\beta, \Afra) g^{-1}$.  
Suppose $h \in H^{1}(\beta, \Afra) \cap g H^{1}(\beta, \Afra) g^{-1}$.  
Therefore we obtain
\begin{eqnarray*}
	\hat{\Phibf}_{i}(g^{-1}hg) & = & \hat{\Phibf}_{i} \left( (g^{-1}bg)(g^{-1}h'g) \right) 	\\
	& = & \hat{\Phibf}_{i}(g^{-1}bg) \hat{\Phibf}_{i}(g^{-1}h'g) = \Phibf_{i}(g^{-1}bg) \psi_{c_{i}}(g^{-1}h'g).  
\end{eqnarray*}
Here, since $\Phibf_{i}$ is a character of $G^{i}(F)=B_{i}^{\times}$ and $g \in B_{i}^{\times}$, we have $\Phibf_{i}(g^{-1}bg) = \Phibf_{i}(b)$.  
Moreover, since $c_{i}$ is an element in $E_{i}$, which is the center of $B_{i}$, we also have
\begin{eqnarray*}
	\psi_{c_{i}}(g^{-1}h'g) & = & \psi \circ \Trd_{A/F}(c_{i}g^{-1}h'g) = \psi \circ \Trd_{A/F}(g^{-1}c_{i}h'g) \\
	& = & \psi \circ \Trd_{A/F}(c_{i}h') = \psi_{c_{i}}(h').  
\end{eqnarray*}
Therefore we obtain $\hat{\Phibf}_{i}(g^{-1}hg) = \Phibf_{i}(b)\psi_{c_{i}}(h') = \Phibf_{i}(bh') = \Phibf_{i}(h)$, which implies (4).  
\end{prf}

\begin{prop}
We have $\prod_{i=0}^{d} \hat{\Phibf}_{i} \in \mathscr{C}(\beta, 0, \Afra)$.  
\end{prop}

\begin{prf}
If $s=-\infty$, then $\Phibf_{d}=1$ and $\hat{\Phibf}_{d}=1$, and then $\prod_{i=0}^{d} \hat{\Phibf}_{i} = 1 \in \mathscr{C}(\beta, 0, \Afra)$.  
Therefore we assume $s \in \Z$.  
If $d=s+1$, then $\Phibf_{d}=1$ and $\hat{\Phibf}_{i} = 1$ and we have $\prod_{j=i}^{d} \hat{\Phibf}_{j} = \prod_{j=i}^{s} \hat{\Phibf}_{j}$ for $i=0, \ldots, s$.  
Thus we show $\bar{\theta}_{i} := \prod_{j=i}^{s} \hat{\Phibf}_{j}|_{H^{t_{j}+1}(\beta, \Afra)} \in \mathscr{C}(\beta_{i}, \lfloor r_{i}/2 \rfloor, \Afra)$ by downward induction on $i=0, \ldots, s$.  

First, suppose $i=s$.  
Since $\beta_{s} = c_{s}$ is minimal over $F$, we need to check (1), (2) and (3) in Definition \ref{defofsimpch}.  
(2) is already shown as Lemma \ref{lemforbeingsimp} (1).  
Since $-v_{\Afra}(c_{s}) = -v_{\Afra}(\beta_{s}) = n$, we have $H^{t_{s}+1}(\beta, \Afra) = \U^{\lfloor n/2 \rfloor +1}(\Afra)$ and (3) is also shown as Lemma \ref{lemforbeingsimp} (2).  
Let $g \in B_{i}^{\times} \cap \Kfra(\Afra)$ and $h \in H^{t_{s}+1}(\beta, \Afra)$.  
Then $g^{-1}hg \in H^{t_{s}+1}(\beta, \Afra)$ by Lemma \ref{lemforbeingsimp} (3), and $\hat{\Phibf}_{i}(g^{-1}hg) = \hat{\Phibf}_{i}(h)$ by Lemma \ref{lemforbeingsimp} (4), which implies (1).  
Therefore $\hat{\Phibf}_{s} \in \mathscr{C}(\beta_{s}, t_{s}, \Afra)$.  

Next, suppose $0<i<s$.  
Since $k_{0}(\beta_{i-1}, \Afra) = v_{\Afra}(c_{i-1}) = -r_{i} > -n = v_{\Afra}(\beta_{i-1})$, the element $\beta_{i-1}$ is not minimal over $F$, and then we need to check (1), (2) and (4) in Definition \ref{defofsimpch}.  

To show (1), let $g \in B_{i-1}^{\times} \cap \Kfra(\Afra)$ and $h \in H^{t_{i-1}+1}(\beta, \Afra)$.  
Then $g^{-1}hg \in H^{t_{i-1}+1}(\beta, \Afra)$ by Lemma \ref{lemforbeingsimp} (3).  
For $j=i-1, \ldots, s$, we have $g \in B_{i-1}^{\times} \cap \Kfra(\Afra) \subset B_{j}^{\times} \cap \Kfra(\Afra)$.  
Therefore by Lemma \ref{lemforbeingsimp} (4) we have $\hat{\Phibf}_{j}(g^{-1}hg) = \hat{\Phibf}_{j}(h)$ and $\bar{\theta}_{i-1}(g^{-1}hg) = \prod_{j=i-1}^{s}\hat{\Phibf}_{j}(g^{-1}hg) = \prod_{j=i-1}^{s} \hat{\Phibf}_{j}(h) = \bar{\theta}_{i-1}(h)$, whence (1) holds.  

For $j=i-1, \ldots, s$, the restriction of $\hat{\Phibf}_{j}$ to $B_{j}^{\times} \cap H^{t_{i-1}+1}(\beta, \Afra)$ factors through $\Nrd_{B_{j}/E_{j}}$.  
Since $\Nrd_{B_{j}/E_{j}} |_{B_{i-1}^{\times}} = \Nrm_{E_{i-1}/E_{j}} \circ \Nrd_{B_{i-1}/E_{i-1}}$, the restriction of $\hat{\Phibf}_{j}$ to $B_{i-1}^{\times} \cap H^{t_{i-1}+1}(\beta, \Afra)$ factors through $\Nrd_{B_{i-1}/E_{i-1}}$.  
Then the character $\bar{\theta}_{i-1}=\prod_{j=i-1}^{s} \hat{\Phibf}_{j}|_{B_{i-1}^{\times} \cap H^{t_{i-1}+1}(\beta, \Afra)}$ also factors through $\Nrd_{B_{i-1}/E_{i-1}}$ and (2) holds.  

We show (4).  
We put $r'_{i-1}=0$ and $r'_{j}=r_{i}$ for $j=i, \ldots, s$.  
Then the sequence $([\Afra, n, r'_{(i-1)+i'}, \beta_{(i-1)+i'}])_{i'=0}^{s-i+1}$ is a defining sequence of $[\Afra, n, 0, \beta_{i-1}]$.  
Since $k_{0}(\beta_{i-1}, \Afra) = r_{i}$, we have $\max \{ \lfloor r_{i-1}/2 \rfloor , \lfloor k_{0}(\beta_{i-1}, \Afra)/2 \rfloor \} = \lfloor r_{i}/2 \rfloor = t_{i-1}$.  
Then $\bar{\theta}_{i-1}|_{H^{t_{i-1}+1}(\beta, \Afra)} = \bar{\theta}_{i} \hat{\Phibf}_{i-1}|_{H^{t_{i-1}+1}(\beta, \Afra)}$.  
The character $\bar{\theta}_{i}$ is an element in $\mathscr{C}(\beta_{i}, \lfloor r_{i}/2 \rfloor, \Afra)$ by induction hypothesis.  
On the other hand, $\hat{\Phibf}_{i-1}|_{H^{t_{i-1}+1}(\beta, \Afra)} = \psi_{c_{i-1}}$ by Lemma \ref{lemforbeingsimp} (2).  
Therefore (4) is shown and we complete the proof.  
\end{prf}

We put $\theta = \prod_{i=0}^{d} \hat{\Phi}_{i}$, and let $\eta_{\theta}$ be the Heisenberg representation of $\theta$.  

\begin{prop}
\label{betaext'}
$\kappa_{0} \otimes \cdots \otimes \kappa_{d}$ is an extension of $\eta_{\theta}$ to $K^{d}$.  
\end{prop}

\begin{prf}
Similar to the proof of Proposition \ref{betaext}.  
\end{prf}

\begin{thm}
\label{Main2}
Let $\left(x, (G^{i})_{i=0}^{d}, (\rbf_{i})_{i=0}^{d}, (\Phibf_{i})_{i=0}^{d}, \rho \right)$ be a Yu datum.  
Then there exists a maximal, tame simple type $(J, \lambda)$ associated with $[\Afra, n, 0, \beta]$ and a maximal extension $(\tilde{J}, \Lambda)$ of $(J, \lambda)$ such that
\begin{enumerate}
\item $\hat{J} :=\hat{J}(\beta, \Afra) = K^{d}$, and
\item $\rho_{d} = \cInd_{\tilde{J}}^{\hat{J}} \Lambda$.  
\end{enumerate}
\end{thm}

\begin{prf}
We can construct a tame simple stratum $[\Afra, n, 0, \beta]$ and a simple character $\theta \in \mathscr{C}(\beta, \Afra)$ as above.  
We take a $\beta$-extension $\kappa$ of $\eta_{\theta}$ and an extension $\hat{\kappa}$ of $\kappa$ to $\hat{J}$ by Lemma \ref{exist_ext_of_b-ext} (1).  
On the other hand, let $\kappa_{i}$ be the representation of $K^{d}$ as in Section 3 for $i=-1, 0, \ldots, d$.  
By Proposition \ref{betaext'}, the representation $\hat{\kappa}'=\kappa_{0} \otimes \cdots \otimes \kappa_{d}$ is an extension of the $\beta$-extension ${}^{\circ} \lambda$ of $\eta_\theta$ to $\kappa_{d}$.  
Then by Lemma \ref{exist_ext_of_b-ext} (2), there exists a character $\chi$ of $\hat{J}/J^{1}(\beta, \Afra)$ such that $\hat{\kappa}' \cong \hat{\kappa} \otimes \chi$.  
The representation $\kappa_{-1}$ is the extension of $\rho$ to $K^{d}$, trivial on $K^{d} \cap G(F)_{x, 0+} = J^{1}(\beta, \Afra)$.  

We construct ``depth-zero part" $\sigma$ of a simple type from $\rho$.  
By Lemma \ref{depth0ofYu}, there exists a depth--zero simple type $(G^{0}(F)_{x}, \sigma^{0})$ of $G^{0}(F)$ and a maximal extension $(\tilde{J}^{0}, \tilde{\sigma}^{0})$ such that $\rho \cong \Ind_{\tilde{J}^{0}}^{G^{0}(F)_{[x]}} \tilde{\sigma}^{0}$.  
We put $\tilde{J} = \tilde{J}^{0}J=\tilde{J}^{0}J^{1}(\beta, \Afra)$.  
Since $J^{1}(\beta, \Afra) \cap G^{0}(F) = G^{0}(F)_{x, 0+}$, we can extend $\tilde{\sigma}^{0}$ to $\tilde{J}$ as $\tilde{\sigma}^{0}$ is trivial on $J^{1}(\beta, \Afra)$.  
Let $\tilde{\sigma}$ be the $\tilde{J}$-representation $\tilde{\sigma}^{0}$, and we put $\sigma = \Res_{J}^{\tilde{J}} \tilde{\sigma}$.  
The representation $\sigma$ is an extension of $\sigma^{0}$ to $J$, trivial on $J^{1}(\beta, \Afra)$.  
Since $(G^{0}(F)_{x}, \sigma^{0})$ is a maximal simple type of depth zero and $\chi$ is a character of $\hat{J}$ trivial on $J^{1}(\beta, \Afra)$, the $J(\beta, \Afra)/J^{1}(\beta, \Afra)$-representation $\sigma \otimes \chi$ is cuspidal, and then $(J, \sigma \otimes \chi \otimes \kappa)$ is a simple type.  
By construction of $\tilde{J}$ and $\tilde{\sigma}$, the pair $(\tilde{J}, \tilde{\sigma} \otimes \Res_{\tilde{J}}^{\hat{J}}(\chi \otimes \hat{\kappa}))$ is a maximal extension of $(J, \sigma \otimes \chi \otimes \kappa)$.  
We put $\Lambda = \tilde{\sigma} \otimes \Res_{\tilde{J}}^{\hat{J}}( \chi \otimes \hat{\kappa})$.  

The representation $\kappa_{-1}$ is the extension of $\rho$ as $\kappa_{-1}$ is trivial on $K_{+}^{0}J^{1} \cdots J^{d} = J^{1}(\beta, \Afra)$, that is, the representation $\kappa_{-1}$ is $\rho$ regarded as a representation of $K^{d}=\hat{J}$ via $K^{0}/K_{+}^{0} = \U(\Bfra)/\U^{1}(\Bfra) \cong K^{d}/(K_{+}^{0}J^{1} \cdots J^{d}) = \hat{J}/J^{1}(\beta, \Afra)$.  
Then we have $\kappa_{-1} \cong \cInd_{\tilde{J}}^{\hat{J}} \tilde{\sigma}$ by Lemma \ref{inf_and_ind} and
\[
	\cInd_{\tilde{J}}^{\hat{J}} \Lambda \cong (\cInd_{\tilde{J}}^{\hat{J}} \tilde{\sigma}) \otimes \chi \otimes \hat{\kappa} \cong \kappa_{-1} \otimes \kappa_{0} \otimes \cdots \otimes \kappa_{d} = \rho_{d}, 
\]
which finishes the proof.  
\end{prf}

\begin{cor}
\label{tame_and_esstame}
The set of essentially tame supercuspidal representations of $G$ is equal to the set of tame supercuspidal representations of $G$.  
\end{cor}

\begin{prf}
Let $\pi$ be an irreducible supercuspidal representation of $G$.  
Since $\cInd_{K^{d}(\Psi)}^{G} \rho^{d}(\Psi)$ is irreducible for any Yu's datum $\Psi$, $\pi$ is tame supercuspidal if and only if $\pi \supset \rho^{d}(\Psi)$ for some $\Psi$.  
However, by Theorem \ref{Main1} and \ref{Main2} it holds if and only if $\pi$ contains some compact induction of a maximal extension $(\tilde{J}, \Lambda)$ of a tame, maximal simple type, which implies $\pi$ is essentially tame by \ref{esstame}.  
\end{prf}

\bigbreak\bigbreak
\noindent Yuki Yamamoto\par
\noindent Graduate School of Mathematical Sciences, 
The University of Tokyo, 3--8--1 Komaba, Meguro-ku, 
Tokyo, 153--8914, Japan\par
\noindent E-mail address: \texttt{yukiymmt@ms.u-tokyo.ac.jp}

\end{document}